\UseRawInputEncoding

\RequirePackage[l2tabu, orthodox]{nag}
\RequirePackage{fixltx2e}
\documentclass[10pt,a4paper,reqno,oneside,final, english]{smfart}

\usepackage{leftidx}
\usepackage[titletoc]{appendix}
\usepackage{tikz}\usetikzlibrary{decorations.markings,matrix,arrows}
\usepackage{tikz-cd}
\usepackage{amsfonts}
\usepackage{amsthm}
\usepackage[T1]{fontenc}
\usepackage[mathscr]{eucal}
\usepackage{setspace}
\usepackage{amssymb,latexsym,amsmath,amscd}
\usepackage{graphicx}
\usepackage{showkeys}
\usepackage{layout}
\usepackage{enumerate}
\usepackage{enumitem}
\usepackage{indentfirst}
\usepackage[varg]{pxfonts}
\usepackage[stretch=10]{microtype}
\usepackage{booktabs}
\usepackage{xspace}
\usepackage{calrsfs}
\usepackage{filecontents}
\usepackage{mathtools}
\usepackage{titlesec}
\usepackage{fixltx2e}
\usepackage{todonotes}
\usepackage{scalerel,stackengine}
\usepackage{enumitem}
\usepackage[backref = page, colorlinks   = true,
     citecolor    = gray, pdfborder={0 0 0}]{hyperref}
\usepackage{nameref}
\usepackage{cleveref}
\usepackage{braket}

\allowdisplaybreaks

\titleformat{\paragraph}
{\normalfont\normalsize\bfseries}{\theparagraph}{1em}{}
\titlespacing*{\paragraph}
{0pt}{3.25ex plus 1ex minus .2ex}{1.5ex plus .2ex}

\makeatletter
\newcommand{\neutralize}[1]{\expandafter\let\csname c@#1\endcsname\count@}
\makeatother

\theoremstyle{plain}
\newtheorem{thmm}{Theorem}[section]

\newtheorem{thm}[thmm]{Theorem}
\newtheorem*{thm*}{Theorem}

\newtheorem{lem}[thmm]{Lemma}

\newtheorem{lem-def}[thmm]{Lemma-Definition}
\newtheorem{claim}[thmm]{Claim}
\newtheorem{pro}[thmm]{Proposition}

\newtheorem{pro-def}[thmm]{Proposition-Definition}
\newtheorem{cor}[thmm]{Corollary}

\newtheorem{conj}[thmm]{Conjecture}
\newtheorem{que}[thmm]{Question}

\newtheorem{assump}[thmm]{Assumption}

\theoremstyle{definition}
\newtheorem{Def}[thmm]{Definition}

\newtheorem{rem}[thmm]{Remark}

\theoremstyle{remark}
\newtheorem{ex}[thmm]{Example}

\newcommand{\ssec}{\subsection}

\newcommand{\sssec}{\subsubsection}

\newcommand{\ol}{\overline}
\newcommand{\ti}[1]{\tilde{#1}}
\newcommand{\ul}{\underline}
\newcommand{\vast}{\bBigg@{4}}
\newcommand{\Vast}{\bBigg@{5}}
\newcommand{\wt}{\widetilde}

\newcommand{\te}{\ti{e}}

\stackMath
\newcommand\reallywidehat[1]{%
\savestack{\tmpbox}{\stretchto{%
  \scaleto{%
    \scalerel*[\widthof{\ensuremath{#1}}]{\kern-.6pt\bigwedge\kern-.6pt}%
    {\rule[-\textheight/2]{1ex}{\textheight}}
  }{\textheight}%
}{0.5ex}}%
\stackon[1pt]{#1}{\tmpbox}%
}
\definecolor{armygreen}{rgb}{0.29, 0.33, 0.13}
\definecolor{ao(english)}{rgb}{0.5, 0.2, 0.0}

\newcommand{\bC}{\mathbf{C}}

\newcommand{\bH}{\mathbf{H}}

\newcommand{\bP}{\mathbf{P}}
\newcommand{\bQ}{\mathbf{Q}}
\newcommand{\bR}{\mathbf{R}}

\newcommand{\bZ}{\mathbf{Z}}

\newcommand{\cA}{\mathcal{A}}
\newcommand{\cAut}{\mathcal{A}ut}
\newcommand{\cB}{\mathcal{B}}
\newcommand{\cC}{\mathcal{C}}
\newcommand{\cD}{\mathcal{D}}
\newcommand{\cE}{\mathcal{E}}
\newcommand{\cF}{\mathcal{F}}

\newcommand{\cI}{\mathcal{I}}
\newcommand{\cJ}{\mathcal{J}}
\newcommand{\cK}{\mathcal{K}}
\newcommand{\cL}{\mathcal{L}}
\newcommand{\cM}{\mathcal{M}}

\newcommand{\cO}{\mathcal{O}}
\newcommand{\cP}{\mathcal{P}}

\newcommand{\cS}{\mathcal{S}}
\newcommand{\cT}{\mathcal{T}}
\newcommand{\cU}{\mathcal{U}}
\newcommand{\cV}{\mathcal{V}}

\newcommand{\cX}{\mathcal{X}}
\newcommand{\cY}{\mathcal{Y}}
\newcommand{\cZ}{\mathcal{Z}}
\newcommand{\sA}{\mathscr{A}}

\newcommand{\sC}{\mathscr{C}}

\newcommand{\sE}{\mathscr{E}}
\newcommand{\sL}{\mathscr{L}}

\newcommand{\sS}{\mathscr{S}}

\newcommand{\cExt}{\mathcal{E}xt}
\newcommand{\cHom}{\mathcal{H}om}

\newcommand{\fm}{\mathfrak{m}}

\newcommand{\fT}{\mathfrak{T}}
\newcommand{\fU}{\mathfrak{U}}
\newcommand{\fV}{\mathfrak{V}}

\newcommand{\gD}{\Delta}

\newcommand{\gO}{\Omega}
\newcommand{\gS}{\Sigma}
\newcommand{\gT}{\Theta}

\newcommand{\ga}{\alpha}
\newcommand{\gb}{\beta}
\newcommand{\gd}{\delta}
\newcommand{\gep}{\varepsilon}
\newcommand{\ep}{\varepsilon}
\newcommand{\gk}{\kappa}
\newcommand{\gl}{\lambda}
\newcommand{\go}{\omega}
\newcommand{\gs}{\sigma}
\newcommand{\gt}{\theta}

\newcommand{\Alb}{\mathrm{Alb}}

\newcommand{\At}{\mathrm{At}}
\newcommand{\Aut}{\mathrm{Aut}}

\newcommand{\can}{\mathrm{can}}

\newcommand{\Coh}{\mathrm{Coh}}
\newcommand{\coker}{\mathrm{coker}}
\newcommand{\Cone}{\mathrm{Cone}}

\newcommand{\Defo}{\mathrm{Def}}
\newcommand{\Div}{\mathrm{Div}}

\newcommand{\Ext}{\mathrm{Ext}}

\newcommand{\Gal}{\mathrm{Gal}}

\newcommand{\Hom}{\mathrm{Hom}}

\newcommand{\Id}{\mathrm{Id}}
\newcommand{\Ima}{\mathrm{Im}}

\newcommand{\KS}{\mathrm{KS}}

\newcommand{\mmin}{\mathrm{min}}

\newcommand{\NS}{\mathrm{NS}}

\newcommand{\ob}{\mathrm{ob}}

\newcommand{\Pic}{\mathrm{Pic}}

\newcommand{\pr}{\mathrm{pr}}

\newcommand{\rank}{\mathrm{rank}}

\newcommand{\rk}{\mathrm{rk}}

\newcommand{\Spec}{\mathrm{Spec \ }}
\newcommand{\Sing}{\mathrm{Sing}}

\newcommand{\torsion}{\mathrm{torsion}}
\newcommand{\Tot}{\mathrm{Tot}}
\newcommand{\tr}{\mathrm{tr}}
\newcommand{\Tr}{\mathrm{Tr}}

\newcommand{\eg}{\emph{e.g.} }

\newcommand{\ie}{\emph{i.e.} }

\newcommand{\bss}{\setminus}

\newcommand{\colonec}{\mathrel{:=}}
\newcommand{\cnec}{\mathrel{:=}}
\newcommand{\ctr}{\lrcorner}
\newcommand{\cupp}{\mathbin{\smile}}

\newcommand{\dr}{\partial}

\newcommand{\vv}{{\vee\vee}}

\renewcommand{\(}{\left(}
\renewcommand{\)}{\right)}

\newcommand{\cto}{\circlearrowleft}
\newcommand{\dto}{\dashrightarrow}
\newcommand{\lto}{\leftarrow}
\newcommand{\tto}{\twoheadrightarrow}
\newcommand*\eto{%
  \xrightarrow[]{\raisebox{-0.25 em}{\smash{\ensuremath{\sim}}}}%
}
\newcommand{\hto}{\hookrightarrow}
\newcommand{\xto}[1]{\xrightarrow{ #1 }}
\makeatletter
\newcommand{\xdto}[2][]{\ext@arrow 0359\rightarrowfill@@{#1}{#2}}
\def\rightarrowfill@@{\arrowfill@@\relax\relbar\rightarrow}
\def\leftarrowfill@@{\arrowfill@@\leftarrow\relbar\relax}
\def\leftrightarrowfill@@{\arrowfill@@\leftarrow\relbar\rightarrow}
\def\arrowfill@@#1#2#3#4{%
	$\m@th\thickmuskip0mu\medmuskip\thickmuskip\thinmuskip\thickmuskip
	\relax#4#1
	\xleaders\hbox{$#4#2$}\hfill
	#3$%
}
\makeatother
\newcommand{\xhto}[1]{\xhookrightarrow{ #1 }}

\makeatletter
\let\orgdescriptionlabel\descriptionlabel
\renewcommand*{\descriptionlabel}[1]{%
  \let\orglabel\label
  \let\label\@gobble
  \phantomsection
  \edef\@currentlabel{#1}%
  \let\label\orglabel
  \orgdescriptionlabel{#1}%
}
\makeatother
\tikzset{node distance=2cm, auto}

\numberwithin{equation}{section}

\title{Algebraic approximations of compact K\"ahler threefolds}

\author{Hsueh-Yung Lin}
\address{Department of Mathematics, National Taiwan University, 
	No. 1, Sec. 4, Roosevelt Rd., Taipei 10617, Taiwan.}
\email{hsuehyunglin@ntu.edu.tw}

\begin{document}

\begin{abstract}
We prove that every compact K\"ahler manifold of dimension 3 has arbitrarily 
small deformations to some projective manifolds, 
thereby solving the Kodaira problem in dimension 3.
As a consequence,
we show that more generally, every 
compact complex threefold in the Fujiki class $\cC$ with at worst rational singularities
admits algebraic approximations.
This implies Peternell's conjecture on the existence of algebraic approximations of minimal K\"ahler varieties in dimension $3$.
\end{abstract}

\maketitle

\section{Introduction}

\ssec{The Kodaira problem in dimension 3}
\hfill

Let $X$ be a compact K\"ahler manifold. An \emph{algebraic approximation of $X$} is a deformation $\cX \to \gD$ of $X$ such that up to shrinking $\gD$, the subset parameterizing projective manifolds in this family is dense in $\gD$. The so-called Kodaira problem asks whether a compact K\"ahler manifold in general admits an algebraic approximation. 
While Kodaira proved that compact K\"ahler surfaces always have algebraic approximations in the early 1960s~\cite[Theorem 16.1]{KodairaSurfaceII}, starting from dimension 4 there exist compact K\"ahler manifolds in each dimension for which the Kodaira problem has a negative answer. Such examples were first constructed by Voisin~\cite{Voisincs}, 
and her examples even show that there exist
compact K\"ahler manifolds which do not have the homotopy type of a projective manifold.
See also~\cite[Proposition 5.3]{OguisoJap} for other examples constructed by Oguiso,
as well as~\cite{VoisinBiratKod, voisin2021compact} for further results.
 
The aim of this article is to answer in the affirmative the Kodaira problem in dimension 3, which has been an open problem since the work of Kodaira and Voisin mentioned above. 
Together with their work,
this completes the study of the Kodaira problem in terms of the dimension of the manifolds.

\begin{thm}\label{thm-main3}
Every compact K\"ahler manifold of dimension 3 has an algebraic approximation.
\end{thm}

So far, non-trivial examples of compact K\"ahler manifolds admitting algebraic approximations can be found in~\cite{KodairaSurfaceII, Buch2, Schrackdefo, CaoJApproxalg, GrafDefKod0, ClaudonToridefequiv, HYLbimkod1, ClaudonHorpi1, HYLkodfibellip, claudon2020kahler} and the list is rather exhaustive at present. 
Obviously, a positive answer to the Kodaira problem for a compact K\"ahler manifold $X$ implies that invariants of $X$ that are preserved under small deformations (such as the fundamental group and the Hodge diamond) can be realized by projective manifolds. Immediate corollaries of Theorem~\ref{thm-main3} of this sort are left to the readers.

The existence of algebraic approximations holds more generally (and naturally) for threefolds in the Fujiki class $\cC$ with rational singularities,
which we formulate now as the main theorem of the article.
By definition, a compact complex variety\footnote{We refer to \S\ref{ssec-defnot} for the convention of complex varieties we work with in this article.} $X$ is in the Fujiki class $\cC$ if it is meromorphically dominated by a compact K\"ahler manifold. Thanks to~\cite[Lemma 4.6.1]{FujikiClosednessDouady} and~\cite[Théorème 3]{VarouchasFPV}, a compact complex variety $X$ is in the Fujiki class $\cC$ if and only if there exist a compact K\"ahler manifold $\ti{X}$ and a bimeromorphic morphism $\ti{X} \to X$. 

Algebraic approximations of compact complex varieties are defined as follows.

\begin{Def}[Algebraic approximation]\label{def-appalg}
Let $X$ be a compact complex variety. An algebraic approximation of $X$ is a (flat) deformation $\Pi : \cX \to \gD$ of $X$ such that up to shrinking $\gD$, the subset of points of $\gD$ parameterizing Moishezon varieties is dense for the Euclidean topology. 
\end{Def}

The main result that we prove in this article is the following. 

\begin{thm}\label{thm-mainC}
	Every compact complex threefold in the Fujiki class $\cC$ with at worst rational singularities has an algebraic approximation.
\end{thm}

Applying Moishezon's criterion for varieties with rational singularities 
due to Namikawa~\cite[Theorem 1.6]{NAMIKAWA2002}, 
Theorem~\ref{thm-mainC} has the following immediate corollary, 
which contains in particular Theorem~\ref{thm-main3}.

\begin{cor}\label{cor-sing}
	Every compact K\"ahler threefold $X$ with at worst rational singularities 
	has a deformation $\cX \to \gD$ such that 
	the subset of points of $\gD$ parameterizing projective varieties is dense for the Euclidean topology.
\end{cor}

We have mentioned that 
there exist compact K\"ahler manifolds
which do not admit any 
algebraic approximation. 
So far all known examples are constructed as the blow-ups of some manifolds,
which are therefore not minimal.
For minimal varieties, Peternell formulated the following conjecture (see \eg~\cite[Conjecture 1.2]{GrafDefKod0}),
which is still wide open.

\begin{conj}[Peternell]\label{conj-Peternell}
	Every minimal K\"ahler variety has an algebraic approximation.
\end{conj}

In dimension 2, since minimal surfaces are smooth,
Conjecture~\ref{conj-Peternell} is covered by Kodaira's theorem~\cite[Theorem 16.1]{KodairaSurfaceII}.
In general, minimal varieties are singular, but they only have terminal (and therefore rational) singularities. 
Thus in dimension 3, 
Conjecture~\ref{conj-Peternell} is covered by
Corollary~\ref{cor-sing}.

\begin{cor}
	Conjecture~\ref{conj-Peternell} holds in dimension 3.
\end{cor}

Let $X$ be a compact complex threefold as in Theorem~\ref{thm-mainC}.
Our approach of constructing algebraic approximations of $X$
consists of two steps. First we find a simpler bimeromorphic model $X \dto X'$ of $X$,
then we construct an algebraic approximation 
$\cX' \to \gD$ of $X'$ with some \emph{locally trivial} properties,
allowing it to induce an algebraic approximation of $X$ 
(see \S\ref{ssec-appgenintro} for a presentation of the strategy).
In view of this approach, as being K\"ahler 
is not invariant under bimeromorphic modifications,
this is the reason why Theorem~\ref{thm-mainC} is considered more natural than
Theorem~\ref{thm-main3}.

In the remainder of the introduction, 
we provide an overview of the proof of Theorem~\ref{thm-mainC}
and formulate the
intermediate results leading to it.

\ssec{A general approach proving the existence of algebraic approximations}\label{ssec-appgenintro}
\hfill

Given a compact complex variety $X$, one way to prove that $X$ has an algebraic approximation is to find a (simpler) bimeromorphic model $\nu : X \dto X'$ of $X$ together with an algebraic approximation $\cX' \to \gD$ of $X'$ which induces a deformation of $\nu$. Note that in general, it is not enough if we only prove that $X'$ has an algebraic approximation, because asking a deformation $\Pi : \cX' \to \gD$  of $X'$ to induce deformations of $X$ usually imposes non-trivial restrictions on $\Pi$. For instance, if $\nu : X \to X'$ is the blow-up of $X'$ along a submanifold $Y \subset X'$, then a necessary condition for $\Pi : \cX' \to \gD$ to induce a deformation of $X$ is that the submanifold $Y$ is preserved along the deformation $\Pi$~\cite{KodStab}.

On the other hand, we observe that if $Y \subset X'$ is a subvariety such that $\nu^{-1}_{|X' \bss Y}$ is biholomorphic onto its image, then to prove that $\Pi : \cX' \to \gD$ is a deformation of $X'$ which induces deformations of $\nu : X \dto X'$, it suffices to show that there exists a neighborhood $U \subset X'$ of $Y$ which deforms trivially along $\gD$. 
If such a neighborhood exists, we will call $\cX' \to \gD$ a \emph{$Y$-locally trivial} deformation (see Definition~\ref{def-loctriv}). 
This leads to the following proposition, which we will prove in  \S\ref{ssec-appgen}.

\begin{pro}\label{pro-red}
Let $X'$ be a normal compact complex variety and let $X$ be a 
compact complex variety with at worst rational singularities 
which is bimeromorphic to $X'$. 
Assume that $X'$ has a $Y$-locally trivial algebraic approximation for every subvariety $Y \subset X'$ satisfying $\dim Y \le \dim X' - 2$. 
Then $X$ has an algebraic approximation.

\end{pro}

Thus finding an algebraic approximation of $X$ is transformed into the problem of finding: 
\begin{itemize}
	\item[$\bullet$] a normal bimeromorphic model $X'$ of $X$;
	\item[$\bullet$] an algebraic approximation of $X'$ 
	which is $Y$-locally trivial for every $Y \subset X'$ with $\dim Y \le \dim X - 2$.
\end{itemize}

This is how we will prove Theorem~\ref{thm-mainC} for uniruled threefolds, 
as well as threefolds of algebraic dimension $a \le 1$.
We formulate the results as follows, which are the two main propositions
we prove in this article.

\begin{pro}\label{pro-main3gk1}
Let $X$ be a compact complex threefold in the Fujiki class $\cC$. If $X$ is non-uniruled and $a(X) \le 1$, then $X$ is bimeromorphic to a normal compact complex variety $X'$
such that $X'$ has an algebraic approximation which is $Y$-locally trivial
for every subvariety $Y \subset X'$ with $\dim Y \le 1$.
\end{pro}

\begin{pro}\label{pro-main3unir}
	Let $X$ be a compact complex threefold in the Fujiki class $\cC$. 
	If $X$ is uniruled, then $X$ is bimeromorphic to a normal compact complex variety $X'$
	such that $X'$ has an algebraic approximation which is $Y$-locally trivial
	for every subvariety $Y \subset X'$ with $\dim Y \le 1$.
\end{pro}

While the conclusions in both Propositions~\ref{pro-main3gk1} and~\ref{pro-main3unir} are identical,
the approaches proving them are different, depending on the geometry of the threefold $X$ (see \S\ref{ssec-introale1} and \S\ref{ssec-glunir}).
This is why, and also for easier reference, we separate the statement into two propositions.

As for threefolds $X$ of algebraic dimension $a(X) \ge 2$, if $X$ is non-algebraic then $a(X) = 2$. 
Theorem~\ref{thm-mainC} for these remaining threefolds is covered
by the following theorem we proved in~\cite{HYLkodfibellip}.

\begin{thm}[{\cite[Corollary 1.4]{HYLkodfibellip}}]\label{thm-HYLkodfibellip}
Let $X$ be a compact complex variety in the Fujiki class $\cC$ with at worst rational singularities. If $a(X) = \dim X - 1$, then $X$ has an algebraic approximation. 
\end{thm}

We will therefore focus on uniruled threefolds and non-uniruled threefolds of algebraic dimension $a \le 1$ in this article. 
Let us first describe the
bimeromorphic models of the threefolds in questions
that we choose to prove Propositions~\ref{pro-main3gk1} and~\ref{pro-main3unir}.

\ssec{Bimeromorphic models of compact K\"ahler threefolds}
\hfill

First we describe the bimeromorphic models that we choose to prove Proposition~\ref{pro-main3gk1}.

\begin{pro}\label{pro-classk01} 
Let $X_0$ be a compact complex threefold in the Fujiki class $\cC$.
If $X_0$ is non-uniruled and has algebraic dimension $a(X_0) \le 1$, 
then $X_0$ is bimeromorphic to one of the following:
\begin{enumerate}[label = \roman{enumi})]
\item $X = \ti{X}/G$ where $G$ is a finite group and 
$\ti{X}$ is either a 3-torus of algebraic dimension 0, 
or the product $S \times B$ of a K-trivial surface $S$ of algebraic dimension 0 with a smooth projective curve $B$. In the latter case, the $G$-action on $\ti{X} = S \times B$ is diagonal.
Here, 
a K-trivial surface is a smooth compact K\"ahler surface $S$ such that $K_S \simeq \cO_S$.
\item The total space of a fibration $f : X \to B$ over a smooth projective curve $B$
such that $f$ has no  multi-section and a general fiber $F$ of $f$ is a 2-torus. 
Moreover, we can choose $X$ to be a smooth compact K\"ahler threefold; 
or when $F$ is non-algebraic, 
we can choose $f$ to be the quotient by the finite group $G$ of a $G$-equivariant 
smooth isotrivial torus fibration 
$\ti{f} : \ti{X} \to \ti{B}$  over 
 a smooth projective curve $\ti{B}$.
\end{enumerate}
In either case, the bimeromorphic model $X$ is normal.
\end{pro}

\begin{rem}\label{rem-KMMP}
	Initially, the Minimal Model Program (MMP) for compact K\"ahler threefolds~\cite{HorPetsurvey} had served as an important guideline to the author to understand the bimeromorphic descriptions of non-algebraic compact K\"ahler threefolds. But at the end, it turns out that the only place where
	we use the MMP in this article is   
	to rule out the existence of simple non-Kummer varieties~\cite[Corollary 1.4]{das2023log}
	from Fujiki's classification~\cite{FujikiStruC} in the proof of Proposition~\ref{pro-classk01}.
\end{rem}

Before we move on to non-algebraic uniruled threefolds, 
let us first introduce some 
ruled surfaces which will appear in our bimeromorphic description of uniruled threefolds.

Let  $\sS \to \sC$ be a ruled surface. 
We say that $\sS$ is \emph{decomposable}
if $\sS = \bP(\cL \oplus \cL')$ for some line bundles $\cL$ and $\cL'$ over $\sC$.
Following~\cite[Chapter V.2]{Hart}, recall that
$\sS = \bP(\cE)$ for some normalized vector bundle $\sE$ over $\sC$
and we define $\gep(\sS) = - \deg \det(\sE) \in \bZ$ (which is denoted by $e(\sS)$ in~\cite[Chapter V.2]{Hart}).
The following trichotomy of ruled surfaces with 
$\gep = 0$ will be relevant in this work.

\begin{Def}\label{Def-typee=0}
	For ruled surfaces $\sS \to \sC$ with $\gep(\sS) = 0$,
	we define the following trichotomy:
	\begin{itemize}
		\item[I.] $\sS \simeq \bP^1 \times \sC$ over $\sC$; 
		\item[II.] $\sS$ is a decomposable ruled surface 
		such that $\sS \not\simeq \bP^1 \times \sC$ over $\sC$; 
		\item[III.] $\sS$ is indecomposable.
	\end{itemize}
\end{Def}

The bimeromorphic description of non-algebraic uniruled compact K\"ahler threefolds 
that we need is the following. 

\begin{pro}\label{pro-classuniregl}
	A non-algebraic uniruled compact complex threefold in the Fujiki class $\cC$ 
	is bimeromorphic to a threefold $X$ satisfying one of the following descriptions:
\begin{enumerate}[label = \roman{enumi})]
\item 
$X$ is the total space of a good $\bP^1$-bundle $f : X \to S$
(see Definition~\ref{Def-bonP1})
over a surface $S$ with $a(S) = 0$.
\item
$X$ is the total space of a $\bP^1$-fibration $f: X \to S$ 
over a non-algebraic elliptic surface $p : S \to B$.
Moreover, 
there exist a Zariski dense open subset $B^\circ \subset B$  and $N \in \{I,II,III\}$
such that for every $b \in B^\circ$, $ S_b\cnec p^{-1}(b)$ is a smooth elliptic curve
and  $X_b \cnec (p \circ f)^{-1}(b)$ 
is a ruled surface over $S_b$ with $\gep(X_b) = 0$ of type $N$ (see Definition~\ref{Def-typee=0}).
\end{enumerate}
In both cases, we can choose $X$ and $S$ to be compact K\"ahler manifolds.
\end{pro}

	The fibrations $X \to S$ in Proposition~\ref{pro-classuniregl} 
	are actually the MRC fibrations of $X$. 
	Descriptions \emph{i)} and \emph{ii)}
	correspond to the case where $a(S) = 0$ and $a(S) = 1$ respectively.

We will prove Propositions~\ref{pro-classk01} and~\ref{pro-classuniregl} in \S\ref{sec-bim}.
These statements follow essentially from
Fujiki's classification of algebraic reductions of compact K\"ahler threefolds~\cite{FujikiStruC}
together with some classical results on conic bundles~\cite{SarkisovConicstr} 
and the MMP for K\"ahler threefolds~\cite{HorPetsurvey,das2023log} 
(see also Remark~\ref{rem-KMMP}). 
They are minor improvements or variants of Fujiki's results adapted to our needs.

The bimeromorphic models in Propositions~\ref{pro-classk01} and~\ref{pro-classuniregl}
are the ones we choose to prove Propositions~\ref{pro-main3gk1} and~\ref{pro-main3unir} respectively.
Now we provide an overview of the proofs of Propositions~\ref{pro-main3gk1} and~\ref{pro-main3unir}
for these models, which are central in this article.
 
\ssec{Algebraic approximations of non-uniruled compact K\"ahler threefolds of algebraic dimension $a \le 1$}\label{ssec-introale1}
\hfill

We start with non-uniruled threefolds of algebraic dimension $a \le 1$. 
The following is a more precise statement of Proposition~\ref{pro-main3gk1} 
taking Proposition~\ref{pro-classk01} into account.

\begin{pro}\label{pro-mainpair}
	Let $X$ be one of the bimeromorphic models of
	non-uniruled threefolds in the Fujiki class $\cC$ of algebraic dimension $a \le 1$
	in Proposition~\ref{pro-classk01}. 
	Then $X$ has an algebraic approximation 
which is $Y$-locally trivial for every subvariety $Y \subset X$ with $\dim Y \le 1$.
\end{pro}

We will prove Proposition~\ref{pro-mainpair} in \S\ref{sec-casparcas}. 
Here we indicate the ingredients of the proof of Proposition~\ref{pro-mainpair}  
for each threefold appearing in Proposition~\ref{pro-classk01}.

Algebraic approximations of smooth torus fibrations are constructed 
by Claudon~\cite{ClaudonToridefequiv}
through the so-called \emph{tautological families} 
(which we will recall in \S\ref{ssec-fibTlissetaut}, together with some refined results).
When $X$ is a finite quotient of a 3-torus or a smooth isotrivial 2-torus fibration,
the proof of Proposition~\ref{pro-mainpair} for $X$ 
is based on a
more precise construction of tautological families; this will be carried out 
in \S\ref{ssec-Gtores} and \S\ref{ssec-2tore}.

If $X$ is the total space of a fibration whose general fiber is an abelian surface, 
then Proposition~\ref{pro-mainpair} is a consequence of the following theorem we proved in~\cite{HYLbimkod1}.

\begin{thm}[{\cite[Corollary 1.3]{HYLbimkod1}}]\label{thm-AbFibDefprec} 
Let  $f : X \to B$ be a fibration from a compact K\"ahler manifold $X$ to a smooth projective curve $B$. Suppose that a general fiber of $f$ is an abelian variety, then $f$ has an algebraic approximation
$$\Pi: \cX \to B \times V \to V  $$
 which is locally trivial over $B$.
\end{thm}

Finally, when $X$ is a finite quotient of the product of a K-trivial surface $S$ with a curve, 
Proposition~\ref{pro-mainpair} is based on 
the existence of algebraic approximations of $S$ with more precision (see \S\ref{ssec-Gsurf}).

\ssec{Algebraic approximations of uniruled threefolds}\label{ssec-glunir}
\hfill

Finally we outline the proof of Proposition~\ref{pro-main3unir} about uniruled threefolds.
This presents the main difficulties of the article,
as the existence of algebraic approximations is unknown for most uniruled compact K\"ahler threefolds, 
even up to bimeromorphic modifications.
We will develop new approaches  
constructing algebraic approximations for these threefolds in \S\ref{sec-a0K3} and \S\ref{sec-a1}.

We will prove the following result, which implies Proposition~\ref{pro-main3unir} by Proposition~\ref{pro-classuniregl}.

\begin{pro}\label{pro-sec}
Let $f: X \to S$ be one of the $\bP^1$-fibrations in Proposition~\ref{pro-classuniregl}
(which cover all uniruled threefolds in the Fujiki class $\cC$, up to bimeromorphic modifications). 
There exists a bimeromorphic modification $f' : X' \to S'$ of $f$ 
with $X'$ normal, together with a deformation 
$$\Pi' : \cX' \to \cS' \to \gD$$ 
of $f'$, such that the underlying deformation of $X'$ 
is an algebraic approximation 
which is $f'^{-1}(C)$-locally trivial for every proper subvariety $C \subsetneq S'$.  
\end{pro}

According to whether $f$ is in the first or the second case of Proposition~\ref{pro-classuniregl}, we will prove Proposition~\ref{pro-sec} in \S\ref{sec-a0K3} and \S\ref{sec-a1}, which corresponds to Proposition~\ref{pro-aaa=0} and  Corollary~\ref{cor-ellipsuraa} respectively. 
In the first case where $X\to S$ is a $\bP^1$-bundle, the main idea consists in regarding $X$ as the projectivization of a twisted vector bundle $E$ of rank 2, and studying the 
deformations of the pair $(S,E)$. 
The deformation theory of twisted vector bundles that we need parallels well 
the classical theory of the untwisted ones~\cite{BuchweitzFlennerACsemireg} (see \S\ref{ssec-deftrod}). 
From the outset, one could argue as if $E$ is a vector bundle in the usual sense. 

Finally in the second case where $f : X \to S$ is a $\bP^1$-fibration over an elliptic surface $p: S \to B$, recall that $p$ has an algebraic approximation $\Pi :  \cS \to B \times V \to V$ by Theorem~\ref{thm-AbFibDefprec}. 
Up to replacing $f : X \to S$ by some bimeromorphic modification of it,
we will show that the algebraic approximation $\Pi$ that we constructed to prove Theorem~\ref{thm-AbFibDefprec} can be 
lifted to an algebraic approximation of $f: X \to B$ 
which is locally trivial over $B$. 
Such a lifting
will be constructed by means of 1-cocycles of local relative automorphisms. 
This will imply Proposition~\ref{pro-sec} for $f : X \to S$.

 \ssec{Organization of the article and remark on the dependence on~\cite{HYLbimkod1, HYLkodfibellip}}
 
 \hfill

We will first introduce some basic terminologies 
and recall or prove various preliminary general results in \S\ref{sec-prelim}. 
 In \S\ref{sec-bim} we will prove Propositions~\ref{pro-classk01} and~\ref{pro-classuniregl},
 which provide bimeromorphic descriptions of 
 non-uniruled compact K\"ahler threefolds of algebraic dimension $a \le 1$, 
 and non-algebraic uniruled compact K\"ahler threefolds respectively.  
 Starting from \S\ref{sec-casparcas}, 
 we will study the existence of algebraic approximations for these bimeromorphic models 
 and prove Propositions~\ref{pro-main3gk1} and~\ref{pro-main3unir}. 
 More precisely, we will prove Proposition~\ref{pro-mainpair} 
 (which implies Proposition~\ref{pro-main3gk1}) in \S\ref{sec-casparcas} 
 and Proposition~\ref{pro-sec} (which implies Proposition~\ref{pro-main3unir}) 
 in \S\ref{sec-a0K3} and \S\ref{sec-a1}. 
 Assembling these results, 
 we will finish the proof of Theorem~\ref{thm-mainC} in \S\ref{sec-conclf}.
 
 The proof of Theorem~\ref{thm-mainC} depends on the main results of two other articles~\cite{HYLbimkod1, HYLkodfibellip} of the author (Theorems~\ref{thm-HYLkodfibellip} and~\ref{thm-AbFibDefprec}).  As we can see from the outline of the proof, the arguments proving Theorem~\ref{thm-mainC} vary between different threefolds and for certain threefolds, the existence of algebraic approximations turns out to be a more general phenomenon. This is the case for fibrations in abelian surfaces over a curve and threefolds of algebraic dimension 2, and the corresponding statements subsequently evolve into Theorems~\ref{thm-AbFibDefprec} and~\ref{thm-HYLkodfibellip}.

\section{Preliminaries and general results}\label{sec-prelim}

In this section, we will recall and prove 
some general results which we need in this article. 
Let us start with some basic terminologies.

\ssec{Basic notions and terminologies}\label{ssec-defnot}
\hfill

In this article, a \emph{complex variety} is a Hausdorff, second-countable, irreducible and reduced complex  space. 
According to this convention, 
complex varieties are locally compact and paracompact. 
A threefold is a compact complex variety of dimension $3$.
A \emph{complex manifold} is a smooth complex variety. 
A subvariety of a complex variety $X$ is a closed reduced analytic subspace of $X$
(which is not necessarily irreducible nor equidimensional).
A \emph{fibration} $f : X \to B$ is a proper holomorphic surjective map with connected fibers. 
The fiber $f^{-1}(b)$ of $f$ over $b \in B$ will often be denoted by $X_b$ if there is no ambiguity.
A fibration $f : X \to B$ is called \emph{isotrivial} if for every $b\in B$,
there exists a neighborhood $U \subset B$ of $b$ (for the Euclidean topology) such that 
$f^{-1}(U) \to U$ is biholomorphic to a 
trivial fibration\footnote{Isotrivial fibrations as we define in this article
are also called \emph{locally trivial fibrations} in some classical references (\eg~\cite{MR184258}). 
In this article, 
since very often we encounter various notions of locally trivial deformations (see Definition~\ref{def-loctriv}), 
in order to avoid possible confusion 
the term \emph{locally trivial} always qualifies a deformation and never a fibration.}. 
We refer to~\cite{VarouchasKS} for a definition of (singular) K\"ahler spaces.
 A \emph{K\"ahler variety} is a complex variety which is K\"ahler. 
A subset $Z \subset X$ of a complex space is \emph{open} (resp. \emph{Zariski open}) if $Z$ is open in $X$ for the Euclidean (resp. Zariski) topology. 
The same convention also applies for dense (resp. Zariski dense) subsets and local (resp. Zariski local) properties.

A \emph{deformation} of a complex space $X$ is a surjective and flat morphism $\Pi : \cX \to \gD$ containing $X$ as a fiber. 
When $X$ is compact, the map $\Pi : \cX \to \gD$ is required to be proper. 
 A \emph{deformation of a holomorphic map $f  : X \to B$} is a composition $\Pi : \cX \xto{q} \cB \xto{\pi} \gD$ where $\Pi$ and $\pi$ are deformations of $X$ and $B$ respectively such that $q_{|\cX_o} : \cX_o \to \cB_o$ equals $f$ for some $o \in \gD$. We say that a deformation of $f$ \emph{preserves $B$} (or \emph{fixes $B$}) if in the above definition, $\pi$ is the projection $B \times \gD \to \gD$. Such a deformation will be denoted by 
$$\Pi : \cX \xto{q} B \times \gD \to \gD.$$

Let $\Pi : \cX \to \gD$ be a deformation of a complex variety $X$. 
Given  a subvariety  $Y \subset X$, we say that $\Pi$ \emph{preserves $Y$} if there exists a subvariety $\cY \subset \cX$ such that $\cY \cap X = Y$ and $\cY$ is isomorphic to $Y \times \gD$ over $\gD$. Similarly, let $\Pi : \cX \xto{q} B \times \gD \to \gD$ be a deformation of $f : X \to B$ fixing $B$ and 
let $Z$ be a subvariety of $B$. We say that $\Pi$ \emph{preserves $Y \cnec f^{-1}(Z) \to Z$} if $q^{-1}(Z \times \gD)$ is isomorphic to $Y \times \gD$ over $B \times \gD$. 

Let $G$ be a group and $X$ a complex space endowed with a $G$-action. We say that a deformation $\Pi : \cX \to \gD$ of $X$ \emph{preserves the $G$-action} (or $\Pi$ is a \emph{$G$-equivariant deformation of $X$}) if there exists a $G$-action on $\cX$ extending the given $G$-action on $X$ such that $\cX \to \gD$ is $G$-invariant (namely, $G$-equivariant with $G$ acting trivially on $\gD$). 
Similarly, let $f : X \to B$ be a $G$-equivariant map. We say that a deformation $\Pi : \cX \xto{q} \cB  \to \gD$ of $f$ \emph{preserves the $G$-action} if there exist  $G$-actions on $\cX$ and on $\cB$ extending the $G$-action on $f:X \to B$ such that $q$ is $G$-equivariant and $\cX \to \gD$ is $G$-invariant.

\begin{Def}[Locally trivial deformations]\label{def-loctriv}
\hfill
\begin{enumerate}[label = \roman{enumi})]

\item\label{def-loctriv0} Let $X$ be a complex space and $Y \subset X$ a subvariety of $X$. A deformation $\cX \to \gD$ of $X$ is called \emph{$Y$-locally trivial}  if there exists a subset $\cU \subset \cX$ such that $U \cnec \cU \cap X$ is a neighborhood of $Y \subset X$ and $\cU \simeq U \times \gD$ over $\gD$.

\item\label{def-loctriv1} A deformation $\cX \to \gD$ of $X$ is called \emph{locally trivial} if it is $\{p\}$-locally trivial for every point $p \in X$.

\item In i), let $G$ be a group acting on a $X$ and assume that $Y$ is $G$-stable. A \emph{$G$-equivariantly $Y$-locally trivial} deformation of $X$ is a $Y$-locally trivial deformation $\cX \to \gD$ of $X$ preserving the $G$-action with the additional property that $\cU \subset \cX$ is $G$-stable and $\cU \simeq U \times \gD$ is $G$-equivariant.
\item\label{def-loctriv2} A deformation $\Pi : \cX \xto{q} B \times \gD \to \gD$ of $f : X \to B$ fixing $B$ is said to be \emph{locally trivial over $B$} if there exists an open cover $\{U_i\}$ of $B$ such that $q^{-1}(U_i \times \gD) \simeq f^{-1}(U_i) \times \gD$ over $U_i \times \gD$ for each $i$. 
\item In iv), let $G$ be a group and let $f : X \to B$ be a $G$-equivariant map. 
We say that $\Pi$ is \emph{$G$-equivariantly locally trivial over $B$} if $\Pi$ preserves the $G$-action and the collection of isomorphisms 
$\left\{q^{-1}(U_i \times \gD) \simeq f^{-1}(U_i) \times \gD\right\}$ above is $G$-invariant for some $G$-invariant open cover $\{U_i\}$ of $B$.
\end{enumerate}
\end{Def}

The name \emph{$Y$-locally trivial} deformation in Definition~\ref{def-loctriv}.i) is justified by the following lemma.

\begin{lem}\label{lem-Yltexpl}
	If $\Pi : \cX \to \gD$ is a $Y$-locally trivial deformation of $X$, 
	then there exist $\cU \subset \cX$ such that $U \cnec \cU \cap X$ is a neighborhood of $Y \subset X$ and an isomorphism $\Phi : \cU \eto U \times \gD$ over $\gD$ such that $\Phi_{|U} : U \to U$ is the identity. In particular, 
	the subvariety $\cY \cnec \Phi^{-1}(Y \times \gD) \subset \cX$ satisfies
	$\cY \cap X = Y$ and $\Phi$ induces an isomorphism of the inclusions  
	\begin{equation}
	\(\cY \hto \cU\) \simeq \(Y \times \Delta \hto U \times \Delta \)
	\end{equation}
	over $\Delta$. Namely, $(\cY \subset \cU) \to \gD$ is a trivial deformation
  of the pair $(Y \subset U)$ over $\gD$ in $\Pi : \cX \to \gD$.
\end{lem}
\begin{proof}
		
	Let $\Phi_0 : \cU \eto U \times \gD$ be the isomorphism in Definition~\ref{def-loctriv}.i) 
	and let ${\Phi_0}_{|U} : U \to U$ be its restriction to $U$. 
	Then the isomorphism 
	$\Phi \cnec ({\Phi_0}_{|U}^{-1} \times \Id_{\gD}) \circ \Phi_0$
	satisfies the property in Lemma~\ref{lem-Yltexpl}.
\end{proof}

The following properties of locally trivial deformations are also easy to see.

\begin{lem}\label{lem-loctrivsm}
\hfill
\begin{enumerate}[label = \roman{enumi})]
\item If $\Pi : \cX \to \gD$ is a smooth (\ie flat and submersive) family, then it is locally trivial.
\item Let $\Pi : \cX \to \gD$ be a deformation of $X$ which is $Y_1$-locally trivial and $Y_2$-locally trivial. If $Y_1 \cap Y_2 = \emptyset$, then up to shrinking $\gD$, $\Pi$ is $(Y_1 \cup Y_2)$-locally trivial.  
\item  If $\Pi : \cX \xto{q} B \times \gD \to \gD$ is a deformation of $f : X \to B$ which is locally trivial over $B$,  then $\Pi : \cX \to \gD$  is a $Y$-locally trivial deformation of $X$ whenever $Y \subset X$ is a subvariety contained in a finite union of fibers of $f$.
\item Let $\Pi : \cX \xto{q} \cB \xto{\pi} \gD$  be a deformation of a map $f : X \to B$.
If $\pi$ is smooth and $q$ is a smooth isotrivial fibration,  
then  for every $b \in B$, up to shrinking $\gD$, $\Pi : \cX \to \gD$ is an $f^{-1}(b)$-locally trivial deformation of $X$. 
\end{enumerate}
\end{lem}
\begin{proof}
	Statements \emph{i)} and \emph{iii)} are trivial.
	For \emph{ii)}, by assumption for $i = 1,2$, we can find $\cU_i \subset \cX$ such that $U_i \cnec \cU_i \cap X$ is a neighborhood of $Y_i \subset X$ and $\cU_i \simeq U_i \times \gD$ over $\gD$. Up to shrinking $\gD$, $\cU_1$, and $\cU_2$, we can assume that $\cU_1 \cap \cU_2 = \emptyset$, so $(\cU_1 \cup \cU_2) \simeq (U_1 \cup U_2) \times \gD$ over $\gD$ and $(\cU_1 \cup \cU_2) \cap X = U_1 \cup U_2$ is a neighborhood of $Y_1 \cup Y_2$. 
	Finally for \emph{iv)}, applying \emph{i)} to $\pi: \cB \to \gD$ there exists  $\cU \subset \cB$ such that 
	$U \cnec \cU \cap B$ is a neighborhood of $b \in B$ and 
	$\cU \simeq U \times \gD$ over $\gD$. Up to shrinking $\gD$ and $\cU$, 
	we can assume that $\cU$ is a sufficiently small neighborhood of $b \in \cB$, 
	so that $q^{-1}(\cU) \to \cU$ and $f^{-1}(U) \to U$ are trivial fibrations in $f^{-1}(b)$ because $q$ is a smooth isotrivial fibration. 
	Hence $q^{-1}(\cU) \cap X = f^{-1}(U)$ is a neighborhood of $f^{-1}(b)$ satisfying
	$$q^{-1}(\cU) \simeq  f^{-1}(b) \times \cU \simeq f^{-1}(b) \times U \times \gD \simeq f^{-1}(U) \times \gD$$
	over $\gD$, which proves \emph{iv)}.
\end{proof}

Another obvious property about locally trivial deformations is that the quotient of a $G$-equivariantly locally trivial deformation is a locally trivial deformation.

\begin{lem}\label{lem-Gquotloctriv}
Let $G$ be a finite group acting on a complex space $X$ and let $Y \subset X$ be a $G$-stable subvariety. If $\Pi : \cX \to \gD$ is a $G$-equivariantly $Y$-locally trivial deformation of $X$, then the quotient $\cX/G \to \gD$ is a $Y/G$-locally trivial deformation of $X/G$. Similarly, if $\Pi : \cX \xto{q} B \times \gD \to \gD$ is a $G$-equivariantly locally trivial deformation over $B$ of a $G$-equivariant map $f : X \to B$, then the quotient $\cX/G \to (B/G) \times \gD \to \gD$ is a deformation of $X/G \to B/G$ which is  locally trivial  over $B/G$.
\end{lem}
\begin{proof}

By assumption, there exists a $G$-stable open subset $\cU$ of $\cX$ containing $Y$ which is $G$-equivariantly isomorphic to $U \times \gD$ over $\gD$ where $U \cnec \cU \cap X$. So the open subset $\cU/G$ of $\cX/G$ is isomorphic to $(U /G) \times \gD$ over $\gD$.  As $Y$ is a $G$-stable subset of $U$, the quotient $Y/G$ is contained in $\cU/G$, which proves that $\cX/G \to \gD$ is a $Y/G$-locally trivial deformation of $X/G$. 
The second statement is proven in~\cite[Lemma 2.2]{HYLbimkod1} with a similar argument.
\end{proof}

\ssec{$G$-equivariant locally trivial deformations}

\hfill

The following lemma shows that if a $G$-equivariant deformation of a complex manifold is $Y$-locally trivial for some $G$-stable subvariety $Y$, then it is $G$-equivariantly $Y$-locally trivial.

\begin{lem}\label{lem-Gloctriv}
Let $X$ be a complex manifold and $G$ a finite group acting on $X$. Let $Y$ be a $G$-stable subvariety of $X$ and let $\Pi : \cX \to \gD$ be a deformation of $X$ which is $G$-equivariant and $Y$-locally trivial. 
Then up to shrinking $\gD$, $\Pi : \cX \to \gD$ is $G$-equivariantly $Y$-locally trivial.
\end{lem}

The proof of Lemma~\ref{lem-Gloctriv} is inspired by the proof of~\cite[Proposition 6.2]{GrafDefKod0}. Before proving Lemma~\ref{lem-Gloctriv}, let us first prove a technical lemma.

\begin{lem}\label{lem-techG}
Let $G$ be a finite group acting on a complex manifold $X$ and let $\Pi: \cX \to \gD$ be a $G$-equivariant deformation of $X$. Let $\cV \subset \cX$ be an open subset such that there exists an isomorphism $\cV \simeq V \times \gD$ over $\gD$ where $V \colonec \cV \cap X$. Let 
$$\cV_G \colonec \bigcap_{g \in G} g(\cV).$$ 
Then for every $G$-stable relatively compact subset $U \subset V_G \colonec \cV_G \cap X$, up to shrinking $\gD$ there exists a $G$-stable subset $\cU \subset\cV_G$ such that $\cU \cap X = U$ and that $\cU$ is $G$-equivariantly isomorphic to $U \times \gD$ over $\gD$.
\end{lem}

\begin{proof}
We may assume that $V_G \ne \emptyset$. Since $\cV_G$ is open by finiteness of $G$, up to shrinking $\gD$ we can assume that the restriction of $\Pi$ to $\cV_G$ is surjective. 
We can also assume that $\gD$ is a polydisc $D(0,\ep_0)^n \subset \bC^n$ centered at $(0,\ldots,0)$ which parameterizes the central fiber $X$ in $\Pi$. Let $t_1,\ldots,t_n$ be the coordinates of $\bC^n$.

For each $i = 1,\ldots,n$, define the homomorphism of Lie algebras
\begin{equation}
\begin{split}
\xi_i : \bC & \to  \Gamma(\cV_G , T_{\cV_G}) \\
z & \mapsto \frac{1}{|G|} \sum_{g \in G} g^* \(\chi_i(z)  _{| \cV_G}\), 
\end{split}
\end{equation}
where $\chi_i (z)$ is the pullback under $\cV \eto V \times \gD$ of the vector field on $V \times \gD$ which projects to $z \cdot \frac{\dr}{\dr t_i}$ in $\gD$ and to $0$ in $V$.  
By~\cite[Satz 3]{Kaup} (see also~\cite[Theorem 4.3]{GrafDefKod0}), there exists a local group action 
$$\Phi_i : \gT_i \to \cV_G$$ 
of $\bC$ on $\cV_G$ inducing $\xi_i$, which means that $\Phi_i$ is a map defined on  a neighborhood  $\gT_i \subset \bC \times \cV_G$ of $\{0\} \times \cV_G$ 
such that
\begin{enumerate}[label = \roman{enumi})]
\item For all $x \in \cV_G$, the subset $\gT_i \cap \(\bC  \times \{x\}\)$ is connected.
\item $\Phi_i(0,\bullet)$ is the identity map on $\cV_G$.
\item $\Phi_i(gh,x) = \Phi_i(g,\Phi_i(h,x))$ whenever it is well-defined.
\item The morphism of Lie algebras $ \bC  \to  \Gamma(\cV_G , T_{\cV_G})$ induced by $\Phi_i$ coincides with $\xi_i$.
\end{enumerate}
Since  the vector field $\xi_i(z)$ is $G$-invariant for all $z \in \bC$, the map $\Phi_i$ is $G$-equivariant (with $G$ acting trivially on $\bC$). Also since $\cV_G \to \gD$ is $G$-invariant, the projection of $\xi_i(z)$ in $\Gamma(\cV_G, \Pi^*T_{\gD})$ is the constant vector field $z \cdot \frac{\dr}{\dr t_i}$. So $\Phi$ descends to a local group action of $\bC$ on $\gD$, namely we have the  commutative diagram
 \begin{equation}\label{diag-locGact}
\begin{tikzcd}[cramped, column sep = 20]
\gT_i \arrow[r, "\Phi_i"] \ar[d] & \cV_G \ar[d, "\Pi"] \\
(\Id_{\bC} \times \Pi)(\gT_i) \ar[r, "\Phi_{\gD,i}"] & \gD
\end{tikzcd}
\end{equation}
where  for every $(z,b) \in (\Id_{\bC} \times \Pi)(\gT_i) \subset \bC \times \gD$, $\Phi_{\gD,i}(z,b) \cnec b + (0,\ldots,0,z,0,\ldots,0)$ with $z$ placed at the $i$-th coordinate. 

In what follows, we define inductively a chain of relatively compact subsets $\cU_0 \subset \cdots \subset \cU_n \subset \cV_G$ and a sequence $\{\gep_i > 0\}_{i = 1,\ldots,n}$  such that $\Pi(\cU_i) = \prod_{j = 1}^i D(0,\gep_j) \times 0^{n-i} \subset \bC^n$ and $\Phi_{i}$ maps $D(0,\gep_{i}) \times \cU_{i-1}$ isomorphically onto its image. Set $\cU_0 \cnec U$, which is relatively compact in $\cV_G$ by assumption. Assume that $\cU_{i-1}$ is constructed. Since $\cU_{i-1}$ is relatively compact in $\cV_G$, there exists $\gep_i >0$ such that 
$$ \ol{D(0,\gep_i)} \times \ol{\cU_{i-1}} \subset \gT_i.$$
We define $\cU_i \cnec \Phi_i (D(0,\gep_i) \times \cU_{i-1})$. Since $\Phi_i$ is a local group action satisfying~\eqref{diag-locGact} and $\Pi(\cU_{i-i}) = \prod_{j = 1}^{i-1} D(0,\gep_j) \times 0^{n-i + 1} $, it follows that $\Pi(\cU_i) = \prod_{j = 1}^i D(0,\gep_j) \times 0^{n-i} $ and  $\Phi_i$ maps  $\ol{D(0,\gep_i)} \times \ol{\cU_{i-1}} $  isomorphically onto its image $\ol{\cU_i}$. The latter shows that $\cU_i$ is relatively compact in $\cV_G$. 

Let $\gep \cnec \min \{\gep_0,\ldots,\gep_n\}$. For every $x \in U$ and $z_1,\ldots,z_n \in D(0,\gep)$, define
$$\Phi(x,z_1,\ldots,z_n) \cnec \Phi_n\(z_n, \Phi_{n-1}(z_{n-1}, \cdots \Phi_1(z_1,x))\). $$
By construction, if $\pr : U \times D(0,\gep)^n \to D(0,\gep)^n$ is the projection, then the diagram 
$$
\begin{tikzcd}[cramped, column sep = 20]
U \times D(0,\gep)^n \arrow[rr, "\Phi"] \ar[dr, "\pr"']& & \cV_G \ar[dl, "\Pi"] \\
& D(0,\gep)^n & 
\end{tikzcd}
$$
commutes and $\Phi$ is an isomorphism onto its image. Since $U$ is $G$-stable and each $\Phi_i$ is $G$-equivariant, $\Phi$ is $G$-equivariant as well. Hence up to shrinking $\gD$ to $D(0,\gep)^n$, the subset $\cU \cnec \Ima(\Phi) \subset \cV_G$ satisfies the desired properties in Lemma~\ref{lem-techG}.
\end{proof}

\begin{proof}[Proof of Lemma~\ref{lem-Gloctriv}]
As the deformation $ \Pi : \cX \to \gD$ is  $Y$-local trivial, 
there exists an open subset  $\cV \subset \cX$ such that $V \colonec \cV  \cap X$ contains $Y$ and that 
 $\cV \simeq V \times \gD$ over $\gD$. 
 Since $Y$ is $G$-stable and since $\cV_G \colonec \bigcap_{g \in G} g(\cV)$, being a finite intersection, is an open subset of $\cX$, it follows that $V_G \colonec \cV_G \cap X$ is a $G$-stable neighborhood of $Y$. Let $U \subset V_G$ be a $G$-stable neighborhood of $Y$ which is relatively compact in $V_G$. Applying Lemma~\ref{lem-techG} to $\cV$ and to $U$, we deduce that up to shrinking $\gD$, there exists a $G$-stable subset $\cU \subset \cV_G$ containing $Y$ which is $G$-equivariantly isomorphic to $U \times \gD$ over $\gD$.
\end{proof}

\ssec{General approaches proving the existence of algebraic approximations}\label{ssec-appgen}

\hfill

In \S\ref{ssec-appgen} we prove Proposition~\ref{pro-red}, which formulates a general approach to finding varieties admitting algebraic approximations. Let us first prove the following.

\begin{lem}\label{lem-defcontrpoint}
	Let $\nu : \ti{X} \to X$ be a map between complex varieties and assume that there exists a subvariety $Y \subset X$ such that $\nu$ maps $\ti{X} \bss  \nu^{-1}(Y)$ isomorphically onto $X \bss Y$.  Then every $Y$-locally trivial deformation $\Pi : \cX \to \Delta$ of $X$ lifts to a deformation $\ti{\cX} \xto{F} \cX \to \gD$ of $\nu$ such that for each $t \in \gD$, the map $\ti{\cX}_t \to \cX_t$ parameterized by $t$ is bimeromorphic.
\end{lem}

\begin{proof}
	
	By assumption, $\nu : \ti{X} \to X$ is the blow-up of $X$ along an ideal sheaf $I$ whose co-support is contained in $Y$.
	Since $\Pi$ is $Y$-locally trivial, by Lemma~\ref{lem-Yltexpl}
	there exists a neighborhood $\cU \subset \cX$ of $Y$ 
	together with an isomorphism $\Phi : \cU \eto U \times \gD$ over $\gD$
	with $U \cnec \cU \cap X$ such that $\Phi_{|U} : U \to U$ is the identity.
	We can thus extend $I$ to the
	the ideal sheaf $\cI$ on $\cX$ satisfying $\cI_{|\cU} = \Phi^*\pr_1^* \cI_{|U}$ where $\pr_1 : U \times \gD \to U$ is the projection.
	As $\cI$ is flat over $\gD$,
	defining $F : \ti{\cX} \to \cX$ as the blow-up of $\cX$ along $\cI$ yields 
	a flat family $\ti{\cX} \to \gD$, which is
	the desired lifting of $\Pi : \cX \to \gD$.
\end{proof}

\begin{rem}
	Lemma~\ref{lem-defcontrpoint} could also be regarded as an 
	immediate corollary of~\cite[Proposition 5.1]{HYLkodfibellip},
	asserting that the same conclusion holds 
	as long as the deformation $\Pi : \cX \to \gD$ 
	preserves the \emph{formal neighborhood} of $Y$.
	While we need~\cite[Proposition 5.1]{HYLkodfibellip}
	to construct algebraic approximations of a compact K\"ahler manifold $X$ 
	with $a(X) = \dim X -1$, 
	in this work Lemma~\ref{lem-defcontrpoint} is already sufficient, 
	and we gave a more elementary proof of it.
\end{rem}

\begin{proof}[Proof of Proposition~\ref{pro-red}]
	Let $\tau : X' \dto X$ be a bimeromorphic map and let
	\begin{equation}\label{map-res1}
	\begin{tikzcd}[cramped, row sep = 0, column sep = 20]
	X' & \ti{X}  \arrow[r, "\eta"] \ar[l, "\nu" , swap] & X \\
	\end{tikzcd}
	\end{equation}
	be a resolution of $\tau$ with $\ti{X}$ smooth. 
	Since $X'$ is normal, the exceptional locus $E \subset \ti{X}$ of $\nu$ 
	satisfies $\dim \nu(E) \le \dim X' - 2$. By assumption, $X'$ has a $\nu(E)$-locally trivial algebraic approximation $\Pi : \cX' \to \gD$. Since $\nu^{-1}_{|X' \bss \nu(E)}$ is an isomorphism onto its image, by Lemma~\ref{lem-defcontrpoint} $\Pi$ lifts to a deformation $\ti{\cX} \to \cX' \to \gD$ of $\nu$. As $\eta : \ti{X} \to X$ is a desingularization of a compact complex variety with at worst rational singularities, up to shrinking $\gD$ the deformation $\ti{\cX} \to \gD$ of $\ti{X}$ induces a deformation $\ti{\cX} \to \cX$ of $\eta$ over $\gD$~\cite[Theorem 2.1]{RanStabMap}. So we have a deformation $\cX' \lto \ti{\cX} \to \cX$ of~\eqref{map-res1} over $\gD$ and for each $t \in \gD$, the fiber $\cX'_t \lto \ti{\cX}_t \to \cX_t$ parameterized by $t$ is bimeromorphic. 
	Since $\cX' \to \gD$ is an algebraic approximation of $X'$, namely Moishezon fibers are dense in the family (see Definition~\ref{def-appalg}), 
	and the property of being Moishezon is invariant under bimeromorphic modifications, 
	$\cX \to \gD$ is an algebraic approximation of $X$.
\end{proof}

Finally, we have the following infinitesimal criterion for the existence of $G$-equivariant algebraic approximations according to~\cite[Theorem 1.11]{GrafDefKod0} (or rather its proof) formulated by Buchdahl and Graf, which is a consequence of Green's density criterion~\cite[Proposition 5.20]{VoisinII}. 

\begin{thm}[{\cite[Proof of Theorem 1.11]{GrafDefKod0}}]\label{thm-densecritG}
	Let $X$ be a compact K\"ahler manifold such that its Kuranishi space is smooth. Let $G$ be a finite group acting on $X$. If there exists a $G$-invariant K\"ahler class $\go \in H^1(X, \gO^1_{X})^G$ such that the contraction $\ctr \go:  H^1(X, T_{X}) \to H^2(X,\cO_X)
	$
	by $\go$ is surjective, then $X$ admits a $G$-equivariant algebraic approximation.
\end{thm}

\ssec{Campana's criterion}\hfill

Let $X$ be a complex variety. We say that $X$ is \emph{algebraically connected} if a general pair of points $x,y \in X$ is contained in a compact connected (but not necessarily irreducible) curve of $X$. 
Since a Moishezon variety always dominates a projective variety of the same dimension, 
Moishezon varieties are algebraically connected.
The converse was proven by Campana for varieties in the Fujiki class $\cC$.

\begin{thm}[Campana {\cite[Corollaire on p.212]{CampanaCored}}]\label{thm-algconn}
Let $X$ be a compact complex variety in the Fujiki class $\cC$. Then $X$ is Moishezon if and only if $X$ is algebraically connected. 
\end{thm}

The following corollaries of Campana's criterion will be useful in this article.

\begin{cor}[Special case of Campana's criterion]\label{cor-multsecMoibase}
Let $f : X \to B$ be a fibration over an algebraically connected variety (\eg a projective curve). Assume that $X$ is in the Fujiki class $\cC$ and the general fiber of $f$ is algebraically connected, then $X$ is Moishezon if and only if $f$ has a multi-section.
\end{cor}

Corollary~\ref{cor-multsecMoibase} implies the following corollaries.

\begin{cor}\label{cor-sltaa}
Let $X$ be a non-algebraic compact complex variety in the Fujiki class $\cC$ and let $f : X \to B$ be a  fibration over a curve with Moishezon fibers.  Let $\Pi : \cX \to B \times \gD \to \gD$ be a deformation of $f$  which is locally trivial over $B$. Then for every subvariety $Y \subset X$ with $\dim Y \le 1$,  the family $\cX \to \gD$ is a $Y$-locally trivial deformation of $X$.
\end{cor}

\begin{proof}

Since $X$ is non-algebraic and since the base and the fibers of $f$ are Moishezon, every subvariety $Y \subset X$ with $\dim Y \le 1$  is contained in a finite number of fibers of $f$ by Corollary~\ref{cor-multsecMoibase}. As $\Pi$ is locally trivial over $B$, by Lemma~\ref{lem-loctrivsm}.\emph{iii)} the underlying deformation $\cX \to \gD$ of $X$ is necessarily $Y$-locally trivial.  
\end{proof}

\begin{cor}\label{cor-critprojP1}
Let $f:X \to B$ be a $\bP^1$-fibration. Assume that $X$ is in the Fujiki class $\cC$ and $B$ is algebraically connected. Then $X$ is Moishezon. 
\end{cor}

\ssec{Good Stein open covers}\label{ssec-gStein}
\hfill

Let $X$ be a complex space and let $\fU \cnec \{U_i\}_{i \in I}$ be an open cover of $X$. Given $i_1,\ldots,i_p$, we define $U_{i_i\cdots i_p} \cnec \cap_{j = 1}^p U_{i_j}$. We say that $\fU$ is a \emph{good} (resp. \emph{Stein}) open cover of $X$ if $\fU$ is locally finite and for every $p \in \bZ_{>0}$ and $i_i,\ldots,i_p \in I$, the open subset $U_{i_i\cdots i_p}$ is contractible (resp. Stein).

\begin{thm}\label{thm-bonStein}
	A good Stein cover exists for any K\"ahler manifold $X$.
\end{thm}

\begin{proof}
	The proof is a refinement of the proof of~\cite[Theorem 5.1]{BottTu}.
	For every $x \in X$, first we choose a Stein neighborhood $V_x \subset X$ of $x$. Let $U_x \subset V_x$ be a geodesically convex neighborhood of $x$ with smooth boundary for the underlying (smooth) Riemannian metric of the K\"ahler manifold $X$. Since $X$ is paracompact, there exists a locally finite refinement $\fU \cnec \{U_i\}_{i \in I}$ of $\{U_x\}_{x \in X}$. As $U_{i_i\cdots i_p} \cnec \cap_{j = 1}^p U_{i_j}$ is geodesically convex for every $p > 0$ and $i_1,\ldots, i_p \in I$, it is contractible. Since $U_{i_i\cdots i_p}$ is a geodesically convex open subset of some neighborhood $V_x$, which is K\"ahler and Stein, it follows from~\cite[Theorem 7]{GreeneWuCurvature} that $U_{i_i\cdots i_p}$ is Stein.
\end{proof}

When $X$ is a Riemann surface, since non-compact open subsets of $X$ are Stein~\cite[p. 134]{GrauertRemmertStein}, we have the following.

\begin{lem}\label{lem-bonStein}
	A good open cover of a Riemann surface $X$ is automatically Stein.
\end{lem}

\ssec{Smooth torus fibrations and tautological families}\label{ssec-fibTlissetaut}
\hfill

In \S\ref{ssec-fibTlissetaut}, we recall the construction of the ($G$-equivariant) tautological family associated to a smooth torus fibration from~\cite{ClaudonToridefequiv} 
and prove some additional results.
Especially, we extend some results in~\cite{ClaudonToridefequiv}
which were stated and proven for compact K\"ahler manifolds,
to manifolds in the Fujiki class $\cC$.

\sssec{Smooth torus fibrations as torsors under the Jacobian fibration}
\hfill

Let $f : X \to B$ be a smooth torus fibration and let $p : J \to B$ be the Jacobian fibration associated to $f$. The sheaf $\cJ_{\bH/B}$ of local holomorphic sections of $J \to B$ lies in the short exact sequence
\begin{equation}\label{SE-Jac}
\begin{tikzcd}[cramped, row sep = 5, column sep = 40]
0 \ar[r] & \bH  \arrow[r]  & \cE_{\bH/B} \ar[r,"\exp_{\bH/B}"] & \cJ_{\bH/B} \ar[r] & 0 
\end{tikzcd}
\end{equation}
where $\bH \colonec R^{2g-1}f_*\bZ$ (with $g$ being the relative dimension of $f$) and 
\begin{equation}\label{def-EHB}
\cE_{\bH/B} \colonec (\bH   \otimes \cO_{B}) / R^{g-1}f_*\gO^g_{X/B} \simeq  R^gf_*\gO^{g-1}_{X/B}  \simeq f_*T_{X/B}.
\end{equation} 
We also use the notations $$\cJ_{X/B} \cnec \cJ_{\bH/B} \ \text{ and } \ \cE_{X/B} \cnec \cE_{\bH/B}.$$
Since all $J$-torsors have the same underlying VHS $\bH$ of weight 1,
up to isomorphisms the sheaves $\cE_{\bH/B},\cJ_{\bH/B}$ and the short exact sequence~\eqref{SE-Jac}
do not depend on the $J$-torsors $f: X \to B$.

As $f : X \to B$ is a $J$-torsor, its isomorphism class determines an element $\eta(f) \in H^1(B,\cJ)$. Conversely, given $\eta \in H^1(B,\cJ)$ and let $(\eta_{ij})$ be a \v{C}ech 1-cocycle with respect to some open cover $\{U_i\}_{i \in I}$ of $B$ representing $\eta$. The fibration $f : X \to B$ obtained by gluing the local fibrations $J_i \cnec p^{-1}(U_i) \to U_i$ along $J_{ij} \cnec p^{-1}(U_{ij}) \to U_{ij} \cnec U_i \cap U_j$ using the translations $\tr(\eta_{ij}) : J_{ij} \to J_{ij}$  by the local sections $\eta_{ij}$ over $U_{ij}$, satisfies $\eta(f) = \eta$.
Under this one-to-one correspondence, 
the $J$-torsors admitting a multi-section correspond bijectively to the torsion elements of $H^1(B,\cJ)$~\cite[Proposition 2.2]{ClaudonToridefequiv}.

Let $p : J \xto{\phi} J' \xto{p'} B$ be a composition of smooth torus fibrations such that both $p$ and $p'$ are Jacobian fibrations. Then $\phi$ induces $\phi_* : \cJ_{J/B} \to \cJ_{J'/B}$, which further induces  $\phi_* : H^1(B,\cJ_{J/B}) \to H^1(B,\cJ_{J'/B})$. 
 
\begin{lem} \label{lem-descent}
	Let $p : J \xto{\phi} J' \xto{p'} B$ be as above.
	Every $J$-torsor $f : X \to B$  
	has a factorization 
	$$f : X \xto{h} X' \xto{f'} B$$ 
	where $f'$ is a $J'$-torsor satisfying 
	$$\eta(f') = \phi_*\eta(f) \in H^1(B,\cJ_{J'/B})$$ 
	and $f^{-1}(b) \xto{h} f'^{-1}(b)$ is isomorphic to $p^{-1}(b) \xto{\phi} p'^{-1}(b)$ for every $b \in B$.
\end{lem}

\begin{proof}
	Let $(\eta_{ij})$ be a \v{C}ech 1-cocycle with respect to some open cover $\{U_i\}_{i \in I}$ of $B$ representing $\eta(f)$. Then $\phi_*\eta(f)$ is represented by $(\phi_*\eta_{ij})$. Since the translations $\tr(\eta_{ij}) : J_{ij} \cnec p^{-1}(U_{ij}) \to J_{ij}$ descend to the translations $\tr(\phi_*\eta_{ij}) : J'_{ij} \cnec p'^{-1}(U_{ij}) \to J'_{ij}$, Lemma~\ref{lem-descent} follows from the reconstruction of the smooth torus fibrations $f$ and $f'$ respectively from $\eta(f)$ and $\phi_*\eta(f)$.
\end{proof}

 Suppose that $f$ is $G$-equivariant for some finite group $G$, then $G$ acts naturally on~\eqref{SE-Jac}, which induces a $G$-action on $J \to B$. The fibration $f$ is a $G$-equivariant $J$-torsor and each isomorphism class of $G$-equivariant $J$-torsors $f$ corresponds bijectively to an element $\eta_G(f) \in H^1_G(B,\cJ_{\bH/B})$ of the $G$-equivariant cohomology group~\cite[\S2.4]{ClaudonToridefequiv}, generalizing the correspondence $f \mapsto \eta(f)$ before.

Following~\cite[Remark 2.3 and \S 2.4]{ClaudonToridefequiv}, here is one way to construct $\eta_G(f) \in H^1_G(B,\cJ_{\bH/B})$.  Recall that the relative Deligne complex is defined as 
$$\ul{D}_{X/B}(g) = \Cone(\bZ \to \gO_{X/B}^{\bullet \le g-1})[-1]$$
where $\bZ$ is regarded as a complex concentrated in degree 0 and the map is induced by the inclusion $ \bZ \xhto{\times (2\pi \sqrt{-1})^g } \gO^\bullet_{X/B}$.
Applying $Rf_*$ to $\ul{D}_{X/B}(g)$, and taking the isomorphisms
$$ \coker \(R^{2g-1}f_*\bZ \to R^{2g-1}f_*\gO^{\bullet \le g-1}_{X/B}\) \simeq \cJ$$
$$ \ker\(R^{2g}f_*\bZ \to R^{2g}f_*\gO^{\bullet \le g-1}_{X/B}\) \simeq \bZ$$
into account, we obtain a short exact sequence
\begin{equation}\label{SE-DeligneJac}
\begin{tikzcd}[cramped]
0  \arrow[r] & \cJ  \arrow[r] & R^{2g}f_*\ul{D}_{X/B}(g) \arrow[r] & \bZ \arrow[r] & 0.
\end{tikzcd}
\end{equation}
Then $\eta_G(f) = \gd_G(1)$ where $\gd_G : \bZ \simeq H_G^0(B,\bZ) \to H^1_G(B,\cJ)$ is the connecting morphism induced by~\eqref{SE-DeligneJac}. Without the $G$-action, $\gd_G$ specializes to 
$$\gd : \bZ \simeq H^0(B,\bZ) \to H^1(B,\cJ)$$ 
and we have $\eta(f) = \gd(1)$.
From the above construction of $f \mapsto \eta_G(f)$, one checks that for every $G$-equivariant map $\phi : B' \to B$, the base change $f' : X \times_B B' \to B'$ is the smooth torus fibration representing  $\phi^*\eta_G(f) \in  H^1_G(B',\cJ_{\phi^{-1}\bH/B'})$  where 
$$\phi^* : H^1_G(B,\cJ_{\bH/B}) \to H^1_G(B,\phi_*\cJ_{\phi^{-1}\bH/B'})  \to H^1_G(B',\cJ_{\phi^{-1}\bH/B'})$$
is the map induced by the map $\phi^* : \cJ_{\bH/B} \to  \phi_*\cJ_{\phi^{-1}\bH/B'}$ pulling back local sections. 

Let $c :  H^1(B,\cJ) \to  H^2(B,\bH)$ 
(resp. $c_G :  H^1_G(B,\cJ) \to  H^2_G(B,\bH)$) be the connecting morphism induced by~\eqref{SE-Jac} (resp. when $f$ is $G$-equivariant for some finite group $G$). 
The following result 
generalizes~\cite[Proposition 2.5 and Proposition 2.11]{ClaudonToridefequiv}
to manifolds in the Fujiki class $\cC$.

\begin{pro}\label{pro-KCequiv}
	Let $f: X \to B$ be a smooth torus fibration over a compact K\"ahler manifold. 
	Assume that $X$ is in the Fujiki class $\cC$. 
	Then $c(\eta(f)) \in H^2(B,\bH)$ is torsion.
	If $f$ is $G$-equivariant for some finite group $G$, 
	then $c_G(\eta_G(f)) \in H_G^2(B,\bH)$ is torsion as well.
\end{pro}
\begin{proof}
	
	The proof of the first statement is similar to~\cite[Proof of Proposition 2.5]{ClaudonToridefequiv}. As $X$ is in the Fujiki class $\cC$, there exists a surjective bimeromorphic morphism $\nu : \ti{X} \to X$ from a compact K\"ahler manifold $\ti{X}$. Let $\ti{f} = f \circ \nu : \ti{X} \to B$. Since $\ti{X}$ is K\"ahler, 
	$$\ti{f}_* : H^{2g}(\ti{X},\bR) \xto{\nu_*} H^{2g}(X,\bR)   \xto{f_*} H^0(B,R^{2g}f_*\bR) \simeq H^0(B,\bR)$$
	is surjective~\cite[Lemma 7.28]{VoisinI}, so $f_*$ is surjective as well. It follows that the Leray differential $d: H^0(B,R^{2g}f_*\bZ) \to H^2(B,R^{2g-1}f_*\bZ) = H^2(B,\bH)$ has torsion image. Since $d = c \circ \gd$~\cite[(2.8)]{ClaudonToridefequiv}, we conclude that $c(\eta(f)) = c(\gd(1))= d(1) \in H^2(B,\bH)$ is torsion.
	
	The same argument proving~\cite[Proposition 2.11]{ClaudonToridefequiv} 
	proves the second statement.
\end{proof}

For later use, we will need the following extension lemma.

\begin{lem}\label{lem-extfibtriv}
	Let $G$ be a finite group and let $B$ be a compact K\"ahler manifold endowed with a $G$-action.
	Let $f : X \to B$ be a $G$-equivariant smooth torus fibration
	such that
	its associated $G$-equivariant Jacobian fibration
	is a trivial fibration $p : B \times T \to B$
	endowed with a diagonal $G$-action on $B \times T$.
	Let $\pi : \cB \to \gD$ be a $G$-equivariant deformation of $B$
	over a contractible Stein space.
	Let $J' \cnec \cB \times T$, endowed with the product $G$-action.
	Assume that $X$ is K\"ahler.
	Then there exists a $G$-equivariant $J'$-torsor $\cX \to \cB$ 
	which lifts $\pi$ to a $G$-equivariant deformation of $f$.
\end{lem}

\begin{proof}
	Let $\bH \cnec R^{2g-1}p_* \bZ$ and $\bH' \cnec R^{2g-1}p'_* \bZ$ 
	where $g = \dim T$ and $p' : \cB \times T \to \cB$ is the projection. 
	As $p'$ is a trivial Jacobian fibration, we have
	$\bH' \simeq \bZ^{2g}$ and 
	$\cE_{\bH'/\cB}  \simeq \cO_{\cB}^g$.
	Consider the commutative diagram
	\begin{equation}\label{dc-extfibtriv}
	\begin{tikzcd}[cramped, row sep = 20, column sep = 20]
	H^1_G(\cB,\cE_{\bH'/\cB}) \arrow[r, "\exp'"] \ar[d, two heads, "r_1"] & 
	H^1_G(\cB ,\cJ_{\bH'/\cB})  \ar[r, "c_G'"] \ar[d, "r_2"]& 
	H_G^2(\cB,\bH') \ar[r]\ar[d, "r_3", "\wr"'] & H_G^2(\cB,\cE_{\bH'/\cB})  \\
	H^1_G(B ,\cE_{\bH/B}) \arrow[r, "\exp"]  & H^1_G(B ,\cJ_{\bH/B})  \ar[r, "c_G"] & H^2_G(B,\bH) 
	\end{tikzcd}
	\end{equation}
	where the rows are the long exact sequences induced by~\eqref{SE-Jac}
	and the vertical arrows are induced by the restriction to $B \subset \cB$.
	Since $\pi : \cB \to \gD$ is a $G$-equivariant deformation of $B$
	and since $G$ is finite, 
	$\cB$ is
	$G$-equivariantly homeomorphic to $B \times \gD$ over $\gD$
	(see \eg \cite[Lemma 2.7]{GalatiusSzucsequivcobord}).
	Since $\bH$ is the restriction of the local system 
	$\bH'$ on $\cB$ to the central fiber of $\pi : \cB \to \gD$ with compatible $G$-actions,  
	the map $r_3$ is an isomorphism
	by homotopy invariance~\cite[Theorem 2.3]{equivariantSS}.
	As $\gD$ is Stein, we have 
	$H^1(\cB,\cE_{\bH'/\cB}) \simeq H^0(\gD,R^1\pi_*\cE_{\bH'/\cB})$. 
	As $R^1\pi_*\cE_{\bH'/\cB}$ is free and
	$h^1(\cB_t, \cE_{\bH'/\cB}) = h^1(\cB_t, \cO^g_{\cB_t})$ is constant in $t \in \gD$,
	Grauert's base change implies that
	$(R^1\pi_*\cE_{\bH'/\cB})_{|\pi(B)} \simeq H^1(B ,\cE_{\bH/B})$, and the restriction
	$$r_1' : H^1(\cB,\cE_{\bH'/\cB}) \simeq H^0(\gD,R^1\pi_*\cE_{\bH'/\cB}) 
	\to (R^1\pi_*\cE_{\bH'/\cB})_{|\pi(B)} \simeq H^1(B ,\cE_{\bH/B})$$ 
	is surjective.
	Since $G$ is finite and $r_1$ in~\eqref{dc-extfibtriv} is isomorphic to 
	$H^1(\cB,\cE_{\bH'/\cB})^G  \xto{r_1'} H^1(B ,\cE_{\bH/B})^G$~\cite[Remark 2.9]{ClaudonToridefequiv},
	is follows that $r_1$ is surjective as well.
	
	As $c_G(\eta_G(f)) \in H^2_G(B,\bH)$ is torsion by Proposition~\ref{pro-KCequiv}, 
	so is $r_3^{-1}(c_G(\eta_G(f))) \in H^2_G(\cB,\bH')$.
	As $H_G^2(\cB,\cE_{\bH'/\cB}) \simeq H^2(\cB,\cE_{\bH'/\cB})^G$ is torsion free, 
	necessarily $r_3^{-1}(c_G(\eta_G(f))) \in \Ima(c_G')$,
	so there exists $\eta \in H^1_G(\cB ,\cJ_{\bH'/\cB})$ 
	such that $c_G(r_2(\eta)) = c_G(\eta_G(f))$.
	Thus there exists $\gb \in H^1_G(B ,\cE_{\bH/B})$ 
	such that $\exp(\gb) = r_2(\eta) - \eta_G(f)$.
	Since $r_1$ is surjective, there exists $\ti{\gb} \in H^1_G(\cB,\cE_{\bH'/\cB})$
	such that $r_1(\ti{\gb}) = \gb$. So
	$$r_2\(\eta - \exp'(\ti{\gb})\) = r_2(\eta) - \exp(r_1(\ti{\gb})) = \eta_G(f).$$ 
	Hence the $G$-equivariant $J'$-torsor $\cX \to \cB$ 
	associated to $\eta - \exp'(\ti{\gb}) \in H^1_G(\cB ,\cJ_{\bH'/\cB})$ 
	is a $G$-equivariant deformation of $f: X \to B$.
\end{proof}

 \sssec{Tautological family}
 \hfill
 
With the notations introduced above, 
let $f: X \to B$ be a $G$-equivariant $J$-torsor 
and let
$$\exp : H^1(B ,\cE_{\bH/B})^G \simeq H^1_G(B ,\cE_{\bH/B}) \to H^1_G(B, \cJ_{\bH/B})$$
 denote the morphism induced by $\exp_{\bH/B} : \cE_{\bH/B} \to \cJ_{\bH/B}$. For any linear map $\gl : V \to H^1(B ,\cE_{\bH/B})^G$ from a finite dimensional vector space $V$, there exists a family
\begin{equation}\label{fam-Jtors}
\Pi : \cX \xto{q}  B \times V \to V 
\end{equation}
of $G$-equivariant $J$-torsors~\cite[Proposition 2.10]{ClaudonToridefequiv} such that  $t \in V$ parameterizes the $G$-equivariant $J$-torsor which corresponds to 
$$\eta_G(f) + \exp(\gl(t)) \in H_G^1(B, \cJ_{\bH/B}).$$ 
The family~\eqref{fam-Jtors} is called the \emph{$G$-equivariant tautological family associated to $f$ parameterized by $V$}. When $B$ is compact and $V = H^1(B ,\cE_{\bH/B})^G$, we call~\eqref{fam-Jtors} the \emph{$G$-equivariant tautological family} associated to $f$. 

Following~\cite[\S 2]{ClaudonToridefequiv}, the family $\Pi$ is constructed as follows.
Let  $\pr_1 : B \times V \to B$ be the first projection and define the $G$-action on $B \times V$ to be the pullback of the $G$-action on $B$ under $\pr_1$. Let
$$\xi \in H^1(B, \cE_{\bH/B})^G \otimes V^\vee \subset H^1(B, \cE_{\bH/B})^G \otimes H^0(V,\cO_V) \simeq  
 H^1(B \times V, \pr_1^*\cE_{\bH/B})^G  \simeq
H^1(B \times V, \cE_{\pr_1^{-1}\bH/B \times V})^G$$ be the element corresponding to $\gl : V \to H^1(B, \cE_{\bH/B})^G$. Let 
$$\gt_G \cnec \pr_1^*\eta_G(f) + \wt{\exp}(\xi) \in  H^1_G(B \times V, \cJ_{\pr_1^{-1}\bH/B \times V})$$
where $\pr_1^* : H^1_G(B, \cJ_{\bH/B}) \to  H^1_G(B \times V, \cJ_{\pr_1^{-1}\bH/B \times V})$ is the map induced by $\cJ_{\bH/B} \to {\pr_1}_*\cJ_{\pr_1^{-1}\bH/B \times V}$ and 
$$ \wt{\exp} : H^1(B, \cE_{\pr_1^{-1}\bH/B \times V})^G \simeq  H^1_G(B, \cE_{\pr_1^{-1}\bH/B \times V}) \to H^1_G(B \times V, \cJ_{\pr_1^{-1}\bH/B \times V})$$
the map induced by $\exp_{\pr_1^{-1}\bH/B \times V} : \cE_{\pr_1^{-1}\bH/B \times V} \to \cJ_{\pr_1^{-1}\bH/B \times V}$. 
Then the smooth torus fibration $q : \cX \to B \times V$ 
in~\eqref{fam-Jtors} defining $\Pi$,
is defined to be the $G$-equivariant $(J \times V)$-torsor corresponding to $\gt_G$.

The tautological families satisfy the following density result. 

\begin{thm}[{Claudon~\cite{ClaudonToridefequiv}}]\label{thm-multsec}
Let $G$ be a finite group and $f : X \to B$ a $G$-equivariant smooth torus fibration. Assume that $X$ is in the Fujiki class $\cC$. Then the base $V = H^1(B ,\cE_{\bH/B})^G$ of the $G$-equivariant tautological family $\Pi$ associated to $f$ contains a dense subset parameterizing fibrations with an étale multi-section. In particular, if the fibers of $f$ and $B$ are Moishezon, 
then $\Pi$ is an algebraic approximation of $f$. 
\end{thm}

\begin{proof}
	Consider the long exact sequence
	\begin{equation}
	\begin{tikzcd}[cramped, row sep = 5, column sep = 20]
	H^1(B,\bH) \ar[r, "\phi"] &  H^1(B ,\cE_{\bH/B}) \arrow[r, "\exp"]  & H^1(B ,\cJ_{\bH/B})  \ar[r, "c"] & H^2(B,\bH)
	\end{tikzcd}
	\end{equation}
	induced by~\eqref{SE-Jac}.
	As $X$ is in the Fujiki class $\cC$,
	$c(\eta(f)) \in H^2(B ,\bH)$ is torsion by Proposition~\ref{pro-KCequiv}, 
	so there exist $m\in \bZ_{>0}$ and 
	$\gb \in H^1(B ,\cE_{\bH/B})$ such that $m \cdot \eta(f) = \exp(\gb)$. 
	As $\eta(f) \in H^1(B ,\cJ_{\bH/B})$ is $G$-invariant and $G$ is finite, 
	up to replacing $m$ by a higher multiple of it,
	we can assume that $\gb \in H^1(B ,\cE_{\bH/B})^G$. 
	Let
	$$W \cnec \left\{t - \frac{\gb}{m} \in H^1(B ,\cE_{\bH/B})^G \ \bigg{|} \ t \in \phi\(H^1(B,\bH)^G\)\otimes \bQ \right\}.$$
	For every $t \in \Ima(\phi) \otimes \bQ$, we have $\exp(t) = 0$, so
	$$m\(\eta(f) + \exp(t - \gb/m)\) = m \cdot \eta(f) - \exp(\gb) = 0,$$ 
	which shows that $\eta(f) + \exp(t - \gb/m) \in H^1(B ,\cJ_{\bH/B})$ is torsion.
	Therefore by~\cite[Proposition 2.2]{ClaudonToridefequiv},
	the torus fibrations parameterized by $W$
	in the tautological family admit a multi-sections.
	
	By~\cite[Lemma 3.2]{ClaudonToridefequiv}, 
	we have $\Ima (\phi) \otimes \bR = H^1(B ,\cE_{\bH/B})$.
	As $G$ is a finite group, we have
	$\phi(H^1(B,\bH)^G)\otimes \bR = H^1(B ,\cE_{\bH/B})^G$, which proves that $W$ is dense in $V = H^1(B ,\cE_{\bH/B})^G$.
	
	The last statement follows from Corollary~\ref{cor-multsecMoibase}.
\end{proof}

\begin{rem}\label{rem-cbfamtaut}
We verify easily from the construction that 
the tautological families are compatible with base change. 
Precisely, if $f' : X \times_B B' \to B'$ 
is the base change of 
$f : X \to B$ by a $G$-equivariant morphism $\phi : B' \to B$, 
then the
$G$-equivariant tautological family associated to $f'$ parameterized by 
\begin{equation}
\begin{tikzcd}[cramped, row sep = 5, column sep = 20]
V \ar[r, "\gl"] &  H^1(B ,\cE_{\bH/B})^G \arrow[r, "\phi^*"]  & H^1(B' ,\cE_{\phi^{-1}\bH/B'})^G 
\end{tikzcd}
\end{equation}
is
$$\Pi' : \cX \times_{B \times V} (B' \times V) \to B' \times V \to V.$$ 	
\end{rem}

Finally, let $f : X \xto{h} X' \xto{f'} B$ be a composition of smooth torus fibrations. Let 
$$
\Pi : \cX \xto{q}  B \times V \to V 
$$
be the tautological family associated to $f$ parameterized by $V \to H^1(B,\cE_{X/B})$.  The morphism $V \to H^1(B,\cE_{X/B}) \to H^1(B,\cE_{X'/B})$ induced by $\cE_{X/B} \simeq f_*T_{X/B} \to f'_*T_{X'/B} \simeq \cE_{X'/B}$ defines a
tautological family
$$
\Pi' : \cX' \xto{q'}  B \times V \to V 
$$
associated to $f'$ parameterized by $V$. The following lemma follows directly from Lemma~\ref{lem-descent} and the construction of tautological families. 
\begin{lem}\label{lem-facttau}
	We have a factorization
	$$ q: \cX \xto{\ti{h}} \cX' \xto{q'} B \times V$$
	of $q$ through $q'$ and for every $(b,v) \in B \times V$, the restriction $q^{-1}(b ,v) \xto{\ti{h}} q'^{-1}(b ,v)$ of $\ti{h}$ is isomorphic to $f^{-1}(b) \xto{h} f'^{-1}(b)$.
\end{lem}

\ssec{Tautological families associated to fibrations in abelian varieties over a curve}\label{ssec-famtautfibellip}
\hfill

Let $f : X \to B$ be a fibration over a smooth projective curve and assume that a general fiber of $f$ is an abelian variety of dimension $g$ (in this article, the case $g = 1$ will be our main interest). We also assume that $X$ is a compact K\"ahler manifold and $f$ has local holomorphic sections at every point of $B$. 
In this setting, 
the construction of tautological families in \S\ref{ssec-fibTlissetaut} 
associated to smooth torus fibrations can be generalized to $f$ 
without the smoothness assumption.
In this paragraph, we
recall the construction of the tautological family associated to $f$ following~\cite[\S4]{HYLbimkod1}. 
This construction will be used in \S\ref{sec-a1} to construct algebraic approximations of
$\bP^1$-fibration over an elliptic surface.

Let us assume, more generally, that there exists a finite group $G$ acting on both $X$ and $B$ such that $f$ is $G$-equivariant.
Let  $ B^\circ \subset B$ be a $G$-stable Zariski dense open subset such that the restriction $f^\circ  : X^\circ  \to B^\circ$ of $f$ to $X^\circ \cnec f^{-1}(B^\circ)$ is smooth.
 Let $J \to B^\circ$ be the Jacobian fibration associated to $f^\circ$. Recall from \S\ref{ssec-fibTlissetaut} that the sheaf $\cJ_{\bH/B^\circ} $ of local holomorphic sections of $J \to B^\circ$ can be identified with the quotient $\cE_{\bH/B^\circ} / \bH$ where $\bH = R^{2g-1}f^\circ_*\bZ$ and $ \cE_{\bH/B^\circ}  \cnec R^gf^\circ_*\gO^{g-1}_{X^\circ/B^\circ} \simeq f^\circ_*T_{X^\circ/B^\circ}$. The Deligne canonical extension of the variation of Hodge structures $\bH$ defines an extension $\bar{\cE}_{\bH/B}$
 of $\cE_{\bH/B^\circ}$, which is  a locally free sheaf on $B$.  
The sheaf $\bar{\cE}_{\bH/B}$ does not depend on the choice of $B^\circ \subset B$, and the induced $G$-action on $\cE_{\bH/B^\circ}$ extends to a $G$-action on $\bar{\cE}_{\bH/B}$.

We recall the construction of the \emph{($G$-equivariant) tautological family} 
\begin{equation}\label{fam-Gtaut}
\Pi :  \cX \xto{q}  B \times V \to V \cnec H^1(B,\bar{\cE}_{\bH/B})^G 
\end{equation}
associated to the $G$-equivariant fibration $f$ following the proof of~\cite[Proposition-Definition 4.18]{HYLbimkod1}. 
Let $\fU \cnec \{U_i\}_{i \in I}$ be a $G$-invariant good Stein open cover of $B$ such that  $f : X \to B$ has local sections over every $U_i$ 
and $U_{ij} \cnec (U_i \cap U_j) \subset B^\circ$ whenever $i \ne j$. Since $V$ is Stein and contractible, $\fU \times V \cnec \{U_i \times V\}_{i \in I}$ is a good Stein open cover of $B \times V$. Let $\pr_1 : B \times V \to B$ be the first projection; by abuse of notation, $\pr_1 : B^\circ \times V \to B^\circ$ will also denote the first projection of $B^\circ \times V$.

For every open subset $U \subset B$ and every section $\gs \in \pr_1^*(\bar{\cE}_{\bH/B}) (U \times V)$, we define for every $t \in V$
$$\gs(t) \cnec \gs_{|U \times \{t\}} \in \bar{\cE}_{\bH/B} (U).$$ 
Since $t \mapsto \gs(t)$ determines $\gs$, we can regard each $\gs \in \pr_1^*(\bar{\cE}_{\bH/B}) (U \times V)$ as a map $\gs : V \to \bar{\cE}_{\bH/B} (U)$. Let
$$\xi \in H^1(B, \bar{\cE}_{\bH/B})^G \otimes H^0(V,\cO_V) \simeq H^1\(B \times V, \pr_1^*\bar{\cE}_{\bH/B}\)^G $$ 
be the element which corresponds to $\Id : V \to H^1(B, \bar{\cE}_{\bH/B})^G$. 
With respect to the good Stein open cover $\fU \times V$, the element $\xi$ can be represented by a $G$-invariant 1-cocycle $\left\{\xi_{ij} : V \to  \bar{\cE}_{\bH/B}(U_{ij}) = \cE_{\bH/B^\circ}(U_{ij})\right\}$ such that $\xi_{ij}(0) = 0$~\cite[Lemma 4.19]{HYLbimkod1}.

Let $X_i = f^{-1}(U_i)$ and $X_{ij} = f^{-1}(U_{ij})$. 
According to the proof of~\cite[Proposition-Definition 4.18]{HYLbimkod1}, 
the map $q : \cX \to B \times V$ in~\eqref{fam-Gtaut} is obtained by gluing the collection of fibrations 
$$(f_{|X_i} \times \Id_V ): X_i \times V \to U_i \times V$$ 
indexed by $i \in I$ along $(f_{|X_{ij}} \times \Id_V) : X_{ij} \times V \to U_{ij} \times V$,  using the $G$-equivariant $1$-cocycle of biholomorphic maps 
$$e_{ij} \cnec \tr\( \exp(U_{ij} \times V)(\xi_{ij}) \) : X_{ij} \times V \eto X_{ij} \times V$$
where 
$$\exp : 
(\pr_1^*\bar{\cE}_{\bH/B})_{|B^\circ \times V} \simeq 
\cE_{\pr_1^{-1}\bH/B^\circ \times V} \xto{\exp_{\pr_1^{-1}\bH/B^\circ \times V}} \cJ_{\pr_1^{-1}\bH/B^\circ \times V}$$ 
and  $\tr(\gs) :  X_{ij} \times V \eto X_{ij} \times V$ denotes the translation of the smooth fibration $ (f_{|X_{ij}} \times \Id_V) : X_{ij} \times V \to U_{ij} \times V$ in abelian varieties by the section $\gs \in \cJ_{\pr_1^{-1}\bH/B^\circ \times V}(U_{ij} \times V)$ of the Jacobian fibration associated to $f_{|X_{ij}} \times \Id_V$. This is how the tautological family~\eqref{fam-Gtaut} is constructed.

We recall the main properties that tautological families satisfy.
\begin{pro}[{\cite[Theorem 1.2 and Proposition-Definition 4.18]{HYLbimkod1}}]\label{pro-def-existfam}
The tautological family 
$$ \Pi :  \cX \xto{q}  B \times V \to V = H^1(B,\bar{\cE}_{\bH/B})^G $$
constructed above satisfies the following properties: 
\begin{enumerate}[label = \roman{enumi})]
\item The central fiber of $\Pi$ is $f$.
\item The family $\Pi$ preserves the $G$-action on $f$ and is $G$-equivariantly locally trivial over $B$.
\item The points parameterizing algebraic members of $\Pi$ form a dense subset of $V$.
\end{enumerate}
\end{pro}

\ssec{Twisted sheaves and Atiyah classes}\label{ssec-twSAtiyah}

\sssec{Twisted sheaves}
\hfill

We will first recall the definition of twisted sheaves. The reader is referred to~\cite[Chapter 1]{CaldararuThese} for more detail on their basic properties.

Let $X$ be a complex space and let $\ga = \{\ga_{ijk}\}$ be a \v{C}ech 2-cocycle with coefficients in $\cO_X^\times$ with respect to an open cover $\{U_i\}_{i \in I}$ of $X$. An \emph{$\ga$-twisted sheaf} $(F,\ga)$ (or simply $F$) is a collection of coherent sheaves $\{F_i \in \Coh(U_i)\}_{i \in I}$ together with $\cO_{U_i \cap U_j}$-linear isomorphisms $g_{ij} : {F_i}_{|U_i \cap U_j} \eto {F_j}_{|U_i \cap U_j}$ such that $g_{ii} = \Id$, $g_{ij} = g_{ji}^{-1}$, and $g_{ki} \circ g_{jk} \circ g_{ij}$ is the multiplication by $\ga_{ijk}$ on ${F_i}_{|U_i \cap U_j \cap U_k}$. 
When $\ga = 0$, an $\ga$-twist sheaf is a coherent sheaf in the usual sense;
we also call it an \emph{untwisted} sheaf.
A morphism $\gt: (F,\ga) \to (F',\ga)$ of $\ga$-twisted sheaves is a collection of morphisms $\gt_i : F_i \to F'_i$ such that $g'_{ij} \circ \gt_i =\gt_j \circ g_{ij}$ over $U_i \cap U_j$ where $F'_i$ are the local pieces and $g'_{ij}$ the gluing isomorphisms defining $(F',\ga)$.
We say that $\gt$ is  $\cO_{X}$-linear if $\gt_i$ is  $\cO_{U_i}$-linear for each $i$.
For every local property $\cP$ attributed to $\cO_X$-modules, we say that a twisted sheaf $(F,\ga)$ satisfies $\cP$ if $F_i$ satisfies $\cP$ for every $i$. For instance if each $F_i$ is locally free, then we say that $(F,\ga)$ is a locally free $\ga$-twisted sheaf.

Coherent $\ga$-twisted sheaves together with $\cO_X$-linear morphisms introduced above form an abelian category $\Coh(X,\ga)$.
For instance, following the above notations,
the kernel of a morphism $\gt: (F,\ga) \to (F',\ga)$ of $\ga$-twisted sheaves 
is the $\ga$-twisted sheaf $\ker(\gt)$ defined by
the coherent sheaves $\ker(\gt)_{|U_i} \cnec \ker(\gt_i)$ over each $U_i$ together
with the isomorphisms $\ker(\gt_i)_{|U_{ij}} \eto \ker(\gt_j)_{|U_{ij}}$
induced by $g_{ij} : {F_i}_{|U_i \cap U_j} \eto {F_j}_{|U_i \cap U_j}$ 
thanks to the condition $g'_{ij} \circ \gt_i =\gt_j \circ g_{ij}$.
If $\ker(\gt) = 0$, namely $\ker(\gt_i) = 0$ for all $i$, 
we say that $(F,\ga)$ is a subsheaf of $(F',\ga)$.
Up to equivalence of categories, $\Coh(X,\ga)$ does not depend on the representative $\ga$ of $[\ga] \in H^2(X,\cO_X^\times)$~\cite[Lemma 1.2.3, Lemma 1.2.8]{CaldararuThese}.
If $\ga'$ is another \v{C}ech 2-cocycle with coefficients in $\cO_X^\times$  with respect to the open cover $\{U_i\}_{i \in I}$  and $F'$ is an $\ga'$-twisted sheaf, then $F \otimes F'$ is well-defined as an $\ga\ga'$-twisted sheaf; 
see~\cite[Proposition 1.2.10]{CaldararuThese} for the detail.

Given an $\ga$-twisted locally free sheaf $(E,\ga)$ of rank $r$ and 
let $g_{ij}$ be the gluing isomorphisms as above. Let $\bar{g}_{ij} : \bP({E_i}_{|U_i \cap U_j}) \to \bP({E_j}_{|U_i \cap U_j})$ be the induced maps on the projectivizations. Since $g_{ki} \circ g_{jk} \circ g_{ij}$ is a scalar multiplication, the collection $\{\bar{g}_{ij}\}$ is a $1$-cocycle of maps and thus defines a $\bP^{r-1}$-bundle over $X$ denoted by $\bP(E,\ga) \to X$, or simply $\bP(E)$. Conversely, every $\bP^{r-1}$-bundle over $X$ is isomorphic to some $\bP(E,\ga) \to X$ and $(E,\ga)$ can be chosen so that $\ga$ is with coefficients in the constant sheaf of $r$th roots of unity $\mu_r \subset \cO_X^\times$~\cite[p.9 and Example 1.1.1]{CaldararuThese}; in particular, $[\ga] \in H^2(X,\cO_X^\times)$ is torsion. 
For every twisted invertible sheaf $L$, local computation shows that $\bP(E) \simeq \bP(E \otimes L)$ over $X$. 
This implies the following result.

\begin{lem}\label{lem-untw}
Assume that an $\ga$-twisted invertible sheaf $L$ exists. 
For every  $\ga$-twisted locally free sheaf $E$,
there exists an  untwisted locally free sheaf $E'$ such that
$\bP(E) \simeq \bP(E')$ over $X$. 
\end{lem}
\begin{proof}
	We have $\bP(E) \simeq \bP(E \otimes L^{\vee})$ over $X$.
	As $E$ is $\ga$-twisted and $L^{\vee}$ is $\ga^{-1}$-twisted, 
	the locally free sheaf $E \otimes L^{\vee}$ is untwisted. 
\end{proof}

\sssec{Jet bundle and Atiyah class}
\hfill

Let $Y \to Z$ be a smooth proper morphism between complex spaces and let
$F$ be a locally free sheaf on $Y$. Recall that the first jet bundle  $J^1(F)$ of $F$ is constructed as follows~\cite[\S4]{Atiyah-Ccnx}. As a sheaf of $\bC$-vector spaces we define $J^1(F) = F \oplus (F \otimes_{\cO_Y} \gO^1_{Y/Z})$. For any local sections $f$ and $(\gs,\tau)$ of $\cO_Y$ and $J^1(F)=  F \oplus (F \otimes_{\cO_Y} \gO^1_{Y/Z})$ respectively, define
\begin{equation}
	f\cdot (\gs,\tau) = (f \cdot \gs, \gs \otimes df + f \cdot \tau).
\end{equation}
This defines an $\cO_Y$-module structure on $J^1(F)$, as well as the short exact sequence of coherent sheaves
\begin{equation}\label{exseq-extJetut}
\begin{tikzcd}[cramped, row sep = 0, column sep = 20]
0 \arrow[r] & F \otimes \gO^1_{Y/Z} \arrow[r] & J^1(F) \arrow[r] & F \arrow[r] & 0.
\end{tikzcd}
\end{equation}
The $\bC$-linear map $j : F \to J^1(F)$ defined by the inclusion $F \oplus 0 \hto F \oplus (F \otimes_{\cO_Y}  \gO^1_{Y/Z})$ is called the \emph{first jet map}. We verify that $j$ satisfies the "Leibniz rule": given local sections $\gs$ and $f$ of $F$ and $\cO_{Y}$ respectively, we have 
\begin{equation}\label{eqn-Leibniz}
j(f \gs) = f \cdot j(\gs) -  \gs  \otimes df.
\end{equation}
where we use~\eqref{exseq-extJetut} to identify $ \gs  \otimes df$ as an element of $J^1(F)$.

We can extend the above notions to the twisted setting. 
Let $\ga$ be a \v{C}ech 2-cocycle with coefficients in $\cO_Y^{\times}$  with respect to the open cover $\{U_i\}$ of $Y$. Let $F$ be an $\ga$-twisted locally free sheaf on $Y$ and let $F_i \cnec F_{|U_i}$ be its (untwisted) local pieces. From the description of the $\cO_{U_i}$-module structure on $J^1(F_i)$, local computation shows that
the $\ga$-twisted sheaf structures on $F$ and $(F \otimes \gO_{Y/Z})$ allow us to glue the $J^1(F_i)$'s and form an $\ga$-twisted sheaf $J^1(F)$, and we call $J^1(F)$ the \emph{first jet bundle} of $F$. The local first jet maps $F_i \to J^1(F_i)$ also glue to a global first jet map $j : F \to J^1(F)$, and $j$ still satisfies the Leibniz rule~\eqref{eqn-Leibniz} over each $U_i$.

By construction, the (twisted) first jet bundle $J^1(F)$ sits in the middle of the short exact sequence of $\ga$-twisted sheaves
\begin{equation}\label{exseq-extJet}
\begin{tikzcd}[cramped, row sep = 0, column sep = 20]
0 \arrow[r] & F \otimes \gO^1_{Y/Z} \arrow[r] & J^1(F) \arrow[r] & F \arrow[r] & 0.
\end{tikzcd}
\end{equation}
and the corresponding element $\At(F/Z) \in \Ext^1(F,F\otimes \gO^1_{Y/Z})$ is called the \emph{Atiyah class} of $F$; 
recall that $\Coh(Y,\ga)$ is an abelian category and
see \eg~\cite[Exercise III.5.2]{GelfandManinHomalg} 
for the interpretation of $\Ext^1$ in terms of extensions in an abelian category.
We write $\At(F) = \At(F/Z)$ if $Z$ is a point. Given two twisted vector bundle $E$ and $F$ on $Y$, 
as in the untwisted case~\cite[Proposition 10]{Atiyah-Ccnx} we verify that $J^1(E \otimes F) \simeq (E \otimes J^1(F)) \oplus (J^1(E) \otimes F)$, so
\begin{equation}\label{eqn-tensAt}
\At(E \otimes F/Z) = \Id_E \otimes \At(F/Z) + \At(E/Z) \otimes \Id_F \in \Ext^1(E \otimes F, E \otimes F\otimes \gO^1_{Y/Z}).
\end{equation}

We will also consider the image $\Tr(\At(F/Z)) \in H^1(Y,\gO^1_{Y/Z})$ of the Atiyah class under the trace map 
$\Tr: \Ext^1(F,F\otimes \gO^1_{Y/Z}) \to H^1(Y,\gO^1_{Y/Z})$. Recall that (see~\cite[p.166]{BuchweitzFlennerACsemireg})  when $F$ is untwisted, we have
\begin{equation}\label{eqn-Atc1}
\Tr(\At(F)) = -c_1(F).
\end{equation}

\ssec{Relative automorphism groups and sheaves of automorphisms}\label{ssec-autrel}
\hfill

Let $\pi : X \to Z$ be a proper surjective morphism between complex varieties. 
Let $\Aut(X/Z)$ be the group of automorphisms of $X$ over $Z$, 
namely automorphisms of $X$ 
mapping each fiber of $X \to Z$ to itself.
Let $\Aut_Z(X)$ be the complex Lie group over $Z$ representing the functor 
\begin{equation}
	\begin{split}
		\text{Complex spaces over } Z  & \to  \text{Sets} \\
		(Z' \to Z) & \mapsto \Aut(X \times_Z Z'/Z');
	\end{split}
\end{equation}
see~\cite{Schuster} or~\cite[\S 1.1]{FujikiholFibdle} for the existence of $\Aut_Z(X)$.  
By definition, we have a structural morphism 
$\pi_A : \Aut_Z(X) \to Z$ and for any $z \in Z$, we have 
$$\pi_A^{-1}(z) = \Aut(X_z)$$ 
where $X_z \cnec \pi^{-1}(z)$.
A point $\ga \in \Aut_Z(X)$ therefore represents an automorphism 
$\ga : X_{\pi_A(\ga)} \eto X_{\pi_A(\ga)} $ 
of the fiber $X_{\pi_A(\ga)}$. 
By definition, we also have 
$$\Aut(X/Z) = \Set{ \text{Sections of } \pi_A : \Aut_Z(X) \to Z}.$$

Let $\Aut^0_Z(X)$ denote the connected component of $\Aut_Z(X)$ containing the identity section $Z \to \Aut_Z(X)$. Let $\cAut_Z(X)$ be the sheaf of local holomorphic sections of the structure morphism $\pi_A : \Aut_Z(X) \to Z$. For every open subset $U \subset Z$, by definition  
$\cAut_Z(X)(U)$ is the group of automorphisms of $X_U \cnec \pi^{-1}(U)$ over $U$.
Similarly, we define $\cAut^0_Z(X)$ to be  the sheaf of local holomorphic sections of the restriction $\pi^0_A : \Aut^0_Z(X) \to Z$ of $\pi_A$ to $\Aut^0_Z(X)$.
 
\begin{lem}
	We have 
	$$\Aut^0_Z(X)_z \cnec (\pi^0_A)^{-1}(z) = \Aut^0(X_z).$$
\end{lem}
\begin{proof}
	
	Since $\Aut^0(X_z)$ is the identity component of $\pi_A^{-1}(z) = \Aut(X_z)$ and $\Aut^0_Z(X)$ is the identity component of the relative Lie group $\pi_A : \Aut_Z(X) \to Z$,
	necessarily
	$\Aut^0_Z(X)_z$ contains $\Aut^0(X_z)$ as a connected component. 
	Thus to prove that $\Aut^0_Z(X)_z = \Aut^0(X_z)$, 
	it suffices to show that $\Aut^0_Z(X)_z$ is connected.
	Let $\Aut^0_Z(X) \to \ti{Z} \to Z$ be the Stein factorization of $\pi_A^0$. The identity section of $\pi_A^0$ induces a section of the finite morphism $\ti{Z} \to Z$, so $\ti{Z}$ contains a connected component which maps isomorphically onto $Z$. Since $\Aut^0_Z(X)$ is connected, $\ti{Z}$ is also connected, so $\ti{Z} \to Z$ is an isomorphism. Hence $\Aut^0_Z(X)_z$ is connected, which implies that 
	$\Aut^0_Z(X)_z  = \Aut^0(X_z)$. 
\end{proof}

\begin{lem}\label{lem-fibAtut}
	Assume that both $\pi : X \to Z$ and $Z$ are smooth, 
	and  $z \mapsto h^0(X_z,T_{X_z})$ is constant in $z \in Z$ where $X_z \cnec \pi^{-1}(z)$. 
	Then $\pi^0_A : \Aut^0_Z(X) \to Z$ is smooth. 
	As a consequence, 
	the Lie algebra over $Z$ associated to the Lie group $\Aut^0_Z(X)$ over $Z$ is the total space $\Tot(\pi_*T_{X/Z})$ of the locally free sheaf $\pi_*T_{X/Z}$. 
\end{lem}

\begin{proof}
	By~\cite[Korollar in p.87]{Kaup},
	the Lie algebra associated to the Lie group $(\pi^0_A)^{-1}(z) = \Aut^0(X_z)$ is $H^0(X_z, T_{X_z})$
	for every $z \in Z$.
	As
	$\dim \Aut^0(X_z)  = h^0(X_z,T_{X_z})$ is constant in $z$
	and $Z$ is smooth, $\pi^0_A : \Aut^0_Z(X) \to Z$ is smooth.
	Finally, since 
	$T_{X/Z}$ is flat over $Z$ (by the smoothness of $\pi$) and $z \mapsto h^0(X_z,T_{X_z})$ is constant,
	the direct image $\pi_*T_{X/Z}$ is locally free,
	which proves the last statement of Lemma~\ref{lem-fibAtut}. 
\end{proof}

Let $Z$ be a complex manifold and let 
\begin{equation}\label{CD-XYZ}
\begin{tikzcd}[cramped, row sep = 15, column sep = 10]
X \ar[dr, swap, "\pi"] \arrow[rr, "f"] & & Y \ar[dl, "p"]   \\
& Z  &    
\end{tikzcd}
\end{equation} 
be a morphism of complex varieties over $Z$.
In the remainder of \S\ref{ssec-autrel}, 
we study the relations between the relative Lie groups of automorphisms arising in~\eqref{CD-XYZ}.
For the purpose of this article, it will be enough to work with
the following assumption.
\begin{assump}\label{assump-autrel} 
	The morphism~\eqref{CD-XYZ} of complex varieties over a complex manifold $Z$ satisfies the following properties.
	\begin{enumerate}[label = \roman{enumi})]
		\item The morphisms $f$ and $p$ are smooth fibrations (with connected fibers) and 
		$$Z \ni z \mapsto h^0(X_z,T_{X_z})$$ 
		$$Z \ni z \mapsto h^0(Y_z,T_{Y_z})$$
		are constant where $X_z \cnec \pi^{-1}(z)$ and $Y_z \cnec p^{-1}(z)$.
		\item For every $z \in Z$, 
		the Lie groups $\Aut^0(X_z)$ and $\Aut^0(Y_z)$ are commutative.
	\end{enumerate}
	
\end{assump}

Until the end of \S\ref{ssec-autrel}, we assume that~\eqref{CD-XYZ} satisfies Assumption~\ref{assump-autrel}. 

By Assumption~\ref{assump-autrel}.\emph{i)} and Lemma~\ref{lem-fibAtut}, $\Tot(\pi_*T_{X/Z})$ (resp. $\Tot(p_*T_{Y/Z})$) is the Lie algebra over $Z$ associated to the Lie group $\Aut^0_Z(X)$ (resp. $\Aut^0_Z(Y)$) over $Z$.
Since $f$ is smooth with connected fibers, we have a descent morphism
$\psi : \Aut^0_Z(X) \to \Aut^0_Z(Y)$
by~\cite[Proposition I.1]{BlanchardAut}.
The induced morphism
$\gT : \pi_*T_{X/Z} \to p_*T_{Y/Z}$ 
coincides with the natural morphism induced by 
$f : X \to Y$.
As $\Aut^0(X_z)$ and $\Aut^0(Y_z)$ are commutative by assumption,
their exponential maps are surjective group homomorphisms.
So we have the commutative diagram 
\begin{equation}\label{dc-explie}
\begin{tikzcd}[cramped, row sep = 20, column sep = 40]
\Tot(\pi_*T_{X/Z}) \ar[d, "\exp_{X/Z}"]\ar[r, ""] & \Tot(p_*T_{Y/Z})\ar[d,"\exp_{Y/Z}"]  \\
\Aut^0_Z(X) \ar[r, "\psi"] & \Aut^0_Z(Y).
\end{tikzcd}
\end{equation}
of morphisms of commutative Lie groups over $Z$ 
with surjective vertical arrows. 
Let
\begin{equation}\label{dc-explief}
\begin{tikzcd}[cramped, row sep = 20, column sep = 40]
\pi_*T_{X/Z} \ar[d]\ar[r, "\gT"] & p_*T_{Y/Z}\ar[d]  \\
\cAut^0_Z(X) \ar[r, "\Psi"] & \cAut^0_Z(Y)
\end{tikzcd}
\end{equation}
be the commutative diagram  of morphisms of sheaves of abelian groups 
which corresponds to~\eqref{dc-explie}.

By the definition of $\psi$,
if $g \in \Aut^0_Z(X)$ (resp. $\psi(g) \in \Aut^0_Z(Y)$)
corresponds to the automorphism $g : X_z \to X_z$ 
(resp. $\psi(g) : Y_z \to Y_z$) where $z = \pi_A(g)$, then the diagram
\begin{equation}\label{dc-descente}
\begin{tikzcd}[cramped, row sep = 20, column sep = 40]
X_z \ar[d,"f"]\ar[r, "g"] & X_z \ar[d,"f"]  \\
Y_z \ar[r, "\psi(g)"] & Y_z
\end{tikzcd}
\end{equation}
commutes. It follows that $\ker \Psi =  p_*\cAut_Y(X) \cap \cAut^0_Z(X)$, 
where we regard $p_*\cAut_Y(X)$ as a subsheaf of $\cAut_Z(X)$.
On the other hand since $\ker \gT = \pi_*T_{X/Y}$, the diagram~\eqref{dc-explief} extends to the commutative diagram
\begin{equation}\label{dc-faisTA}
\begin{tikzcd}[cramped, row sep = 20, column sep = 40]
0  \ar[r] & \pi_*T_{X/Y} \ar[d] \ar[r] & \pi_*T_{X/Z} \ar[d]\ar[r, "\gT"] & p_*T_{Y/Z}\ar[d]  \\
0  \ar[r] & p_*\cAut_Y(X) \cap \cAut^0_Z(X)\ar[r] & \cAut^0_Z(X) \ar[r, "\Psi"] & \cAut^0_Z(Y)
\end{tikzcd}
\end{equation}
of morphisms of sheaves of abelian groups with  exact rows.

\begin{lem}\label{lem-pPsisurj}
	Assume that the morphism~\eqref{CD-XYZ} over $Z$ satisfies Assumption~\ref{assump-autrel}. 
	If the descent morphism
	$\psi:\Aut^0_Z(X) \to \Aut^0_Z(Y)$ induced by~\eqref{CD-XYZ} is surjective, then both $\Psi$ and $\gT$ in~\eqref{dc-faisTA} are surjective.
\end{lem}

\begin{proof}
	
	For simplicity, let $\sA \cnec \Aut^0_Z(X)$ and $\sA' \cnec \Aut^0_Z(Y)$. 
	As $\psi: \sA \to \sA'$ is a surjective morphism of Lie groups over $Z$,
	the associated morphism
	$\Tot(\pi_*T_{X/Z}) \to \Tot(p_*T_{Y/Z})$
	of Lie algebras over $Z$ is also surjective. 
	So $\gT : \pi_*T_{X/Z} \to p_*T_{Y/Z}$ is surjective, 
	and the induced morphism of relative tangent sheaves $T_{\sA/Z} \to \psi^*T_{\sA'/Z}$ is surjective as well.
	Since both $\sA \to Z$ and $\sA' \to Z$ are smooth by Lemma~\ref{lem-fibAtut},
	it follows that  $d\psi : T_{\sA} \to \psi^*T_{\sA'}$ is also surjective.
	Thus $\psi : \sA \to \sA'$ is a submersion
	between complex manifolds, so local sections of $ \sA' \to Z$ 
	lift to local sections of  $\sA \to Z$, which proves that $\Psi$ is surjective.
\end{proof}

\begin{ex}\label{ex-fibab}
	If $f : X \to B$ is a smooth fibration in complex tori of dimension $g$, then $\Aut^0_B(X) \to B$ is isomorphic to the Jacobian fibration $J \to B$ associated to $f$ as Lie groups over $B$. So $\cAut^0_B(X)$ is isomorphic to the sheaf $\cJ$ of local holomorphic sections of $J \to B$. If $\bH = R^{2g-1}f_*\bZ$ and $\cE_{\bH/B} = R^gf_*\gO^{g-1}_{X/B}$, then the quotient map $\cE_{\bH/B} \to \cE_{\bH/B}/\bH = \cJ$ defining $\cJ$ (see \S\ref{ssec-fibTlissetaut}) is isomorphic to the exponential map $\exp : f_*T_{X/B} \to \cAut^0_B(X)$.
\end{ex}

\section{Bimeromorphic models of non-algebraic compact K\"ahler threefolds}\label{sec-bim}

In this section we will prove Propositions~\ref{pro-classk01} and~\ref{pro-classuniregl}, 
which describe 
compact K\"ahler threefolds of algebraic dimension $a \le 1$ 
and non-algebraic uniruled compact K\"ahler threefolds, up to bimeromorphic modifications. 

\ssec{Non-uniruled compact K\"ahler threefolds with $a \le 1$}\label{ssec-gk1}
\hfill

First we prove Proposition~\ref{pro-classk01},
which results essentially from Fujiki's classification~\cite{FujikiStruC}.

\begin{proof}[Proof of Proposition~\ref{pro-classk01}]

Let $X_0$ be a non-uniruled compact K\"ahler threefold with $a(X_0) \le 1$.

 If $a(X_0) = 0$, then by~\cite[Theorem in p.236]{FujikiStruC}, 
 $X_0$ is bimeromorphic to either the quotient $T/G$ of a 3-torus $T$
  by a finite group $G$, or to a simple threefold.  
 As a consequence of the abundance conjecture for K\"ahler threefolds, 
 a simple K\"ahler threefold is also bimeromorphic to a finite quotient of a 3-torus $T$~\cite[Corollary 1.4]{das2023log}.
  In both cases, since $a(X_0) = 0$, necessarily $a(T) = 0$. 
Hence $X_0$ is bimeromorphic to a threefold in Proposition~\ref{pro-classk01}.\emph{i)}.
Note that this is the only place where the MMP for K\"ahler threefolds
occurs in this article.

Assume that $a(X_0) = 1$, then by~\cite[Theorem 3]{FujikiStruC}, 
$X_0$ is bimeromorphic to a threefold $X$ satisfying one of the following descriptions:
\begin{itemize}
	\item  $X = (\ti{S} \times B)/G$ where $\ti{S}$ is a smooth compact K\"ahler surface, 
	$B$ is a smooth projective curve, and $G$ is a finite group acting diagonally on both $\ti{S}$ and $B$. Moreover, a general fiber $F$ of $X \to B/G$
	is a surface with $a(F) = 0$.
	\item $X$ is  the total space of a fibration $f : X \to B$ over a smooth projective curve such that a  general fiber $F$ of $f$ a 2-torus,
	and $f$ satisfies "Property (A)" (namely, for any fibration $g: X' \to B$ and any bimeromorphic map $\phi : X \dto X'$ over $B$, there exists a Zariski open $U \subset B$ such that $\phi$ induces an isomorphism $f^{-1}(U) \simeq g^{-1}(U)$). 
\end{itemize}

Assume that $X = (\ti{S} \times B)/G$ as in the first case above. 
As $F$ is a fiber of $X \to B/G$, we have a finite  map
$\ti{S} \to F$. Since $a(F) = 0$, we have $a(\ti{S}) = 0$. 
As $\ti{S}$ is K\"ahler, by~\cite[Theorem VI.1.1]{Barth}
$a(\ti{S}) = 0$ implies that $\gk(\ti{S}) = 0$ and that the minimal model $S$ of $\ti{S}$ is K-trivial.
By~\cite[Claim in p.99]{Barth}, the $G$-action on $\ti{S}$ induces a holomorphic $G$-action on $S$. 
Hence $X_0$ is bimeromorphic to $(S \times B)/G$,
which is a threefold described in Proposition~\ref{pro-classk01}.$i)$. 

 Given a fibration $f : X \to B$  as in the second case.  By~\cite[Lemma 11.1]{FujikiStruC},  the fact that $f$ has "Property (A)" implies that $f$ has no multi-section. 
 It also implies that if $X' \to X$ is a K\"ahler desingularization of $X$
 then a general fiber of the composition $X' \to X \to B$
 is still a 2-torus.
 Therefore up to replacing $X$ by $X'$, 
 we can assume that $X$ is a compact K\"ahler manifold.
 Assume that $F$ is not algebraic, then by~\cite[Remark 13.1]{FujikiStruC},  
 $X_0$ is bimeromorphic to the quotient the total space of a $G$-equivariant smooth isotrivial torus fibration $\ti{f} : \ti{X} \to \ti{B}$ by $G$. 
 This shows that $f$ satisfies the description in Proposition~\ref{pro-classk01}.$ii)$. 
 Hence the bimeromorphic descriptions of $X$ in Proposition~\ref{pro-classk01} is exhaustive.

Finally, since a variety $X$ in Proposition~\ref{pro-classk01} is 
a finite quotient of a smooth variety, $X$ is normal.
\end{proof}

\ssec{Standard conic bundles}\label{ssec-conicb}
\hfill

In the remainder of \S\ref{sec-bim} we will prove Proposition~\ref{pro-classuniregl}, which describes non-algebraic uniruled K\"ahler threefolds up to bimeromorphic modifications. 
We will first prove some auxiliary results about conic bundles
before we conclude in \S\ref{ssec-prop27}.

We fix some terminology. A \emph{conic bundle} is a $\bP^1$-fibration $f : X \to B$ such that $  f = \pi_{|X}$ where $\pi : \bP(\cE) \to B$ is the projectivization of a locally free sheaf $\cE$ of rank 3 on $B$ and $X$ is the zero locus of a non-trivial section $\gs \in H^0\(\bP(\cE),\cO_{\bP(\cE)}(2) \otimes \pi^*\cL\)$  for some invertible sheaf $\cL$ on $B$. The section  $\gs$ defines a map $\gs : \cE^\vee \to \cE \otimes \cL$, and induces  $\det \gs : \det(\cE)^\vee \to \det(\cE) \otimes \cL^{\otimes 3}$. The divisor in $B$ defined by $\det \gs \in H^0( \det(\cE)^{\otimes 2} \otimes \cL^{\otimes 3})$ is called the \emph{discriminant locus} of the conic bundle $f$. As a set, this is the locus where the quadratic form defined by the restriction of $\gs$ to the fibers of $S^2\cE \otimes \cL$ does not have maximal rank. A flat conic bundle $f : X \to B$ is called \emph{standard} if both $X$ and $B$ are compact complex manifolds, the discriminant locus $D(f) \subset B$ of $f$ is a simple normal crossing divisor, and $f^{-1}(D)$ is irreducible for every irreducible divisor $D \subset B$.

Given an algebraic $\bP^1$-fibration $f : X \to B$, Sarkisov proved that it is always birational to a standard conic bundle~\cite[Proposition 1.13]{SarkisovConicstr}. 
We  extend 
Sarkisov's theorem to non-algebraic $\bP^1$-fibrations over surfaces.

\begin{pro}[\emph{cf.} {Miyanishi~\cite[Theorem on p.89]{MiyanishiConic}} or {Sarkisov~\cite[Proposition 1.13]{SarkisovConicstr}}]\label{pro-bimcon}
Let $f : X \to B$ be a $\bP^1$-fibration over a compact complex surface $B$. There exists a bimeromorphic modification
$$
\begin{tikzcd}[cramped, row sep = 25, column sep = 25]
X \ar[d, swap,"f"] \arrow[r, dashed, "\sim"] & X' \ar[d, "f'"]   \\
B \arrow[r, dashed, "\sim"] & B'   
\end{tikzcd}
$$
of $f$ to a standard conic bundle $f' :X' \to B'$.
\end{pro}

\begin{rem}
When $B$ is a projective surface, Proposition~\ref{pro-bimcon} was already proven in~\cite[Theorem on p.89]{MiyanishiConic} by Miyanishi and Zagorskih. A similar statement was proven by Sarkisov when $B$ is any complete algebraic variety~\cite{SarkisovConicstr}. 
It is likely that
Sarkisov's result continues to hold for any compact complex manifold $B$.
We use the assumption $\dim B = 2$ 
in the proof  to first find a conic bundle (over a compact base) which is bimeromorphic to the original one. Once we obtain a conic bundle bimeromorphic to $f:X \to B$, 
the rest of the argument follows faithfully the proof of Sarkisov's theorem mentioned above. 
\end{rem}

\begin{proof}[Proof of Proposition~\ref{pro-bimcon}]
Up to base-changing $f$ with a desingularization of $B$ and resolving the singularities, we can assume that both $X$ and $B$ are smooth. 
Let $\cE \colonec (f_*K_{X/B}^{\vee})^{\vee\vee}$ where $K_{X/B} \cnec K_X \otimes f^*K_B$. Since $B$ is a smooth surface and $\cE$ is reflexive, $\cE$ is  locally free. As the natural map $f^*f_*K_{X/B}^{\vee} \to K_{X/B}^{\vee}$ is generically surjective, it induces a meromorphic map $\jmath : X \dto \bP(f_*K_{X/B}^{\vee}) \dto \bP(\cE)$ over $B$. Over a Zariski open $U \subset B$ where $f$ is smooth, $\jmath$ is the relative anti-canonical embedding of the $\bP^1$-bundle $f^{-1}(U) \to U$. 
It follows that  $\jmath$  is generically injective  and the proper transform $X_0 \subset \bP(\cE)$ of $X$ by $\jmath$ is the zero locus of a section $\gs \in H^0\(\bP(\cE), \cO_{\bP(\cE)}(2) \otimes \pi^*\cL \)$, where $\pi : \bP(\cE) \to B$ is the projection and $\cL$ is some invertible sheaf over $B$. Therefore up to replacing $X$ by $X_0$, we can assume that $f : X \to B$ is a conic bundle embedded into $\bP(\cE)$,
with $X$ being irreducible. 
 
Starting from the conic bundle  $f:X \to B$ which is embedded in $\bP(\cE)$, the rest of the proof
follows \emph{mutatis mutandis}
the proof of~\cite[Theorem 1.13]{SarkisovConicstr} or~\cite[Theorem on p.89]{MiyanishiConic}.
Let $D(f) \subset B$ be the discriminant divisor of $f$. 
Up to base-changing $f : X \to B$ with a log-resolution of the pair $(B,D(f))$, 
we can assume that $D(f)$ is a simple normal crossing divisor. 
We identify the section $\gs \in  H^0\(\bP(\cE), \cO_{\bP(\cE)}(2) \otimes \pi^*\cL\)$ defining $X$ with a family of quadratic forms 
$q(\gs) \in  H^0(B, S^2\cE \otimes \cL)$ parameterized by $B$. 
We also write $q(\gs) : \cE^\vee \to \cE \otimes \cL$.

\begin{claim}[{\cite[Lemma 1.14]{SarkisovConicstr}}]\label{claim-SC1}
	Up to replacing $f : X \to B$ by a bimeromorphic model $X' \to B$ of it,
	every irreducible component of $D(f)$ is reduced.
\end{claim}

\begin{proof}
Assume that $D(f)$ is not reduced, 
and let $C \subset D(f)$ be the underlying reduced curve of a 
non-reduced irreducible component of $D(f)$.
It suffices to construct a conic bundle $f' : X' \to B$ (with $X'$ irreducible)
bimeromorphic to $f : X \to B$
such that $D(f) - D(f')$ is a nonzero effective divisor;
we will use $D(f) > D(f')$ to denote this property.
We then proceed by induction until every irreducible component of $D(f)$ becomes reduced.

Since $C \subset D(f)$, 
the image of $q(\gs)_{|C} : \cE_{|C}^\vee \to (\cE \otimes \cL)_{|C}$ has rank $r \le 2$.
As $X$ is irreducible, we have $r \ne 0$, so either $r = 1$ or $2$.
Since the image of the composition
$$q' : \cE^\vee \xto{q(\gs)}  \cE \otimes \cL \to (\cE \otimes \cL)_{|C}$$
is a locally free sheaf of rank $r$ over $C$,
the kernel $K^\vee \cnec \ker(q')$ is locally free over $B$; we set $K \cnec K^\vv$.
Consider the composition
\begin{equation}\label{eqn-Kdual}
q_K : K^\vee \hto \cE^\vee \xto{q(\gs)}  \cE \otimes \cL \to  K \otimes \cL
\end{equation}
and also the induced map
$$\gd_K: \cO_B \hto \det(K) \otimes \det(\cE)^\vee \to  \det(K) \otimes \det(\cE) \otimes \cL^{\otimes 3} \to  \det(K)^{\otimes2} \otimes \cL^{\otimes 3}$$
by taking the determinant of~\eqref{eqn-Kdual} and tensorizing with $\det(K)$.

Suppose that $r = 1$.
Since $K^\vee \cnec \ker(q')$, 
$q_K : K^\vee \to  K \otimes \cL$
factorizes through $K \otimes \cL(-C)$ and 
$\gd_K$ factorizes through $\det(K)^{\otimes2} \otimes \cL^{\otimes 3}(-3C)$. 
Let $\gs' \in H^0\(\bP(K), \cO_{\bP(K)}(2) \otimes \pi^*\cL(-C) \)$
be the section such that $q(\gs') :  K^\vee \to K \otimes \cL(-C)$
is the morphism in the factorization of $q_K$.
The discriminant divisor of
the conic bundle $f': X' \to B$ defined by $\gs'$ is the zero locus
$$D(f') = Z\(\cO_B \xto{\gd_K(-3C)} \det(K)^{\otimes2} \otimes \cL^{\otimes 3}(-3C)\)
= Z(\gd_K) -3C$$
of the morphism $\gd_K(-3C)$ in the factorization of $\gd_K$.
Since $\cE^\vee/K^\vee$ is a locally free sheaf of rank $1$ over $C$, we have
$$Z\(\det(K) \otimes \det(\cE) \otimes \cL^{\otimes 3} \to  \det(K)^{\otimes2} \otimes \cL^{\otimes 3}\) = Z\(\cO_B \hto \det(K) \otimes \det(\cE)^\vee\) = C.$$
Thus 
$$D(f') = Z(\gd_K) -3C = C + D(f) + C - 3C = D(f) - C.$$
Finally, removing the irreducible components of $X'$ 
which do not dominate $B$ only decreases $D(f')$.
After removing all such irreducible components, 
the resulting conic bundle $f' : X' \to B$ is bimeromorphic to $f : X \to B$ by construction
with $D(f) > D(f')$.

Now assume that $r = \rk(q(\gs)_{|C}) = 2$.
Then for a general point $o \in C$, 
we can find a neighborhood $U \subset B$ of $o$ and trivializations $\cE \simeq \cO_U^3$ and $\cL \simeq \cO_U$
such that the quadratic form $q(\gs)_{|U}$ on  $\cO_U^3$ is of the form
$$q(\gs) = X_0^2 +  X_1^2 + u^m \cdot X_2^2$$
where $u$ is a holomorphic function on $U$ which defines $C$
and $\{X_i\}_{i = 0,1,2}$ is the standard basis of $\cO_U^3$. 
A section of $K^\vee_{|U} =  \ker(q')_{|U}$ (viewed as a subsheaf of $\cO_U^3$) is of the form
$$s = a_0u \cdot X_0 + a_1u \cdot X_1  + a_2 \cdot X_2$$
where $a_0,a_1,a_2$ are holomorphic functions on $U$.
Over $U$ and through the trivializations $\cE \simeq \cO_U^3$ and $\cL \simeq \cO_U$,
we have
$$q_K (s) = a_0u \cdot X_0 + a_1u \cdot X_1  + a_2u^m \cdot X_2 
= u^2\(\frac{a_0}{u} \cdot X_0 + \frac{a_1}{u} \cdot X_1  + a_2u^{m-2} \cdot X_2\)$$
As $D(f)$ is non-reduced along $C$, we have $m \ge 2$,
showing that  $q_K : K^\vee \to  K \otimes \cL$
factorizes through $K \otimes \cL(-2C)$.
Let $\gs' \in H^0\(\bP(K), \cO_{\bP(K)}(2) \otimes \pi^*\cL(-2C) \)$
be the section such that $q(\gs') :  K^\vee \to K \otimes \cL(-2C)$
is the morphism in the factorization of $q_K$,
and let $f': X' \to B$ be 
the conic bundle defined by $\gs'$.

The determinant $\gd_K$ factorizes through $\det(K)^{\otimes2} \otimes \cL^{\otimes 3}(-6C)$. 
Since $\cE^\vee/K^\vee$ is a locally free sheaf of rank $2$ over $C$, we have
$$Z\(\cO_B \hto \det(K) \otimes \det(\cE)^\vee\) = 2C.$$
The same argument as above then shows that 
$$D(f') = D(f) - 2C,$$
and we finish the construction of $f' : X' \to B$ as before
by removing the irreducible components of $X' $ which do not dominate $B$.
\end{proof}

\begin{claim}[{see also~\cite[Proposition 1.16]{SarkisovConicstr}}]\label{claim-SC2}
	Up to replacing $f : X \to B$ by a bimeromorphic model $f': X' \to B'$ of it,
	the discriminant locus $D(f) \subset B$ satisfies the following.
	 \begin{itemize}
		\item $D(f)$ is a simple normal crossing divisor with reduced irreducible components.
		\item For every $o \in B$, we have $\rk(q(\gs)_o) > 0$ and
		$$
		o \in
		\begin{cases}
			D(f) \setminus \Sing(D(f)) \ \ \text{ if  } \ \ \rk(q(\gs)_o) = 2\\
			\Sing(D(f)) 	   \ \ \text{ if  } \ \ \rk(q(\gs)_o) = 1.
		\end{cases}
		$$
	\end{itemize} 
\end{claim}

\begin{proof}
	By Claim~\ref{claim-SC1}, we can assume that $D(f)$ is a simple normal crossing divisor with reduced irreducible components.
	
	Let $o \in B$. 
	We can find a neighborhood $U \subset B$ of $o$ and trivializations $\cE \simeq \cO_U^3$ and $\cL \simeq \cO_U$
	such that $q(\gs)_{|U} : \cO_U^3 \to \cO_U^3$ is of the form
	$$q(\gs) = \sum_{0 \le i \le  j \le 2} a_{ij} \cdot X_iX_j$$
	where  $\{X_i\}_{i = 0,1,2}$ is the standard basis of $\cO_U^3$
	and each $a_{ij}$ is a holomorphic function on $U$. 
	If $\rk(q(\gs)_o) = 0$, then $a_{ij} = 0$ for all $i$ and $j$,
	so $\det(q(\gs))$ has vanishing order at least $3$ at $o$.
	This contradicts the assumption that $Z\(\det(q(\gs))\) = D(f)$ is a normal crossing curve with reduced irreducible components.
	
	If $\rk(q(\gs)_o) = 1$, then up to changing the trivialization, we have 
	$$q(\gs) = X_0^2 + a_{11} \cdot X_1^2 + a_{12} \cdot X_1X_2 + a_{22} \cdot X_2^2$$
	with $a_{ij}(o) = 0$. So $\det(q(\gs))$ has vanishing order at least $2$ at $o$,
	which implies that $o \in \Sing(D(f))$.
	
	Finally assume that $\rk(q(\gs)_o) = 2$.
	Suppose to the contrary that $o \in \Sing(D(f))$, 
	and let $C_1 , C_2 \subset D(f) \cap U$ be the two branches of $D(f) \cap U$
	passing through $o$,
	defined by $u=0$ and $v= 0$ respectively.
	Then up to changing the trivialization, we have 
	$$q(\gs) = X_0^2 +  X_1^2 + uv \cdot X_2^2.$$
	Let $\tau : \ti{B} \to B$ be the blow-up of $B$ at $o$ and let 
	$q(\ti{\gs}) \cnec \tau^*q(\gs)  \in  H^0(\ti{B}, S^2\tau^*\cE \otimes \tau^*\cL)$.
	Then in the local coordinates $(\ti{u},\ti{v})$ of  $\ti{B}$ which satisfy
	$\tau(\ti{u},\ti{v}) = (\ti{u}\ti{v}, \ti{v})$, we have
	\begin{equation}\label{eqn-qtiloc}
	q(\ti{\gs}) = X_0^2 +  X_1^2 + \ti{u}\ti{v}^2 \cdot X_2^2.
	\end{equation}
	Since the exceptional divisor 
	$E \subset \ti{B}$ of $\tau$ is defined by $\ti{v} = 0$ in that local chart, it follows that
	the discriminant locus $D(\ti{f})$ of the conic bundle $\ti{f} : \ti{X} \to \ti{B}$ defined by $q(\ti{\gs})$ is
	$$D(\ti{f}) = D(f)' + 2E$$
	where $D(f)' \subset \ti{B}$ is the strict transform of $D(f) \subset B$.
	We also see that $q(\ti{\gs})$ has rank $2$ generically along $E$, 
	therefore repeating the same argument as in the proof of Claim~\ref{claim-SC1}
	shows that there exists a bimeromorphic model $\ti{f}' : \ti{X}' \to \ti{B}$ of $\ti{f} : \ti{X} \to \ti{B}$ such that
	$$D(\ti{f}') = D(\ti{f}) - 2E = D(f)'.$$
	Replacing $f: X \to B$ by $\ti{f}' : \ti{X}' \to \ti{B}$ and
	repeating the same procedure until there is no double point $o \in \Sing(D(f))$ such that $\rk(q(\gs)_o) = 2$,
	we finish the proof of Claim~\ref{claim-SC2}.
\end{proof}

Up to replacing $f : X \to B$ by a bimeromorphic model of it,
we can assume that $f$ satisfies the properties in Claim~\ref{claim-SC2}.
These properties imply that $X$ is smooth and that $f$ is flat 
by~\cite[Corollary 1.11]{SarkisovConicstr}, whose proof is of local nature and works also analytically locally.
It remains to show that $f :X \to B$ can be contracted to a standard conic bundle~\cite[\S 1.17]{SarkisovConicstr}. 
Let $C \subset B$ be an irreducible divisor.
If $C \nsubset D(f)$, then $X_{C} \colonec f^{-1}(C) \to C$ is a flat $\bP^1$-fibration, 
so $X_{C}$ is irreducible.
Now suppose that $C \subset D(f)$ is an irreducible component of $D(f)$. 
We claim that if $C \cap \Sing(D(f)) \ne \emptyset$, then 
$X_{C}$ is irreducible. 
Indeed, if $o \in C \cap \Sing(D(f))$, then 
since $f : X \to B$ satisfies the properties in Claim~\ref{claim-SC2},
we can find a neighborhood $U \subset B$ of $o$ trivializing $\cE$ and $\cL$
such that $q(\gs)$ is of the form 
$$q(\gs) = a_0 \cdot  X_0^2 + a_1 \cdot X_1^2 + a_2 \cdot X_2^2$$
where each $a_i$ is a holomorphic function on $U$ (see~\cite[(1.17.1)]{SarkisovConicstr}).
Moreover since $\rk(q(\gs)_o) = 1$, 
we can further assume that $a_0(o) = a_1(o) = 0$, $a_2(o) \ne 0$, 
and $C$ is locally defined by $a_0 = 0$.
 From this local expression of $\gs$, the monodromy action around $o \in C$ on the double cover of $C \setminus \Sing(D(f))$ parameterizing lines in the fibers of $f$ exchanges the two lines in $f^{-1}(p)$ for every $p \in C$, so the total space $X_{C}$ is irreducible. 

Thus if $X_{C} = f^{-1}(C)$ is reducible, 
then $C$ is an irreducible component of $D(f)$ and
 $C \cap \Sing(D(f)) = \emptyset$. 
 In this case, the divisor $X_{C}$ has two irreducible components $E_1$ and $E_2$: both $E_1$ and $E_2$ are ruled surfaces over $C$, and $E_1 \cap E_2$ is transversal and is a section of both $E_1 \to C$ and $E_2 \to C$. It follows that $F \cdot E_1 = -1$ where $F$ is a fiber of $E_1 \to C$. So one can blow down the divisor $E_1$ in $X$ to a curve isomorphic to $C$~\cite[Theorem 2]{FujikiContr} and obtain a new conic bundle which is now smooth along  $C \subset B$. After contracting all such ruled surfaces 
 we obtain a standard conic bundle which is bimeromorphic to the original conic bundle. 
\end{proof}

For later use, let us mention that
the discriminant divisor of a standard conic bundle has the following  property.
 
\begin{lem}[{see \eg \cite[Lemma 5.3]{prokhorov2020finite}}]\label{lem-revd}
	Let $C \subset B$ be the discriminant divisor of a standard conic bundle $f : X \to B$.
	For every irreducible component $D \subset C$ such that $D \simeq \bP^1$, 
	we have $D \cap \ol{C \bss D} \ne \emptyset$ and $\#(D \cap \ol{C \bss D})$ is an even number.
\end{lem}

 \ssec{Curves in a surface of algebraic dimension 0}
 \hfill 
 
Another result we need for proving Proposition~\ref{pro-classuniregl} 
(and also results in \S\ref{sec-casparcas})
is the description of curves of a compact K\"ahler surface $S$ with $a(S) = 0$.
The following result is well-known to experts, see \eg~\cite[Lemma 2.17]{MR2824964} and its proof.
 
  \begin{lem}\label{lem-fini-2}
 	Let $S$ be a smooth compact K\"ahler surface. Assume that $S$ is minimal and $a(S) = 0$, then $S$ contains only finitely many curves. Moreover, the union of curves in $S$ is a finite disjoint union of A-D-E curves, and $S$ has no curve if $S$ is a 2-torus. As a consequence, every connected curve of a  smooth compact K\"ahler surface $\ti{S}$ with $a(\ti{S}) = 0$ is a tree of $\bP^1$.
 \end{lem}

\begin{cor}\label{cor-scba0}
	Let $f : X \to S$ be a standard conic bundle over a smooth compact K\"ahler surface $S$.
	If $a(S) = 0$, then $f$ is smooth, namely $f$ is a $\bP^1$-bundle.
\end{cor}

\begin{proof}
	
	Suppose that $f$ is not smooth.
	Since $a(S) = 0$,
	the discriminant divisor $C \subset S$ 
	is a disjoint union of trees of $\bP^1$ by Lemma~\ref{lem-fini-2}. 
	So there exists an irreducible component $D \subset C$ 
	such that $\#(D \cap \ol{C \bss D}) \le 1$. 
	This contradicts Lemma~\ref{lem-revd} because $D \simeq \bP^1$. 
\end{proof}

 \ssec{Proof of Proposition~\ref{pro-classuniregl}}\label{ssec-prop27}
 \hfill
 
 Before we prove Proposition~\ref{pro-classuniregl},
 we need to define good $\bP^1$-bundles in Proposition~\ref{pro-classuniregl}.i).
 Let $S$ be a smooth compact K\"ahler surface with $a(S) = 0$.
 By Lemma~\ref{lem-fini-2}, 
 $S$ has finitely many curves; 
 let $C_1,\ldots,C_m$ be the connected components of the union of all curves of $S$.
 Consider the composition
 $$\eta : S \xto{\mu} S_\mmin \xto{\nu} S_{\can}$$
 where $\mu$ is the blow-down of $S$ to its minimal model $S_\mmin$,
 and $\nu$ is the contraction of all the A-D-E curves of $S_{\mmin}$. 
 By Lemma~\ref{lem-fini-2}, $\eta(C_i)$ is a point $p_i \in S_\can$ for each $i$,
 and $p_i$ is either a smooth point or a rational double point of $S_\can$.
 
 Let $U_i \subset S_\can$ be a Stein neighborhood of $p_i$. 
 For every reflexive sheaf $\cF$ on $U_i$,
 let 
 $$\cF' \cnec (\nu^*\cF)/\torsion \ \ \ \text{ and } \ \ \ \ti{\cF} \cnec \mu^* \cF',$$ 
 which we regard as sheaves on 
 $V_i \cnec \nu^{-1}(U_i)$ and
 $\ti{V}_i \cnec \eta^{-1}(U_i)$ respectively.
 The sheaf $\ti{\cF}$ satisfies the following properties.

 \begin{lem}\label{lem-Esnault}
 	$\ti{\cF}$ is a globally generated locally free sheaf and satisfies 
 	$R^1\eta_*(\ti{\cF}^\vee) = 0$.
 \end{lem}
 \begin{proof}
 	Since $\cF'$ is locally free and globally generated
 	by~\cite[Lemma and definition (2.2). i) and ii)]{EsnaultRefl}, 
 	so is $\ti{\cF}$.
 	It also follows that $\ti{\cF}^\vee = \mu^* (\cF'^\vee)$.
 	Since $R^1\mu_*\cO_{\ti{V}_i} = 0$, we have
 	$\mu_*(\ti{\cF}^\vee) \simeq \cF'^\vee$ and $R^1\mu_*(\ti{\cF}^\vee) =0$.
 	Thus 
 	$$R^1\eta_*(\ti{\cF}^\vee) \simeq R^1\nu_*(\mu_*(\ti{\cF}^\vee)) \simeq R^1\nu_*(\cF'^\vee) 
 	\simeq R^1\nu_*(\cF'^\vee \otimes \go_{V_i}) = 0,$$
 	where the third isomorphism results from $\go_{V_i} \simeq \cO_{V_i}$
 	because $V_i \to U_i$ is the contraction of the A-D-E curve in $V_i$,  
 	and the last vanishing follows from~\cite[Lemma and definition (2.2). iii)]{EsnaultRefl}.
 \end{proof}

 \begin{Def}\label{Def-bonP1}
 	In the above setting, a \emph{good $\bP^1$-bundle} $f : X \to S$ over $S$
 	is a $\bP^1$-bundle such that for every $i \in \{1,\ldots,m\}$,
 	there exist a neighborhood $U_i \subset S_\can$ of $p_i$
 	and a reflexive sheaf $\cF$ on $U_i$ such that
 	$f^{-1}(\ti{V}_i) \simeq \bP(\ti{\cF})$ over $\ti{V}_i$.
 \end{Def}

\begin{lem}\label{lem-bonP1}
	Every $\bP^1$-bundle $f : X \to S$ over $S$ 
	is bimeromorphic to a good $\bP^1$-bundle over $S$.
\end{lem}
\begin{proof}
	We continue to use the notations introduced above.
	
	We have $H^2(\ti{V}_i,\cO^\times_{\ti{V}_i}) = 0$. Indeed,
	as $U_i$ has only rational singularities, 
	we have 
	$H^2(\ti{V}_i , \cO_{\ti{V}_i}) \simeq H^2(U_i,\eta_*\cO_{\ti{V}_i})$, which vanishes because $U_i$ is Stein.
	As $\ti{V}_i$ deformation retracts to $C_i$, we have 
	$H^3(\ti{V}_i, \bZ) \simeq H^3(C_i, \bZ) = 0$. 
	Hence $H^2(\ti{V}_i,\cO^\times_{\ti{V}_i}) = 0$ by the exponential exact sequence.	
	
	Thus there exists a locally free sheaf $\cE_i$ on $\ti{V}_i$ such that
	$\bP(\cE_i) \simeq f^{-1}(\ti{V}_i)$ over $\ti{V}_i$.
	Now let $\cF_i \cnec (\eta_*\cE_i)^\vv$.
	Since the two natural morphisms
	$$\cE_i \lto \eta^*\eta_*\cE_i \to \ti{\cF_i}$$
	are isomorphisms over $\ti{V}_i \bss C_i$,
	they induce a bimeromorphic map
	$$\phi_i : f^{-1}(\ti{V}_i) \simeq \bP(\cE_i) \dto \bP(\ti{\cF_i})$$
	over $\ti{V}_i$ which is an isomorphism over $\ti{V}_i \bss C_i$.
	As $f^{-1}(\ti{V}_i) \subset X$ is a neighborhood of $f^{-1}(C_i)$
	and we can assume that these $f^{-1}(\ti{V}_i)$ are disjoint up to shrinking $U_i$,
	these bimeromorphic modifications $\phi_i$ extend to 
	a bimeromorphic modification $X \dto X'$ of $X$ over $S$
	which is an isomorphism over $S \bss \cup_{i} C_i$.
	By construction,
	$X' \to S$ is a good $\bP^1$-bundle, which proves the lemma.
\end{proof}

 \begin{proof}[Proof of Proposition~\ref{pro-classuniregl}] 

 By~\cite[Proposition 14.1]{FujikiStruC},  a non-algebraic uniruled compact K\"ahler threefold 
 is bimeromorphic to a threefold $X$ satisfying one of the following descriptions:
 \begin{itemize}
 	\item $X$ is a $\bP^1$-fibration over a surface $S$ with $a(S) = 0$. 
 	\item $X$ is the total space of a fibration $\pi : X \to B$ over a smooth projective curve with the following property: 
 	there exists a Zariski dense open subset $U \subset B$ such that
 	 the fiber $X_b = \pi^{-1}(b)$ is a ruled 
 	surface over an elliptic curve for every $b \in U$ 
 	and that, either $X_b$ is of type $N \in \{I,II\}$ (see Definition~\ref{Def-typee=0} for the terminology) for every $b \in U$, or $X_b$ is of type III for every $b \in U$.
 	Moreover the relative Albanese map $X \xdto{ \ \ \ f \ \ \ } S \xto{p} B$ 
 	of $\pi: X \to B$ exists.
 \end{itemize}

In the first case, we can further assume that $f:X \to S$ is a standard conic bundle by Proposition~\ref{pro-bimcon}, 
so $f$ is a $\bP^1$-bundle by Corollary~\ref{cor-scba0}.
It follows from Lemma~\ref{lem-bonP1} that
$X$ is bimeromorphic to a good $\bP^1$-bundle as in Proposition~\ref{pro-classuniregl}.$i)$. 
Since $S$ is a smooth compact K\"ahler surface, $X$ is K\"ahler. 

In the second case, 
by the definition of the relative Albanese map $f$  
(see~\cite[p. 242]{FujikiStruC}), 
$S$ is smooth and
the restriction of $f$ to a general fiber $X_t$ of $\pi : X \to B$
is holomorphic and is the Albanese map of $X_t$.
Since $X_t$ is a ruled surface over an elliptic curve,
$p: S \to B$ is an elliptic fibration.
Let $U \subset B$ be a Zariski dense open subset as above. 
Up to shrinking $U$, we can assume that the restriction of $p$ 
(resp. $\pi$) to $p^{-1}(U)$ (resp. $\pi^{-1}(U)$) is smooth and that
the restriction of $f$ to $\pi^{-1}(U)$ is holomorphic.
Then there exists a resolution $f' : X' \to S$ of $f : X \dto S$
by a desingularization $\nu : X' \to X$ of $X$ 
such that fundamental locus of $\nu^{-1}$ 
(namely the smallest Zariski closed subset $Z \subset X$ such that $\nu^{-1}_{|X \bss Z}$ is an isomorphism onto its image)
is contained in $X\bss \pi^{-1}(U)$.
Let $\nu' : X'' \to X'$ be a bimeromorphic morphism 
from a compact K\"ahler manifold $X''$ and 
let $E \subset X''$ be the exceptional divisor of $\nu'$.
Then $Z \cnec f'(\nu'(E)) \subset S$ satisfies $\dim Z \le 1$.
Since $X''$ is K\"ahler and non-algebraic,  
$S$ is not algebraically connected by Corollary~\ref{cor-critprojP1},
so $p(Z) \subset B$ is finite.
It follows that the bimeromorphic morphism 
$\nu \circ \nu' : X'' \to X' \to X$ only modifies
$X$ along a finite union of fibers of $\pi : X \to B$.
Therefore up to shrinking $U$ and up to replacing $X$ by the bimeromorphic model $X''$, 
we can assume that $X$ is a compact K\"ahler manifold, and that
$\pi$ has a factorization $\pi : X \xto{f} S \xto{p} B$
such that $p$ is an elliptic surface with $S$ smooth and non-algebraic,
and $X_b$ is still a ruled surface over $S_b \cnec p^{-1}(b)$ 
of type $N \in \{I,II\}$ or of type III uniformly for every $b \in U$.

	It remains to show that there exists $N \in \{I,II,III\}$ 
	such that up to shrinking $U \subset B$, $X_b$ is of type $N$ for every $b \in U$.
	We can assume that $X_b$ is of type $N \in \{I,II\}$ for every $b \in U$ (otherwise $X_b$ is of type III already for every $b \in U$). 
	Since $T_{X/B} = (\gO^1_{X/B})^\vee$ is torsion free and $B$ is a smooth curve, 
	$T_{X/B}$ is flat over $B$~\cite[Exercise 11.8]{MatsumuraCRT}. 
	It follows that
	$$\gS \cnec \Set{b \in B | h^0(X_b,{T_{X/B}}) \ge 4 } \subset B.$$
	is Zariski closed in $B$ by upper semi-continuity. 
	We have the following description of $h^0(X_b,T_{X_b})$.
	
	\begin{lem}[{\cite[Lemma 10]{MaruyamaRuled}}]\label{lem-SeilerDefruled} 
		Let $\sS$ be a ruled surface over an elliptic curve with $\gep(\sS) = 0$. 
		We have
		$$
		h^0(\sS,T_{\sS}) =
		\begin{cases}
		4 \text{ if  } \sS \text{ is of type I,  }\\
		2 \text{ if  } \sS \text{ is of type II or III}.
		\end{cases}
		$$
	\end{lem}
The statements in~\cite[Lemma 10]{MaruyamaRuled}
depend on the invariant $N(\sS)$ of a ruled surface $\sS$ instead of $\gep(\sS)$. 
In order to see that~\cite[Lemma 10]{MaruyamaRuled} contains Lemma~\ref{lem-SeilerDefruled},
here we recall how $N(\sS)$ is defined and show that $N(\sS) = -\gep(\sS)$.
Suppose that $\sS = \bP(\sE)$ for some vector bundle $\sE$ over a smooth projective curve $\sC$,
then $N(\sS) \cnec \deg \det(\sE) - 2\deg(\sL)$ where $\sL \subset \sE$ is a line subbundle such that $\deg(\sL)$ is maximal.
The invariant $N(\sS)$ is independent of the choices of $\sE$ and $\sL$. 
If we choose $\sE$ to be a normalized vector bundle defined in~\cite[Proposition V.2.8 and Notation V.2.8.1]{Hart},
namely $\sE$ satisfies $H^0(\sC,\sE) \ne 0$ and $H^0(\sC,\sE \otimes \sL) = 0$ for any line bundle $\sL$ on $\sC$ with $\deg(\sL)<0$,
then $\cO_{\sC} \subset \sE$ (from $H^0(\sC,\sE) \ne 0$) is a line subbundle of maximal degree, which is $0$.
Thus $N(\sS) = \deg \det(\sE) = -\gep(\sS)$.
In particular, $\gep(\sS) = 0$ is equivalent to $N(\sS) = 0$.

	Assume that $\gS = B$. Then by Lemma~\ref{lem-SeilerDefruled}, $X_b$ is of type I for every $b \in U$. 
	Assume that $\gS \subsetneq B$, then once again by Lemma~\ref{lem-SeilerDefruled} and the assumption that $X_b$ is of type I or II for every $b \in U$,  $U \bss \gS$ is a Zariski open subsets of $B$ parameterizing ruled surfaces $X_b$ of type II. 
	Hence $X$ is bimeromorphic to a $\bP^1$-fibration described in Proposition~\ref{pro-classuniregl}.$ii)$.	
\end{proof}

\section{Fiber bundles in K-trivial surfaces and 3-tori}\label{sec-casparcas}

The main result that we prove in this section is Proposition~\ref{pro-mainpair}.
The aim is to study the existence of algebraic approximations of 
non-uniruled compact K\"ahler threefolds of algebraic dimension $a \le 1$. 
Recall that the good bimeromorphic models 
we choose to approach this problem
are described in Proposition~\ref{pro-classk01}.

Let us start with algebraic approximations of 3-tori.

\ssec{Algebraic approximations  of tori}\label{ssec-Gtores}

\hfill

 It is well-known that complex tori have algebraic approximations. The same proof also works in the $G$-equivariant setting.
\begin{lem}\label{lem-3toriaa}
	A complex torus $X$ endowed with a finite group action $G$ has a $G$-equivariant algebraic approximation. 
\end{lem}
\begin{proof} 
	We can assume that $G$ acts faithfully on $X$.
	Since $X$ is a torus, its Kuranishi space is smooth~\cite[Corollary 2]{MR1144440}. 
	Let $\go \in H^1(X,\gO_X^1)$ be a $G$-invariant K\"ahler class. Since $K_X \simeq \cO_X$,  the contraction $\ctr \go : H^1(X,T_X) \to H^2(X,\cO_X)$ by $\go$ is isomorphic to
	$$\cupp \go : H^1(X,\gO^{n-1}_X) \to H^2(X,\gO^n_X),$$
	which is surjective by the hard Lefschetz theorem. We conclude by Theorem~\ref{thm-densecritG} that $X$ has an algebraic approximation  
	$\pi : \cX \to \gD$
	preserving the $G$-action.
\end{proof}

The main result that we prove in \S\ref{ssec-Gtores} is the following.

\begin{pro}\label{pro-aa3tore}
	Let $X$ be a 3-torus with $a(X) = 0$ and $G$ a finite group acting on $X$. There exists a $G$-equivariant algebraic approximation of $X$ which is $C$-locally trivial for every curve $C \subset X$.
\end{pro}

Let us first prove some auxiliary results.

\begin{lem}\label{lem-aa2torelt}
	Let $S$ be a 2-torus and $G$ a finite group acting on $S$. There exists a $G$-equivariant algebraic approximation of $S$ which is $C$-locally trivial for every proper subvariety $C \subsetneq S$.
\end{lem}
\begin{proof}
	We may assume that $S$ is non-algebraic, so either $a(S) = 0$ or $1$. If $a(S) = 0$, then $S$ has no curve, so Lemma~\ref{lem-aa2torelt} follows from Lemma~\ref{lem-3toriaa} and Lemma~\ref{lem-loctrivsm}.\emph{i)}. 
	If $a(S) = 1$, then the algebraic reduction of $S$ is a smooth elliptic fibration $S \to E$ over an elliptic curve $E$. In this case Lemma~\ref{lem-aa2torelt} follows from Theorem~\ref{thm-multsec} and Corollary~\ref{cor-sltaa}.
\end{proof}

\begin{lem}\label{lem-3tsv}
	Let $X$ be a 3-torus with $a(X) = 0$. If $C \subset X$ is an irreducible curve in $X$, then $C$ is a subtorus of $X$, and every irreducible curve of $X$ is a fiber of $X \to X/C$.
\end{lem}

\begin{proof}
	Choose a point $o \in C$ to be the origin of $X$ and let $T \subset X$ be the subtorus generated by $C$. Since $T$ is a quotient of $\Alb(\ti{C})$ where $\ti{C}$ is the normalization of $C$, $T$ is algebraic, so $T \subsetneq X$. We also have $\dim T \ne 2$, because otherwise $X/T$ is an elliptic curve and $a(X) \ge a(X/T) = 1$. Hence $\dim T = 1$, so necessarily $T = C$. Finally since  $a(X) = 0$, we have $a(X/C) = 0$, so the 2-torus $X/C$ contains no curve by Lemma~\ref{lem-fini-2}. Hence every irreducible curve of $X$ is a fiber of $X \to X /C$.
\end{proof}

\begin{proof}[Proof of Proposition~\ref{pro-aa3tore}]
	If $X$ has no curve, then Proposition~\ref{pro-aa3tore} follows from  Lemma~\ref{lem-3toriaa}. Assume that $X$ contains a curve.
	By Lemma~\ref{lem-3tsv}, there exists a 1-dimensional subtorus $T$ of $X$ such that 
	every irreducible curves of $X$ is a fiber of the 
	smooth isotrivial fibration $f: X \to B \cnec X/T$ in $T$. 
	The $G$-action on $X$ descends to $B$ through $f$. 	
	
	Recall the notations introduced in \S\ref{ssec-fibTlissetaut}.
	As $f : X \to B = X/T$ is the quotient of a torus by the subtorus $T$,
	the Jacobian fibration associated to $f$ is isomorphic to the projection $p : B \times T \to B$, so we have $\bH \cnec R^1f_*\bZ \simeq \bZ^2$  and 	$\cE_{\bH/B} \simeq f_*T_{X/B} \simeq \cO_{B}$.  
	Let $\Phi_g : B \times T \to B \times T$ denote the map defining 
	the induced $G$-action on the Jacobian fibration $p$ 
	for every $g \in G$. The fibration $p$ is $G$-equivariant;
	for every $b \in B$, we define
	$$\Phi_{g,b} : T = \{b\}  \times T \xto{\Phi_g}  \{g(b)\}  \times T = T.$$ 
	Since $a(B) \le a(X) = 0$ and $\dim \Aut(T) = 1$, necessarily the map $B \to \Aut(T)$ defined by $b \mapsto \Phi_{g,b}$ is constant for every $g \in G$. It follows that the $G$-action on $B \times T$ is diagonal; 
	let $\phi$ denote the induced $G$-action on $T$.
	
	By Lemma~\ref{lem-3toriaa}, the torus $B$ has a
	$G$-equivariant algebraic approximation $\pi : \cB \to \gD$. 
	We can assume that $\gD$ is contractible and Stein. 
	Define the $G$-action  
	on $\cB \times T$ as the product of the $G$-action on $\cB$ with $\phi$.
	By Lemma~\ref{lem-extfibtriv}, there exists a $G$-equivariant torsor 
	$\ti{f} : \cX \to \cB$ under $\cB \times T \to \cB$
	which lifts $\pi$ to a $G$-equivariant deformation of $f$.
	
	By definition (see \S\ref{ssec-fibTlissetaut}) 
	 we have $\ti{\cE} \cnec \cE_{\cB \times T/\cB} \simeq \cO_{\cB}$ 
	 and it inherits a $G$-sheaf structure. As $\gD$ is Stein, we have $H^1(\cB, \ti{\cE}) \simeq H^0(\gD,R^1\pi_*\cO_{\cB})$. 
	 Since $R^1\pi_*\cO_{\cB}$ is free and 
	 $h^1(\cB_t,\ti{\cE}) = h^1(\cB_t,\cO_{\cB_t})$ is constant in $t \in \gD$,
	 by Grauert's base change
	 $R^1\pi_*\cO_{\cB}|_t \simeq H^1(\cB_t,\ti{\cE})$ 
	 for every $t \in {\gD}$, 
	 and there exists a finite dimensional subspace 
	 $V \subset H^0(\gD,R^1\pi_*\cO_{\cB})$ such that 
	 the restriction map $V \to H^0(\gD,R^1\pi_*\cO_{\cB}|_{t}) \eto H^1(\cB_t,\ti{\cE})$ is surjective for every $t \in {\gD}$. 
Up to replacing $V$ by $\sum_{g \in G} g \cdot V$, we can assume that $V$ is $G$-stable. As $G$ is finite, the $G$-invariant part 
	$\gl_t:  V^G \to H^1(\cB_t,\ti{\cE})^G$ of the restriction map 
	is still surjective for every $t \in \gD$.
	
	Let (see \S\ref{ssec-fibTlissetaut})
	$$\Pi : \ti{\cX} \xto{q} \cB \times V^G \to V^G$$
	be the $G$-equivariant tautological family associated to  
	$\ti{f} : \cX \to \cB$ and parameterized by 
	$V^G \subset  H^0(\gD,R^1\pi_*\cO_{\cB})^G \simeq H^1(\cB, \ti{\cE})^G$. 
	As tautological families are compatible with base change (Remark~\ref{rem-cbfamtaut}), for every $t \in \gD$ the restriction 
	$$\Pi_t : \ti{\cX_t} \cnec q^{-1}(\cB_t \times V^G) \xto{q} \cB_t \times V^G \to V^G$$
	of $\Pi$ over $\cB_t \times V^G$ is the $G$-equivariant tautological family 
	associated to the $G$-equivariant 
	elliptic fibration $\cX_t \cnec \ti{f}^{-1}(\cB_t) \to \cB_t$ 
	parameterized by $\gl_t : V^G \hto H^1(\cB, \ti{\cE})^G \to H^1(\cB_t,\ti{\cE})^G$. 
	In particular 
	 the family
	$$\Pi': \ti{\cX} \xto{q} \cB \times V^G \xto{\pi \times \Id_{V^G}} \gD \times V^G$$
	of $G$-equivariant elliptic fibrations parameterized by $\gD \times V^G$ contains $f : X \to B$ as a member. 
	As $\gl_t : V^G \to H^1(\cB_t,\ti{\cE})^G$ is surjective for every $t\in \gD$, by Theorem~\ref{thm-multsec} 
	$\Pi_t$ is an algebraic approximation of $\cX_t \to \cB_t$
	for every $t\in \gD$ such that $\cB_t$ is algebraic. Since $\cB \to \gD$ is an algebraic approximation of $B$, it follows that $\Pi'$ is a $G$-equivariant algebraic approximation of $X$.
	
	It remains to show that $\Pi'$ is $C$-locally trivial over a neighborhood of $\Pi'(X) \in \gD \times V^G$. 
	Since $\Pi$ is the tautological family associated to $\ti{f} : \cX \to \cB$, which is a fiber bundle in $T$, the fibration $q : \ti{\cX} \to \cB \times V^G$ is also a smooth isotrivial fibration in $T$. Recall that $C \subset X$ is a finite union of fibers of $f : X \to B$, so $q(C)$ is finite and there exists an open neighborhood $U \subset \cB \times V^G$ of $q(C)$ such that $q^{-1}(U) \simeq U \times T$ over $U$. As $\pi: \cB \to \gD$ is smooth and the image of $q(C)$ under $\cB \times V^G \to V^G$ is $\Pi(X)$, up to shrinking $U$ we can assume that 
	$U$ is of the form $U' \times U''$ where 
	$U' \subset \cB$ is an open subset such that $U' \simeq (U' \cap B) \times \gD$ over $\gD$ and $U'' \subset V^G$ is a neighborhood of $\Pi(X)$.
	 Thus
	$$q^{-1}(U) \simeq U \times T = U' \times U'' \times T \simeq   (U' \cap B) \times \gD \times U'' \times T$$
	over $\gD \times U''$. As $q^{-1}(U)$ is a neighborhood of $C$, 
	it follows that the restriction of $\Pi'$ to  $\Pi'^{-1}(\gD \times U'')$ is a $C$-locally trivial deformation of $X$. 
\end{proof}

\ssec{Fiber bundles in 2-tori}\label{ssec-2tore}

\hfill

Next we study the existence of algebraic approximations for smooth isotrivial 2-torus fibrations.

\begin{lem}\label{lem-isotriv2tore}
	Let $G$ be a finite group and let $f : X \to B$ be a $G$-equivariant  
	smooth isotrivial fibration
	in a non-algebraic 2-torus $F$ over a smooth projective curve $B$. 
	Assume that $X$ is
	in the Fujiki class $\cC$ 
	and $f$ has no multi-section. 
	Then $X$ has a $G$-equivariant algebraic approximation, which is
	$C$-locally trivial for every curve $C \subset X$.
\end{lem}

\begin{proof}
	
	As $F$ is non-algebraic, we have $a(F) = 0$ or $1$. 
	
	If $a(F) = 0$, then $F$ has no curve, so $X$ has no curve neither 
	because $f$ has no multi-section. 
	Thus given a bimeromorphic morphism $\nu: \ti{X} \to X$ from a compact K\"ahler manifold,
	there exists a finite subset $P \subset X$ such that
	$X \bss P$ is isomorphic to $\nu^{-1}(X \bss P)$. 
	As $X \bss P \simeq \nu^{-1}(X \bss P)$ is K\"ahler,
	$X$ is K\"ahler as well~\cite[Corollary of Proposition A]{MiyaokaExtKah}.
	Once again, as $X$ has no curve,
	Lemma~\ref{lem-isotriv2tore} in this case 
	follows from~\cite[Theorem 1.2]{ClaudonToridefequiv}. 
	
	Assume that $a(F) = 1$,  
	then the algebraic reduction of $F$ is a smooth isotrivial fibration in $T$ 
	where $T$ is an elliptic curve, 
	and every irreducible curve of $F$ is a fiber of this fibration. 
	It follows that the fibration $f : X \to B$ admits a fiber-wise $T$-action 
	and by defining $S \cnec X/T$,
	it induces 
	a factorization $X \xto{g} S \xto{h} B$ of $f$ 
	by a smooth isotrivial fibration $g : X \to S$ in $T$. 
	Equivalently, $g : X \to S$ is the universal family of irreducible curves of $X$ contained in the fibers of $f: X \to B$.
	By construction, the $G$-action on $X$ descends to $S$ and both $g$ and $h$ are $G$-equivariant
	smooth isotrivial elliptic fibrations. 
	
	Recall the notations in \S\ref{ssec-fibTlissetaut}.	
     Let $\Pi : \cX \xto{q} S \times V \to V$ 
	be the $G$-equivariant tautological family associated to $g : X \to S$. 
	As $g$ is a smooth isotrivial fibration in $T$, so is $q$. 
	By Theorem~\ref{thm-multsec}, there exists a dense subset $Z \subset V$ such that the fibration $g_v : \cX_v \to S$ parameterized by $v \in Z$ in the family $\Pi$ has a multi-section. It follows from Corollary~\ref{cor-multsecMoibase} that for every $v \in Z$, the fibers of the composition $F_v : \cX_v \xto{g_v} S \xto{h} B$ are algebraic.
	
	Consider the fibration $F: \cX \xto{q} S \times V \xto{h \times \Id}  B \times V$  in 2-tori and let $\pi : B \times V \to V$ be the projection. 
	Let $\cE \cnec \cE_{\cX/B\times V} \cnec F_{*}T_{\cX/B \times V}$ (see \S\ref{ssec-fibTlissetaut}), which inherits a $G$-sheaf structure from the $G$-action on $\cX \to B\times V$. Since $G$ acts trivially on $V$, it acts on $(R^1\pi_*\cE)(U)$ for each open subset $U \subset V$; let $(R^1\pi_*\cE)^G$ denote the subsheaf of $R^1\pi_*\cE$ of local sections fixed by the $G$-action. As $\pi_*\cE$ is coherent and $V$ is Stein, we have $H^1(B \times V, \cE) \simeq H^0(V,R^1\pi_*\cE)$. 
	Since $(R^1\pi_*\cE)^G$ is also coherent, there exists a finite dimensional subspace $W \subset H^0(V,(R^1\pi_*\cE)^G)$ such that the evaluation map $e: W \otimes \cO_V \to (R^1\pi_*\cE)^G$ is surjective. 
	
	Let $\ti{\cX} \xto{\ti{q}} B \times V \times W \to W$ be the $G$-equivariant tautological family associated to $F$ parameterized by 
	$$\gl: W \hto H^0(V,(R^1\pi_*\cE)^G) \simeq H^0(V,R^1\pi_*\cE)^G \simeq H^1(B \times V, \cE)^G.$$ Let 
	$$\ti{\Pi} : \ti{\cX}  \xto{\ti{q}} B \times V \times W \to V \times W $$ 
	be the composition of $\ti{q}$ with the projection $B \times V \times W \to V \times W$. 
	By construction, $\ti{\Pi}$ is a $G$-equivariant deformation of $X$.
	Since tautological families are compatible with base change (Remark~\ref{rem-cbfamtaut}),
	for every $v \in V$ the restriction 
	$$\ti{\Pi}_v : \ti{\cX}_v \to B \times \{v\} \times W \to \{v\} \times W$$ 
	of $\ti{\Pi}$ to  
	$\ti{\cX}_v \cnec \ti{\Pi}^{-1}(\{v\} \times W)$ is isomorphic to the $G$-equivariant tautological family associated to $F_v : \cX_v \to B$ 
	parameterized by 
	$$\gl_v : W \xhto{\gl} H^1(B \times V, \cE)^G  
	\to H^1(B \times \{v\}, \cE)^G \simeq H^1(B,F_{v*}T_{\cX_v/B})^G =  H^1(B, \cE_{\cX_v/B})^G.$$
	By Grauert's theorem, there exists a dense Zariski open subset $V' \subset V$ such that for every $v \in V'$, we have 
	$$(R^1\pi_*\cE)^G|_{v}  \simeq (R^1\pi_*\cE|_{v})^G \simeq H^1(B \times \{v\}, \cE)^G \simeq H^1(B, \cE_{\cX_v/B})^G.$$ 
	Since the evaluation map $e: W \otimes \cO_V \to (R^1\pi_*\cE)^G$ is surjective,
	$W \to (R^1\pi_*\cE)^G|_{v} \simeq H^1(B, \cE_{\cX_v/B})^G$ is surjective for every $v \in V'$. It follows from Theorem~\ref{thm-multsec} applied to the tautological family $\ti{\Pi}_v$
	 that there exists a dense subset $Z_v \subset W$ such that for every $w \in Z_v$, the fibration $\ti{\cX}_{(v,w)} \to B$ parameterized by $w$ in the family $\ti{\Pi}_v$  has a multi-section. 
	 Since the fibers of $\ti{\cX}_{(v,w)} \to B$ are isomorphic to 
	 the fibers of $F_v : \cX_v \to B$, which are algebraic whenever $v \in Z$, it follows from Corollary~\ref{cor-multsecMoibase} that the dense subset 
	 $$\Set{(v,w) \in V \times W | v \in V' \cap Z , w \in Z_v } \subset V \times W$$ 
	 parameterizes fibrations $\ti{\cX}_{(v,w)} \to B$ in the family $\ti{\Pi}$ such that $\ti{\cX}_{(v,w)}$ is algebraic. Thus $\ti{\Pi}$ is a $G$-equivariant algebraic approximation of $X$.
	
	Finally, a similar argument as in the last paragraph of the proof of Proposition~\ref{pro-aa3tore}
	shows that  $\ti{\Pi}$ is 
	$C$-locally trivial over a neighborhood of $\ti{\Pi}(X) \in V \times W$ 
	for every curve $C \subset X$. 
\end{proof}

\ssec{Some precision on the algebraic approximations of $K$-trivial surfaces}\label{ssec-Gsurf}

\hfill

In this paragraph, we study the existence of algebraic approximations of
the product of a K-trivial surface and a curve 
endowed with a diagonal finite group action. 
The main result that we will prove is Proposition~\ref{pro-aaK02}.

To this end, we need to prove a more precise statement about the existence of 
algebraic approximations for $K$-trivial surfaces (see Lemma~\ref{lem-2toreaa}).
Let us start with the following well-known result, 
which we prove due to the lack of an appropriate reference.

\begin{lem}\label{lem-aaK02}
	Let $\Pi : \cS \to B$ be a smooth deformation of a K3 surface or a 2-torus $S$ over a complex manifold $B$. If the Kodaira-Spencer map $\KS : T_{B,b} \to  H^1(S, T_{S})$ at $b = \Pi(S)$ is not zero, then $\Pi$ is an algebraic approximation of $S$.
\end{lem}
\begin{proof}
	
	Lemma~\ref{lem-aaK02} follows from the density criterion~\cite[Proposition 1]{Buch2} 
	which we now explain.
	As $K_S \simeq \cO_S$, the contraction $\lrcorner {\go} : H^1(S,T_S) \to H^2(S,\cO_S)$ by a class $\go \in H^{1}(S,\gO_S)$ is isomorphic to the cup product $H^1(S,\gO^1_S) \to H^2(S,\gO^2_S) \simeq \bC$ with $\go$. 
	As the pairing $H^1(S,\gO^1_S) \otimes H^1(S,\gO^1_S) \to H^2(S,\gO^2_S)$ is non-degenerated, the annihilator of the image of $T_{B,b} \xto{\KS} H^1(S,T_S) \simeq H^1(S,\gO^1_S)$ (which is non-zero by assumption) is a proper subspace of $H^1(S,\gO^1_S)$. 
	As the K\"ahler cone $\cK(S)$ of $S$ is Zariski dense in $H^1(S,\gO^1_S)$, 
	there exists $\go \in \cK(S)$ such that 
	$(\lrcorner {\go}) \circ \KS : T_{B,b} \to H^2(S,\cO_S) \simeq \bC$ is nonzero, hence surjective. 
	Thus Lemma~\ref{lem-aaK02} follows from~\cite[Proposition 1]{Buch2}.
\end{proof}

Let $X$ be a compact K\"ahler manifold and 
let $\NS(X) \subset H^2(X,\bZ)$ be the image of $c_1 : \Pic(X) \to H^2(X,\bZ)$. 
Let $\Pi : \cX \to B$ be a smooth deformation of $X$ over a contractible base $B$. 
We say that $\Pi$ \emph{preserves a subset $V \subset \NS(X)$} if the parallel transport of each class of $V$ under the Gauss-Manin connection remains of type $(1,1)$ along $B$. 
Let $L$ be a line bundle on $X$ with $c_1(L) \in V$.
Assume that $B$ is contractible and $\Pi : \cX \to B$ preserves $V$. Then  for every $v \in V \subset H^2(X,\bZ) \simeq H^2(\cX,\bZ)$, since the map $\Pic^v(\cX/B) \to B$ from the relative Picard variety $\Pic^v(\cX/B)$ parameterizing line bundles $L_b$ on $\cX_b = \Pi^{-1}(b)$  such that  
$c_1(L_b) = v$ under the isomorphism $H^2(\cX_b,\bZ) \simeq H^2(\cX,\bZ)$ 
is surjective and smooth, 
every deformation $\Pi : \cX \to B$ of $X$ preserving $V$ lifts
to a deformation of the pair $(X,L)$,
namely there exists a line bundle $\cL$ on $\cX$ such that $\cL_{|X} \simeq L$.

Now let $S$ be a K3 surface or a 2-torus. Let $\Pi : \cS \to \gD$ be a semi-universal deformation of $S$ and let $o \in \gD$ be the point parameterizing $S$. 
As $K_S \simeq \cO_S$, the base $\gD$ is smooth
by the Bogomolov-Tian-Todorov theorem~\cite[Corollary 2]{MR1144440}
and $\dim \gD = h^{1}(S,T_S) = h^{1,1}(S)$.
Let 
$$\KS : T_{\gD,o} \to H^1(S,T_S)$$
denote the Kodaira-Spencer map.
For every subgroup $V \subset \NS(S)$, let $\gD_{V} \subset \gD$ be the largest locus over which $V$ is preserved under the induced deformation of $S$. 
By~\cite[Lemma 10.19, Theorem 10.21]{VoisinI} and~\cite[Lemma 5.16]{VoisinII},
the image of the Kodaira-Spencer map of the induced deformation of $S$ over $\gD_V$ 
is 
\begin{equation}\label{eqn-KSgDV}
\KS(T_{\gD_{V},o}) = \Set{\gt \in H^1(S,T_S) | \gt \lrcorner v = 0 
\in H^2(S,\cO_S) \text{ for every } v \in V \subset H^1(S,\gO^1_S) }.
\end{equation}
Here $T_{\gD_{V},o}$
is the tangent space of the analytic subspace $\gD_{V}$ at $o$, namely
$$T_{\gD_{V},o} = 
\Set{v \in T_{\gD,o} | df(v) = 0, \text{ for all } f \in \cI_{\gD_{V},o} }$$
where $\cI_{\gD_{V},o} \subset \cO_{\gD,o}$ is the ideal of germs of holomorphic functions at $o$ defining $\gD_V$. 

\begin{lem}\label{lem-lisseDV}
	The locus $\gD_V$ is smooth of dimension $h^{1,1}(S) - \rk (V)$.
\end{lem}
\begin{proof}
	Our argument follows~\cite[\S 6.2.4]{HuybK3book} where the case $\rk (V) = 1$ is proven.
 Let $\cP : \gD \to \bP(H^2(S,\bC))$\footnote{In this proof, the projectivization $\bP(E)$ of a vector space $E$ is defined to be the space of 1-dimensional subspaces of $E$.} be the period map  sending $t \in \gD$ to the line $H^{2,0}(\cS_t) \subset 
	H^2(S,\bC)$ where $\cS_t = \Pi^{-1}(t)$ and we identify $H^2(\cS_t,\bC)$ with $H^2(S,\bC)$ using the Gauss-Manin connection. By the local Torelli theorem~\cite[Proposition 6.2.8]{HuybK3book} (which also holds when $S$ is a 2-torus with the same proof), $\cP$ is an open immersion into the period domain
	$$\cD = \Set{[\bC \eta] \subset \bP(H^2(S,\bC)) | Q(\eta,\eta) = 0, Q(\eta,\bar{\eta}) > 0}$$
	where $Q$ is the intersection product on $H^2(S,\bC)$.
	For every $t \in \gD$, let $\gs_t$ be a generator of $H^{2,0}(\cS_t)$, 
	viewed as an element in $H^2(S,\bC)$ under the parallel transport by the Gauss-Manin connection. 
	By definition of $\gD_V$, we have $t \in \gD_V$ if and only if $Q(\gs_t , v) = 0$ for every $v \in V$, so $\cP$ maps $\gD_V$ isomorphically onto an open subset of the intersection of $\cD$ with  $\bP(V^\perp) \subset \bP(H^2(S,\bC))$ where $V^\perp \subset H^2(S,\bC)$ is defined with respect to $Q$. 
	For every $t \in \gD_V$, since $\ol{\gs_t} \in H^{0,2}(\cS_t) \subset  H^{1,1}(\cS_t)^\perp \subset V^\perp$ and $Q(\gs_t, \ol{\gs_t}) \ne 0$, the linear form $Q(\gs_t, \bullet) $ is non-zero on $V^\perp$. As $\cD$ is defined locally by the quadratic form $Q$, the non-vanishing of $Q(\gs_t, \bullet)$ on $V^\perp$ implies that $\cD \cap \bP(V^\perp)$, and thus $\gD_V$, is smooth at $t$. 
	
	Let $\gs : H^1(S,T_S) \to H^1(S,\gO^1_S)$ be the map defined by the
	contraction with a holomorphic symplectic form $\gs$ on $S$, which is an isomorphism.
	By~\eqref{eqn-KSgDV}, the isomorphism $\gs \circ \KS$ 
	sends $T_{\gD_{V},o}$ to $V^\perp \cap H^1(S,\gO_S^1)$. 
	Since $\gD_V$ is smooth, it follows that 
	$\dim \gD_V = \dim T_{\gD_{V},o}= 
	\dim\(V^\perp \cap H^1(S,\gO_S^1)\) = h^{1,1}(S) - \rk (V)$. 
\end{proof}

Now let $G$ be a finite group acting on the surface $S$. By~\cite[Theorem 5.2]{doan2020equivariant}, up to shrinking $\gD$ the $G$-action on $S$ extends to $\cS$ and descends to a $G$-action on $\gD$ through $\Pi$. Moreover, the Kodaira-Spencer map $\KS : T_{\gD,o} \to H^1(S,T_S)$ at $o \cnec \Pi(S)$ is $G$-equivariant with respect to the induced $G$-actions. 

Assume that $V \subset \NS(S)$ is $G$-stable, then $\gD_V \subset \gD$ is also $G$-stable. When $V = \NS(S)$ we have the following.

\begin{lem}\label{lem-Gtriv}
Assume that $S$ is a non-algebraic $K$-trivial surface. The induced $G$-action on $\gD_{\NS(S)}$ is trivial.
\end{lem}

\begin{proof}
		
		As $\gD_{\NS(S)}$ is smooth by Lemma~\ref{lem-lisseDV} and $G$ is finite, the $G$-action on $\gD_{\NS(S)}$ is linearizable~\cite[Lemme 1]{CartanLineaire}, so it suffices to show that $G$ acts as the identity on $T_{\gD_{\NS(S)},o}$.
		
		 Since $G$ is finite,  the $G$-action on $S$ is symplectic
		 by~\cite[Corollary 15.1.10.(ii)]{HuybK3book} 
		 (which also holds when $S$ is a 2-torus by the same proof), 
		 so the contraction $T_S \to \gO_S^1$ with a holomorphic symplectic form induces a $G$-equivariant isomorphism 
		$$\gs : H^1(S,T_S) \eto H^1(S,\gO^1_S).$$
		Let $\KS : T_{\gD,o} \eto H^1(S,T_S)$ be the ($G$-equivariant) Kodaira-Spencer map. 
		By~\eqref{eqn-KSgDV},
		we have 
		$$\gs(\KS(T_{\gD_{\NS(S)},o})) =  \NS(S)^\perp \subset H^1(S,\gO^1_S)$$
		where 
		$\NS(S)^\perp$ is defined with respect to the intersection pairing on $H^1(S,\gO^1_S)$. 
		Since $\gs \circ \KS$ is a $G$-equivariant isomorphism, 
		it suffices to show that $G$ acts as the identity on $\NS(S)^\perp$ to finish the proof. 
		When $S$ is a K3 surface, this is proven as follows.
		By~\cite[Lemma 3.3.1]{HuybK3book}, $\NS(S)^\perp$ is equal to the transcendental lattice $T(S) \subset H^2(S,\bZ)$.
		Since the $G$-action on $S$ is symplectic, $G \cto T(S) = \NS(S)^\perp$ is the identity by~\cite[Remark 15.1.2]{HuybK3book}.
		When $S$ is a 2-torus, 
		both results~\cite[Lemma 3.3.1 and Remark 15.1.2]{HuybK3book} continue to hold 
		since $H^2(S,\bZ)$ is also a Hodge structure of K3 type (namely $h^{2,0}(S) = 1$). 
		The same argument shows that $G \cto T(S) = \NS(S)^\perp$ is the identity.
\end{proof}

\begin{lem}\label{lem-2toreaa}
	Let $S$ be a $K$-trivial surface with $a(S) = 0$ and $G$ a finite group acting on $S$. Then the deformation $\Pi_{\NS(S)} : \cS_{\NS(S)} \to \gD_{\NS(S)}$ of $S$ preserving $\NS(S)$ in a semi-universal deformation $\Pi: \cS \to \gD$ of $S$ is a $G$-equivariant algebraic approximation of $S$ which is $Y$-locally trivial for every proper subvariety $Y \subsetneq S$. 	
\end{lem}

\begin{proof}
	
 By assumption, $S$ is either a K3 surface or a 2-torus.
Since $S$ is non-algebraic, we have $\NS(S) \otimes \bC \subsetneq H^{1,1}(S)$, so $\rk(\NS(S)) < h^{1,1}(S)$. It follows from Lemma~\ref{lem-aaK02} and Lemma~\ref{lem-lisseDV} that the family $\Pi_{\NS(S)} : \cS_{\NS(S)} \to \gD_{\NS(S)}$ is an algebraic approximation of $S$. 
By Lemma~\ref{lem-Gtriv}, the deformation $\Pi_{\NS(S)}$ is $G$-equivariant. It remains to show that $\Pi_{\NS(S)}$ is $Y$-locally trivial for every proper subvariety $Y \subsetneq S$.

Since $\Pi_{\NS(S)}$ is locally trivial by Lemma~\ref{lem-loctrivsm}.\emph{i)}, 
it suffices by Lemma~\ref{lem-loctrivsm}.\emph{ii)} to show that 
$\Pi_{\NS(S)}$ is $C$-locally trivial for every connected curve $C \subset S$. 
As $a(S) = 0$, by Lemma~\ref{lem-fini-2} we can assume that
 $C$ is a connected component of the (finite) union of all curves of $S$, 
 which is an A-D-E curve. 
  Since $\cS_{\NS(S)} \to \gD_{\NS(S)}$ preserves $\NS(S)$, 
  for each irreducible component $C_i$ of $C$, 
  the line bundle $\cO_S(C_i)$ deforms with $S$ and extends to a line bundle $\cL$ over $\cS_{\NS(S)}$.
  Since $h^i(S,\cO_S(C_i)) = 0$ for $i = 1,2$, 
  up to shrinking $\gD_{\NS(S)}$ we have $h^i(\cS_t,\cL_t) = 0$ for all $t \in \gD_{\NS(S)}$
  		where $\cS_t \cnec \Pi^{-1}(t)$ and $\cL_t \cnec \cL_{|\cS_t}$ by upper semicontinuity. 
  	So $h^0(\cS_t, \cL_t) = \chi(\cL_t)$, which is constant equal to $1$, up to further shrinking $\gD_{\NS(S)}$.
  		Hence each irreducible component $C_i$ of $C$ deforms with $S$ along $\gD_{\NS(S)}$ by Grauert's base change theorem.
 
  Let $U \subset S$ be a strictly pseudoconvex neighborhood of $C$~\cite[p. 357]{Grauert} and 
 let $\cU_{\NS(S)} \subset \cS_{\NS(S)}$ be an  open subset such that $\cU_{\NS(S)} \cap S = U$. 
By~\cite[Theorem 8]{LauferDefV}, the strictly pseudoconvex neighborhood 
$U$ admits a versal deformation $\Pi_{\mathrm{vers}} : \cU \to \gD_U$ over a smooth base $\gD_U$ and the
tangent space of $\gD_U$ at the point parameterizing $U$ is isomorphic to $H^1(U,T_U)$.
The versal property of $\Pi_{\mathrm{vers}} : \cU \to \gD_U$ asserted in~\cite[Theorem 8]{LauferDefV}
implies that up to shrinking $\gD_{\NS(S)}$, 
there exist open subsets $\cU'_{\NS(S)} \subset \cU_{\NS(S)}$ and $\cU' \subset \cU$ with $C \subset \cU'_{\NS(S)} 
$ such that
$\cU'_{\NS(S)} \xto{\Pi_{\NS(S)}} \gD_{\NS(S)}$ is the pullback of $\cU' \xto{\Pi_{\mathrm{vers}}} \gD_U$
by a holomorphic map $\varphi: \gD_{\NS(S)} \to \gD_U$.
 On the one hand, since $U$ is a neighborhood of the A-D-E curve $C$, 
 we have $\dim \gD_U = \dim H^1(U,T_U) = m$ where $m$ is the number of irreducible components of $C$~\cite[\S 1.8]{BurnsWahlloccontri}.
On the other hand, 
the image of $\varphi : \gD_{\NS(S)} \to \gD_U$ has codimension at least $m$
by~\cite[Theorem 3.7.i)]{LauferAmbient}. 
Thus $\varphi: \gD_{\NS(S)} \to \gD_U$ is constant, which implies that the family 
$\cU'_{\NS(S)} \to \gD_{\NS(S)}$ is trivial and shows that $\Pi_{\NS(S)}$ is $C$-locally trivial.
 \end{proof}

\begin{pro}\label{pro-aaK02}
	Let $G$ be a finite group acting on a K-trivial surface $S$ and on a smooth projective curve $B$. Assume that $a(S) = 0$. 
	Then $S \times B$, endowed with the diagonal $G$-action,
	has a $G$-equivariant algebraic approximation,
	which is $C$-locally trivial for every curve $C \subset S \times B$. 
\end{pro}

\begin{proof}
	
	By Lemma~\ref{lem-2toreaa}, there exists a 
	$G$-equivariant algebraic approximation $\cS \to \gD$ of $S$ 
	which is $Y$-locally trivial for every 
	proper subvariety $Y \subsetneq S$. 
	In particular, the product $\Pi :  \cS \times B \to \gD$ is 
	a $G$-equivariant algebraic approximation of $S \times B$.
	
	Let $C \subset S \times B$ be a curve and 
	let $\pr_1 : S \times B \to S$ be the first projection. 
	Since $\cS \to \gD$ is $\pr_1(C)$-locally trivial, there exists 
	an open subset $\cU \subset \cS$ containing $\pr_1(C)$ 
	such that $\cU \simeq U \times \gD$ over $\gD$ where $U \cnec \cU \cap S$.
	So the open subset $\cU \times B \subset  \cS \times B$ contains $C$ and is isomorphic to  $(U \times B) \times \gD$ over $\gD$, which shows that $\Pi$ is $C$-locally trivial. 
\end{proof}

For later use (in \S\ref{sec-a0K3}), let us prove a related result about deformations of compact K\"ahler surfaces of algebraic dimension 0.

\begin{lem}\label{lem-lifta=0frommin}
	Let $S$ be a smooth compact K\"ahler surface with $a(S) = 0$ and let 
	$$\eta : \tilde{S} = S_n \to \cdots \to S_1 \to S_0 = S$$
	be a sequence of blow-ups of $S_i$ at one point (\eg the contraction of $\ti{S}$ to its minimal model).
	Let $\Pi : \cS \to \gD$ be a deformation of $S$ such that for every subvariety $Y \subsetneq S$, up to shrinking $\gD$, $\Pi$ is $Y$-locally trivial.	
	Then there exists a sequence of blow-ups
	$$\ti{\eta} : \tilde{\cS} = \cS_n \to \cdots \to \cS_0 = \cS$$ 
	along sections of $\cS_i \to \gD$ such that the composition
	$$\ti{\Pi} : \tilde{\cS} \xto{\ti{\eta}} \cS \xto{\Pi} \gD$$ 
	is a deformation of $\eta$ and
	for every  subvariety $\ti{Y} \subsetneq \ti{S}$,
	 up to shrinking $\gD$, $\ti{\Pi}$ 
	is a $\ti{Y}$-locally trivial deformation of $\tilde{S}$.
\end{lem}

\begin{proof}

By induction, it suffices to prove 
Lemma~\ref{lem-lifta=0frommin} for $n = 1$, 
namely when $\eta : \ti{S} \to S$ is the blow-up at a point $p \in S$. By Lemma~\ref{lem-fini-2}, $S$ contains only finitely many curves. Let $Z_0 \subset S$ be the union of all curves of $S$ and let $Z \cnec Z_0 \cup \{p\}$.
By assumption, up to shrinking $\gD$,  
$\Pi : \cS \to \gD$ is $Z$-locally trivial, 
so by Lemma~\ref{lem-Yltexpl}
there exist 
 a neighborhood  $\cU \subset \cS$ of $Z$ and an isomorphism  
\begin{equation}\label{eqn-isomp}
\Phi :  \cU  \eto  U  \times \gD \ \text{ over } \ \gD
\end{equation}
with $U \cnec \cU \cap X$, such that the restriction of $\Phi$ to the central fiber $U$ is the identity. 
Let $\cP \cnec \Phi^{-1}(\{p\} \times \gD)$ and 
let $\ti{\eta} : \ti{\cS} \to \cS$ be the blow-up of $\cS$ along $\cP$. 
Then $\ti{\Pi} : \ti{\cS} \xto{\ti{\eta}} \cS \xto{\Pi} \gD $ is a deformation of $\eta : \ti{S} \to S$ which lifts $\Pi$, and $\ti{\Pi} : \ti{\cS} \to \gD$ is $\ti{Z}$-locally trivial by construction and~\eqref{eqn-isomp} where $\ti{Z} \cnec \eta^{-1}(Z)$. 
Let $\ti{Y} \subsetneq \ti{S}$ be a subvariety of $\ti{S}$ and  
we write $\ti{Y} = \ti{Y}' \sqcup \ti{Y}''$ where $\ti{Y}' \cnec \ti{Y} \cap \ti{Z}$ and $\ti{Y}'' \cnec \ti{Y} \bss \ti{Z}$.
Since $\ti{Y}' \subset \ti{Z}$ and  $\ti{\Pi}$ is $\ti{Z}$-locally trivial, $\ti{\Pi}$ is  $\ti{Y}'$-locally trivial. 
As $\dim \ti{Y}'' = 0$ because $\ti{Y}'' \cnec \ti{Y} \bss \ti{Z}$ 
and $\ti{Z}$ contains every curve of $\ti{S}$,  by Lemma~\ref{lem-loctrivsm}.\emph{i), ii)} $\ti{\Pi}$ is also $\ti{Y}''$-locally trivial. 
Since $\ti{Y} = \ti{Y}' \sqcup \ti{Y}''$, 
it follows from Lemma~\ref{lem-loctrivsm}.\emph{ii)} that
 $\ti{\Pi} $ is $\ti{Y}$-locally trivial. 
\end{proof}

\ssec{Proof of Propositions~\ref{pro-main3gk1} and~\ref{pro-mainpair}}
\hfill

We finish \S\ref{sec-casparcas} with the proof of Propositions~\ref{pro-main3gk1} and~\ref{pro-mainpair}.

\begin{proof}[Proof of Propositions~\ref{pro-main3gk1} and~\ref{pro-mainpair}]
First we prove Proposition~\ref{pro-mainpair}.
Let $X$ be a variety as in Proposition~\ref{pro-classk01}. If $X$ is the total space of a fibration $f : X \to B$ whose general fiber is an abelian surface, 
then Proposition~\ref{pro-mainpair} follows from Theorem~\ref{thm-AbFibDefprec} and Corollary~\ref{cor-sltaa}. Otherwise $X$ is the quotient $\ti{X}/G$ by a finite group $G$ of a threefold $\ti{X}$ described in Proposition~\ref{pro-classk01}, 
and it suffices by Lemma~\ref{lem-Gquotloctriv} and Lemma~\ref{lem-Gloctriv} to show that 
$\ti{X}$ has a $G$-equivariant algebraic approximation which is $Y$-locally trivial 
for every subvariety $Y \subset \ti{X}$ satisfying $\dim Y \le 1$. 
Since $\ti{X}$ is smooth, every small deformation of $\ti{X}$ is smooth, so by Lemma~\ref{lem-loctrivsm}.\emph{i),ii)}, 
it suffices to show that $\ti{X}$ has a $G$-equivariant algebraic approximation 
which is $C$-locally trivial for every curve $C \subset \ti{X}$.
If $X$ is in the first case of Proposition~\ref{pro-classk01}, 
then we use either Proposition~\ref{pro-aa3tore} or Proposition~\ref{pro-aaK02} to conclude. If $X$ is in the second case of Proposition~\ref{pro-classk01}, we conclude by Lemma~\ref{lem-isotriv2tore}. 

Finally, Proposition~\ref{pro-main3gk1} 
follows from Proposition~\ref{pro-classk01} proven in \S\ref{sec-bim}
and Proposition~\ref{pro-mainpair}.
\end{proof}

\section{Proof of Proposition~\ref{pro-sec} for $\bP^1$-bundles}\label{sec-a0K3}

In \S\ref{sec-a0K3}, we will study the existence of algebraic approximations for
$\bP^1$-bundles over a compact K\"ahler surface of algebraic dimension 0.
The main result that we will prove is Proposition~\ref{pro-aaa=0}.
Since a $\bP^1$-bundle is the projectivization of a twisted vector bundle of rank 2, their deformation problem is directly related to the deformation-obstruction theory of twisted vector bundles that we shall study now.

\ssec{Semi-regularity maps and deformations of twisted sheaves}\label{ssec-deftrod}\hfill

We refer to \S\ref{ssec-twSAtiyah} for a reminder of twisted coherent sheaves and Atiyah class.

Let $\Pi : \cX \to \gD$ be a smooth deformation of 
a compact K\"ahler manifold $X$ over a contractible manifold $\gD$.
Let $\{U_i\}_{i \in I}$ be a sufficiently fine good open cover of $X$. 
By Ehresmann's fibration theorem~\cite[Theorem 9.3]{VoisinI}, 
we have a diffeomorphism $\Psi : X \times \gD \to \cX$ over $\gD$. 
Set $\Set{\cU_i \cnec \Psi(U_i \times \gD) }_{i \in I}$, which forms a good open cover of $\cX$. 
Let $r \in \bZ_{> 0}$ and let
$\ga$ be a \v{C}ech 2-cocycle with respect to $\{U_i\}_{i \in I}$ with coefficients in  the locally constant subsheaf $\mu_r \subset \cO_X^\times$ of $r$th roots of unity. 
The 2-cocycle $\ga$ extends uniquely to a 2-cocycle $\ti{\ga}$ with coefficients in $\mu_r \subset \cO_{\cX}^\times$  with respect to $\{\cU_i \}_{i \in I}$. 

Let $o \in \gD$ be the point parameterizing $X$ and let $\fm \subset \cO_{\gD,o}$ be its ideal sheaf. For every $n \in \bZ_{\ge 0}$, let $\gD_n = \Spec(\cO_{\gD,o}/\fm^{n+1})$ be the $n$-th infinitesimal neighborhood of $o$ and let $\cX_{n} \to  \gD_{n}$ be the base change of $\Pi : \cX \to \gD$ by $\gD_n \to \gD$. Let $\ga_n \cnec \ti{\ga}_{|\cX_n}$. The following lifting criterion involving the injectivity of semi-regularity maps is classical in the untwisted case (see for example~\cite[Theorem 5.1]{BuchweitzFlennerACsemireg}) and can be proven with a similar argument. 

\begin{pro}\label{pro-semireg}  
 Let $E$ be an $\ga$-twisted locally free sheaf on $X$.
Assume that the trace map $\tr : \Ext^2(E,E) \to H^2(X,\cO_X)$ is injective 
and the trace of the Atiyah class (see \S\ref{ssec-twSAtiyah})
$$\Tr\(\At(E)\) \in H^1(X,\gO^1_X) \subset H^2(X,\bC)$$ 
remains of type $(1,1)$ along $\gD$ under the parallel transport with respect to the Gauss-Manin connection. 
Then starting from $\cE_0 = E$, for each $n \in \bZ_{\ge 0}$ we can lift $\cE_n$ inductively to an $\ga_{n+1}$-twisted locally free sheaf $\cE_{n+1}$ on $\cX_{n+1}$.
\end{pro}

\begin{proof}

Assume that $\cE_n$ is constructed. 
Let $\gt \in \Ext^1(\gO_{\cX_{n}/\gD_n},I)$ be the element	corresponding  to the square-zero extension  $\cX_{n+1} \to \gD_{n+1}$  of $\cX_n \to \gD_n$ where $I$ is the ideal sheaf of $\cX_n$ in $\cX_{n+1}$  considered as an $\cO_{\cX_n}$-sheaf. In other words, $\gt$ is the class representing the conormal exact sequence
$$
\begin{tikzcd}[cramped, row sep = 0, column sep = 30]
0  \arrow[r] & I  \arrow[r]  & (\gO_{\cX_{n+1}/\gD_{n+1}})_{|\cX_n} \arrow[r]  & \gO_{\cX_{n}/\gD_n}\arrow[r]  &  0.
\end{tikzcd}
$$

\begin{lem}\label{lem-lift}
Define
$$
\begin{tikzcd}[cramped, row sep = 0, column sep = 30]
\ob_{\cE_n} : \Ext^1(\gO_{\cX_{n}/\gD_n},I)  \arrow[r,  "\Id_{\cE_n} \otimes"] &  \Ext^1( \cE_n \otimes\gO_{\cX_{n}/\gD_n},\cE_n \otimes I)  \arrow[r]  & \Ext^2(\cE_n,\cE_n  \otimes I)
\end{tikzcd}
$$
where the last arrow is the Yoneda product with $\At(\cE_n/\gD_n) \in \Ext^1(\cE_n, \cE_n \otimes\gO_{\cX_{n}/\gD_n}) $.
Then $\cE_n$ can be extended to an $\ga_{n+1}$-twisted locally free sheaf $\cE_{n+1}$ on $\cX_{n+1}$ if and only if $\ob_{\cE_n}(\gt) = 0$.
\end{lem}

\begin{proof}
	Recall that $\At(\cE_n/\gD_n) \in \Ext^1(\cE_n, \cE_n \otimes\gO_{\cX_{n}/\gD_n})$ 
	is the extension class of the first jet bundle $J^1(\cE_n)$  (see \S\ref{ssec-twSAtiyah}).
	Applying~\cite[Exercise III.5.3]{GelfandManinHomalg} to the abelian category $\Coh(\cX_n,\ga_n)$ and the extension classes $\At(\cE_n/\gD_n) \in \Ext^1(\cE_n, \cE_n \otimes\gO_{\cX_{n}/\gD_n}) $ and $\Id_{\cE_n} \otimes \gt \in \Ext^1( \cE_n \otimes\gO_{\cX_{n}/\gD_n},\cE_n \otimes I)$, the vanishing $\ob_{\cE_n}(\gt) = 0$ is equivalent to the existence of an $\ga_n$-twisted sheaf $\cF$ sitting in the middle of the following commutative diagram with exact rows and columns.
\begin{equation}\label{diag-ext2}
\begin{tikzcd}[cramped, row sep = 15, column sep = 15]
& 0 \ar[d] & 0 \ar[d] \\
& \cE_n\otimes I  \ar[d] \ar[r, equal] & \cE_n\otimes I \ar[d] \\
0 \ar[r] & \cE_n \otimes{ \gO_{\cX_{n+1}/\gD_{n+1}}}_{|\cX_n} \arrow[r] \ar[d] & \cF \arrow[r] \ar[d, "r"] & \cE_n  \arrow[r] \ar[d,equal] & 0 \\
0 \ar[r] & \cE_n \otimes{ \gO_{\cX_{n}/\gD_n}} \ar[d] \arrow[r] &J^1(\cE_n)\ar[d]  \arrow[r] & \cE_n  \arrow[r] & 0 \\
& 0  & 0  
\end{tikzcd}
\end{equation}
If $\cE_n$ can be extended to an $\ga_{n+1}$-twisted locally free sheaf $\cE_{n+1}$ on $\cX_{n+1}$, then $\cF \cnec J^1(\cE_{n+1})_{|\cX_n}$ fits in~\eqref{diag-ext2}, 
where the horizontal exact sequence in the middle is the restriction to $\cX_n$ of the analogue extension~\eqref{exseq-extJet} for $J^1(\cE_{n+1})$ and $r : J^1(\cE_{n+1})_{|\cX_n} \to J^1(\cE_n)$ is the map induced by the restriction to $\cX_n$.

Conversely, assume that~\eqref{diag-ext2} exists. Recall that we have an open cover $\{\cU\}_{i \in I}$ of $\cX$ and we 
define $\Set{U_i^{(n)} \cnec \cU_i \cap \cX_n}_{i \in I}$ as an open cover of $\cX_n$. Let $\cF_i \cnec \cF_{|U_i^{(n)}}$ and $\cE_{n,i} \cnec {\cE_n}_{|U_i^{(n)}}$, which are untwisted sheaves. 
First we recall how the extension $\cE_{n+1,i}$ of $\cE_{n,i}$ is constructed as an (untwisted) $\cO_{U_i^{(n+1)}}$-sheaf~\cite[Proposition 4.4]{BuchweitzFlennerACsemireg}. 
Let $j : \cE_n \to J^1(\cE_n)$ be the first-jet map (see \S\ref{ssec-twSAtiyah}). As a sheaf of $\bC$-vector spaces over $U_i^{(n)}$, define 
\begin{equation}\label{def-En}
\cE_{n+1,i} \colonec \ker\(r\circ \pr_1 + j \circ \pr_2 : \cF_i \oplus \cE_{n,i} \to J^1(\cE_{n,i})\)
\end{equation}
We define the $\cO_{U_i^{(n+1)}}$-module structure on $\cE_{n+1,i} $ as follows. First we identify the underlying sheaf of $\bC$-vector spaces of $\cO_{\cX_{n+1}}$ with
\begin{equation}\label{def-On}
\cO_{\cX_{n+1}} = \ker\(r'\circ \pr_1 - d \circ \pr_2 : { \gO_{\cX_{n+1}/\gD_{n+1}}}_{|\cX_n} \oplus \cO_{\cX_n} \to  \gO_{\cX_n/\gD_n} \)
\end{equation}
where $r' :  {\gO_{\cX_{n+1}/\gD_{n+1}}}_{|\cX_n} \to \gO_{\cX_n/\gD_n}$ is the restriction map. Under this identification, the ring structure on $\cO_{\cX_{n+1}}$ is locally defined by 
\begin{equation}\label{eqn-LeibO}
(\gb,f)(\gb',f') = (f\gb' + f'\gb,ff')
\end{equation}
where $(\gb,f)$ and $(\gb',f')$ are local sections of $\cO_{\cX_{n+1}} \subset {\gO_{\cX_{n+1}/\gD_{n+1}}}_{|\cX_n} \oplus \cO_{\cX_n}$. Given local sections $(\gb, f)$ of $\cO_{\cX_{n+1}} \subset {\gO_{\cX_{n+1}/\gD_{n+1}}}_{|\cX_n} \oplus \cO_{\cX_n}$ and $(\gamma, \gs)$ of $  \cE_{n+1,i} \subset \cF_i \oplus  \cE_{n,i}$ over an open subset $U \subset U_i^{(n)}$, 
we define
\begin{equation}\label{eqn-onstru}
(\gb, f)\cdot (\gamma, \gs) \colonec (f \gamma + \gs \otimes \gb, f \gs) \in \cF_i(U) \oplus \cE_{n,i}(U)
\end{equation}
where $\gs \otimes \gb$ is regarded as an element of $\cF_i(U)$ using~\eqref{diag-ext2}.
It follows from the Leibniz rules~\eqref{eqn-Leibniz} and~\eqref{eqn-LeibO} together with~\eqref{diag-ext2},~\eqref{def-En}, and~\eqref{def-On} that $(\gb, f)\cdot (\gamma, \gs) \in \cE_{n+1,i}(U)$. We verify easily that~\eqref{eqn-onstru} defines an $\cO_{U_i^{(n+1)}}$-module structure on $\cE_{n+1,i}$. 

To show that $\cE_{n+1,i}$ is locally free, let $U \subset U_i^{(n)}$ be  a sufficiently small open subset over which $\cE_{n,i}$ is $\cO_{U_i^{(n)}}$-free and let $\gs_1,\ldots,\gs_r$ be a basis of  $\cE_{n,i}(U)$. Let $\gamma_1,\ldots,\gamma_r \in \cF(U)$ such that $(\gamma_i,\gs_i) \in \cE_{n+1,i}(U)$ under the identification~\eqref{def-En}. 
For every $(\gamma,\gs) \in \cE_{n+1,i}(U)$, there exist unique $f_1,\ldots,f_r \in \cO_{\cX_n}(U)$ 
such that $\sum_i f_i\gs_i = \gs$. 
As 
\begin{equation}
	r(\gamma) = -j(\gs) = -\sum_i j(f_i\gs_i)  = \sum_i (-f_i\cdot j(\gs_i) + \gs_i \otimes df_i) = \sum_i f_i \cdot r(\gamma_i) + \sum_i \gs_i \otimes df_i,
\end{equation}
it follows from~\eqref{diag-ext2} that $\gamma - \sum_i f_i\gamma_i \in  \cE_n \otimes {\gO_{\cX_{n+1}/\gD_{n+1}}}_{|\cX_n}(U)$. So there exist unique  $\gb_1,\ldots,\gb_r \in { \gO_{\cX_{n+1}/\gD_{n+1}}}_{|\cX_n}(U)$ such that 
$(\sum_i \gs_i \otimes \gb_i, 0) = (\gamma,\gs) - \sum_i (f_i\gamma_i, f_i\gs_i)$. 
Hence there exist unique $(\gb_1, f_1),\ldots,(\gb_r, f_r) \in \cO_{\cX_{n+1}}(U)$ such that
$$\sum_i (\gb_i, f_i)\cdot (\gamma_i, \gs_i)= (\gamma,\gs),$$
which shows that $\cE_{n+1,i}$ is locally free over $\cO_{U_i^{(n+1)}}$.

Finally we glue the local pieces $\cE_{n+1,i}$ into an $\ga_{n+1}$-twisted sheaf. Let 
$$ \gt_{ij} : {\cF_i}_{|U_i^{(n)} \cap U_j^{(n)}} \eto {\cF_j}_{|U_i^{(n)} \cap U_j^{(n)}} \ \text{ and } \  \gt'_{ij} : {\cE_{n,i}}_{|U_i^{(n)} \cap U_j^{(n)}} \eto {\cE_{n,j}}_{|U_i^{(n)} \cap U_j^{(n)}}$$ 
be the gluing isomorphisms defining the $\ga_n$-twisted sheaves $\cF$ and $\cE_n$. These morphisms are compatible with the morphisms in~\eqref{diag-ext2} and both $\gt_{ki} \circ \gt_{jk} \circ \gt_{ij}$ and $\gt'_{ki} \circ \gt'_{jk} \circ \gt'_{ij}$ are multiplications by $\ga_{ijk} \in \mu_r$.   
Under the identification~\eqref{def-En}, let
$$\gT_{ij} : {\cE_{n+1,i}}_{|U_i^{(n+1)} \cap U_j^{(n+1)}} \to  {\cE_{n+1,j}}_{|U_i^{(n+1)} \cap U_j^{(n+1)}}$$ 
be the restriction of $\gt_{ij} \oplus \gt'_{ij}$ to $ {\cE_{n+1,i}}_{|U_i^{(n+1)} \cap U_j^{(n+1)}} $. 
Using~\eqref{eqn-onstru} and the compatibility of the 
gluing isomorphisms $\gt_{ij}$, $\gt'_{ij}$ with the morphisms in~\eqref{diag-ext2},
we verify that $\gT_{ij}$ is $\cO_{\cX_{n+1}}$-linear. By construction, the composition 
$$\gT_{ki} \circ \gT_{jk} \circ \gT_{ij} : {\cE_{n+1,i}}_{|U_i^{(n+1)} \cap U_j^{(n+1)} \cap U_k} \to {\cE_{n+1,i}}_{|U_i^{(n+1)} \cap U_j^{(n+1)} \cap U_k^{(n+1)}}$$ 
is the multiplication by $\ga_{ijk} \in \mu_r$. Hence $\{\gT_{ij}\}$ is a collection of gluing isomorphisms which defines ${\cE_{n+1}}$ as an $\ga_{n+1}$-twisted sheaf on $\cX_{n+1}$.
\end{proof}

To finish the proof, it suffices by Lemma~\ref{lem-lift} to verify that $\ob_{\cE_n} (\gt) =0$.
We can compose $\ob_{\cE_n}$ with the trace map $\tr :  \Ext^2(\cE_n,\cE_n  \otimes I) \to H^2(\cX_n,I)$. As the trace maps commute with the Yoneda products, we have
$$
\tr \circ \ob_{\cE_n} (\gt) = \gt \cupp \Tr(\At(\cE_n/\gD_n)) \in H^2(\cX_n, I).
$$
Since the parallel transports of $\Tr(\At(E))$ by the Gauss-Manin connection remain of type $(1,1)$ by assumption,  we have $ \gt \cupp \Tr(\At(\cE_n/\gD_n)) =0$  by~\cite[Lemma 5.7 and 5.8]{BuchweitzFlennerACsemireg}. 
Since $\tr : \Ext^2(E,E) \to H^2(X,\cO_X)$ is injective, by~\cite[Lemma 5.10]{BuchweitzFlennerACsemireg} with the same proof adapted for twisted sheaves, the trace map $\tr :  \Ext^2(\cE_n,\cE_n  \otimes I) \to H^2(\cX_n,I)$ is also injective. Hence $\ob_{\cE_n} (\gt) =0$.
\end{proof}

\begin{cor}\label{cor-semireg}
	In Proposition~\ref{pro-semireg}, up to shrinking $\gD$, there exists a $\bP^{r-1}$-bundle $\cP \to \cX$ (with $r = \rank (E)$) whose restriction to $X \subset \cX$ coincides with $\bP(E) \to X$. If $E$ is untwisted, then 
	 the deformation $\cX \to \gD$ lifts to a deformation of the pair  $(X,E)$. 
\end{cor}
\begin{proof}
	The second statement is well-known and we refer to~\cite[p.179, Proof of Theorem 5.1]{BuchweitzFlennerACsemireg} for a proof. 
	
	We prove the first statement; the argument is inspired by the proof of~\cite[Theorem 8.1]{HorikawaIII}.
	Let $\mathbb{P} \to \Defo_{\bP(E)}$  be a versal deformation of the compact complex manifold $\bP(E)$. 
	Let $\cD\to \Defo_{\bP(E)} \times \gD $ be the Douady--Pourcin space~\cite[Th\'eor\`eme 2]{Pourcin} of closed complex subspaces of 
	$\mathbb{P} \times \cX$ over $\Defo_{\bP(E)} \times \gD$:
	this is the complex space over $\Defo_{\bP(E)} \times \gD$ which represents the functor
	\begin{equation*}
		\begin{split}
			\text{Complex spaces} / \Defo_{\bP(E)} \times \gD  & \to  \text{Sets} \\
			(S \to \Defo_{\bP(E)} \times \gD) & \mapsto 
			\Set{ \text{Closed complex subspaces } \cZ \subset (\mathbb{P} \times \cX) \times_{\Defo_{\bP(E)} \times \gD} S \text{ flat over } S}. 
		\end{split}
	\end{equation*}
Let $(\pi_1,\pi_2) : \cD \to \Defo_{\bP(E)} \times \gD$ denote the structural morphism.

	Let $\cE_n$ be the $\ga_n$-twisted sheaf on $\cX_n$ constructed  
	in Proposition~\ref{pro-semireg}. 
	Then $\bP(\cE_n) \to \cX_n \to \gD_n$ defines a deformation of $\bP(E) \to X$ over $\gD_n$;
	let $r_n : \gD_n \to \Defo_{\bP(E)}$ be a map which induces the deformation $\bP(\cE_n) \to \gD_n$.
	Let $\imath_n : \gD_n \hto \gD$ be the inclusion.
	Since the graph of $\bP(\cE_n) \to \cX_n$, 
	viewed as a closed complex subspace of $\bP(\cE_n) \times_{\gD_n} \cX_n$, 
	is flat over $\gD_n$,
	by the universal property of $\cD$ it defines a factorization
	$$(r_n, \imath_n) : \gD_n \xto{\gs_n} \cD \to \Defo_{\bP(E)} \times \gD$$
	of $(r_n, \imath_n)$ with $\gs_n(o) \in \cD$ being the point parameterizing the graph of $\bP(E) \to X$
	for all $n \in \bZ_{\ge 0}$. 
	Composing $(r_n, \imath_n)$ with the projection $\Defo_{\bP(E)} \times \gD \to \gD$,
	we obtain a factorization
	$$\imath_n : \gD_n \xto{\gs_n} \cD \xto{\pi_2} \gD$$
	of $\imath_n : \gD_n \hto \gD$.
	Applying Wavrik's generalization of
	Artin's approximation theorem~\cite[Theorem 2.4]{Wavrik} to the collection of commutative diagrams
	$$
	\begin{tikzcd}[cramped, row sep = 25, column sep = 25]
		\gD    \ar[dr, swap,"\Id"] &  \gD_n \ar[l,hook'] \ar[r,"\gs_n"] & \cD \ar[dl, "\pi_2"]   \\
		& \gD  &   
	\end{tikzcd}
	$$ 
	for each $n \in \bZ_{\ge 0}$,
	we obtain a section $\gs : \gD \to \cD$ of $\pi_2 : \cD \to \gD$ such that $\gs(o) \in \cD$ is the point parameterizing the graph of $\bP(E) \to X$.
	By the universal property of $\cD$, the section $\gs : \gD \to \cD$ defines a deformation 
	$$\cZ \subset \cP \times_\gD \cX \to \gD$$
	of the graph $\Gamma \subset \bP(E) \times X$ of $\bP(E) \to X$ over $\gD$, 
	where $\cP \to \gD$ is the deformation of $\bP(E)$ induced by $ \pi_1 \circ \gs :  \gD \to \Defo_{\bP(E)}$.
	Since the projection $\Gamma \to \bP(E)$ is an isomorphism, so is $\cZ \to \cP$ up to shrinking $\gD$~\cite[Proposition 10.1]{EspDouady}. 
	Thus $\cZ \subset \cP \times_\gD \cX$ 
	is the graph of a morphism $\cP \to \cX \to \gD$ over $\gD$, which is a deformation of $\bP(E) \to X \to \gD$ lifting $\cX \to \gD$.
	Since small deformations of $\bP^{r-1}$ remain $\bP^{r-1}$ and since $\cX \to \gD$ is proper, up to shrinking $\gD$ it follows that $\cP \to \cX$ is a $\bP^{r-1}$-bundle.
\end{proof}

\begin{rem}
	The first statement of Corollary~\ref{cor-semireg} would follow more naturally from the existence of a deformation of the pair $(X,E)$ lifting the deformation of $X$.
	This will follow from Proposition~\ref{pro-semireg} if we know that 
	versal deformations of twisted sheaves exist, and we expect that it is the case.
\end{rem}

\ssec{Deformations of $\bP^1$-bundles}\label{ssec-P1defs}
\hfill

In \S\ref{ssec-P1defs}, we prove the existence of algebraic approximations of 
$\bP^1$-bundles over surfaces  $\ti{S}$ with $a(\tilde{S}) = 0$ 
(see Proposition~\ref{pro-aaa=0} for the precise statement).
Let us start with the following lifting result.

\begin{pro}\label{pro-liftP1}
Let $\tilde{S}$ be a (smooth) compact K\"ahler surface with $a(\tilde{S}) = 0$ and let $f: X\to \tilde{S}$ be a $\bP^1$-bundle over $\tilde{S}$.
Let $\eta : \ti{S} \to S$ be the map from $\ti{S}$ to its minimal model 
(which is either a K3 surface or a 2-torus).
 Let $\Pi : \cS \to \gD$ be a deformation of $S$ preserving $\NS(S)$ and 
 assume that $\Pi$ is $\gS$-locally trivial for every subvariety $\gS \subsetneq S$. 
 Let 
$$\ti{\Pi} : \ti{\cS} \xto{\ti{\eta}} \cS \to \gD$$ 
be the deformation of $\ti{S}$ constructed in Lemma~\ref{lem-lifta=0frommin} which lifts $\Pi$. 
Then up to shrinking $\gD$, $\ti{\Pi}$ can be lifted to a deformation 
$$\cX \xto{\ti{f}} \ti{\cS} \xto{\ti{\Pi}} \gD$$
of the $\bP^1$-bundle $f : X \to \ti{S}$ with $\ti{f}$ being a $\bP^1$-bundle as well.
\end{pro}

\begin{proof}
 
As we mentioned in \S\ref{ssec-twSAtiyah}, the $\bP^1$-bundle $f$
is isomorphic to the projectivization $\bP(\ti{E}) \to \tilde{S}$
of an $\ga$-twisted locally free sheaf $\ti{E}$ of rank 2
where $\ga$ is a \v{C}ech 2-cocycle 
with coefficients in $\mu_2 \subset \cO_{\ti{S}}^\times$.

First we assume that $\ti{E}$ contains an $\ga$-twisted subsheaf $\ti{L}_0$ of rank 1. 
Up to replacing $\ti{L}_0$ with $\ti{L}_0^{\vee \vee}$,
we can assume that $\ti{L}_0$ is reflexive,
so $\ti{L}_0$ is invertible because $\ti{S}$ is a smooth surface. 
Since $\bP(\ti{E}) \simeq \bP(\ti{E} \otimes \ti{L}_0^\vee)$ over $\ti{S}$, up to replacing $\ti{E}$ by $\ti{E} \otimes \ti{L}_0^\vee $ we can assume that $\ti{E}$ is an untwisted locally free sheaf containing $\cO_{\ti{S}}$ (see also Lemma~\ref{lem-untw}).
Up to further tensoring $\ti{E}$ with another (untwisted) invertible sheaf,
we can assume that $\cO_{\ti{S}} \subset \ti{E}$ is saturated. 
To prove Proposition~\ref{pro-liftP1} in this case, it suffices to show that up to shrinking $\gD$, $\ti{\Pi}$ lifts to a deformation of the pair $(\ti{S}, \ti{E})$. By~\cite[Lemma 4.7]{Schrackdefo}, it suffices to show that the locally free sheaf $E \cnec (\eta_*\ti{E})^{\vee \vee}$ deforms with ${S}$ over $\gD$\footnote{
	In~\cite[Lemma 4.7]{Schrackdefo}, 
	we can replace "algebraic approximation" by "deformation"
	everywhere in the statement without changing the proof. 
	Here we apply~\cite[Lemma 4.7]{Schrackdefo} through this statement.}. 

As we assume that $\cO_{\ti{S}} \subset \ti{E}$ is saturated, 
we have $\cO_S \subset E$ as a saturated subsheaf. 
So we have an extension
\begin{equation}\label{suitex-extE}
\begin{tikzcd}[cramped, row sep = 0, column sep = 30]
0  \arrow[r] & \cO_S  \arrow[r]  & E \arrow[r]  & L \otimes I_Y \arrow[r]  &  0
\end{tikzcd}
\end{equation}
where $L$ is a line bundle on $S$ and $I_Y$ is the ideal sheaf of a 0-dimensional 
complex subspace $Y \subset S$.
Since $\Pi : \cS \to \gD$ preserves $\NS(S)$, 
 every line bundle on $S$ deforms with ${S}$ over $\gD$.
 In particular, there exists an invertible sheaf $\cL$ on $\cS$ which is a deformation of $L$. 
 By~\cite[Proposition 4.3 and 4.4]{Schrackdefo} 
 we can choose $\cL$ such that up to shrinking $\gD$, 
 $h^i(\cS_t,\cL)$ is constant in $t \in \gD$ for every $i$ where $\cS_t \cnec \Pi^{-1}(t)$.
 Let $Z$ be the union of all the curves of $S$, which has only finitely many irreducible components by Lemma~\ref{lem-fini-2}.
 Since $\Pi$ is $(Y \cup Z)$-locally trivial by assumption, 
 there exists by Lemma~\ref{lem-Yltexpl}
 a neighborhood  $\cU \subset \cS$ of $Y \cup Z$ together with an isomorphism  
 \begin{equation}\label{iso-paireYC}
\Phi : \cU \eto U \times \gD \ \text{ over } \ \gD,
 \end{equation}
with $U \cnec \cU \cap X$ such that the restriction of~\eqref{iso-paireYC} to the central fiber $U$ is the identity. 
We fix such an isomorphism $\Phi$ and 
define  
  $$\cY \cnec \Phi^{-1}(Y \times \gD) \subset \cS, \ \ \ \ \ \cZ \cnec \Phi^{-1}(Z \times \gD) \subset \cS,$$ namely 
 the trivial deformations of $Y$ and $Z$ in $\cS \to \gD$ induced by~\eqref{iso-paireYC}.

Now to show that $E$ deforms with $S$ over $\gD$ up to shrinking $\gD$, 
it suffices to prove the following.
\begin{claim}\label{claim-constext1}
	Up to shrinking $\gD$,  
	$$t \mapsto \dim \Ext^1_{\cO_{\cS_t}}((\cL \otimes I_{\cY})_{|\cS_t}, \cO_{\cS_t}) = h^1\({\cS}_t,\cL \otimes I_\cY\)$$ 
	is constant in $t \in \gD$ 
	where $I_\cY$ is the ideal sheaf of $\cY \subset \cS$.
\end{claim}  
Indeed, as both $\cL \otimes I_{\cY}$ and $\cO_{\cS}$ are flat over $\gD$,
Claim~\ref{claim-constext1} together with~\cite[Satz 3]{GrauertExt} 
implies that the relative ext-sheaf
$\cExt^1_{\Pi}(\cL \otimes I_{\cY},\cO_{\cS})$ is locally free and $\cExt^1_{\Pi}(\cL \otimes I_{\cY},\cO_{\cS})|_{o} \simeq  \Ext^1_{\cO_S}(L \otimes I_{Y},\cO_{S})$ where $o \in \gD$ is the point parameterizing $S$. 
Up to shrinking $\gD$, we can find a section $\gs$ of $\cExt^1_{\Pi}(\cL \otimes I_{\cY},\cO_{\cS})$ such that $\gs|_o \in \Ext^1_{\cO_S}(L \otimes I_{Y},\cO_{S})$ corresponds to the extension~\eqref{suitex-extE}. As $\gD$ is Stein up to shrinking $\gD$, Leray spectral sequence and the vanishing of higher cohomology groups show that
$$\gs \in H^0(\gD, \cExt^1_{\Pi}(\cL \otimes I_{\cY},\cO_{\cS})) 
\simeq \Ext_{\cO_{\cS}}^1(\cL \otimes I_{\cY},\cO_{\cS}),$$
and the sheaf $\cE$ on $\cS$ defined by $\gs$ as an extension of $\cL \otimes I_{\cY}$ by $\cO_{\cS}$ defines a deformation of the pair $(S,E)$ over $\gD$.

\begin{proof}[Proof of Claim~\ref{claim-constext1}]
Note that the equality in the statement follows from Serre duality because $\cS_t$ is $K$-trivial.
Since $\cY_t \cnec \cY \cap \cS_t$ has dimension $0$ and flat over $\gD$, 
$h^0(\cY_t,\cL)$ is constant in $t$ and $h^1(\cY_t,\cL) = 0$.
By assumption,
$h^i(\cS_t,\cL)$ is also constant in $t$,
so by the exact sequence
$$
\begin{tikzcd}[cramped, row sep = 0, column sep = 20]
0 \arrow[r] & H^0( {\cS}_t,  {\cL} \otimes I_{\cY} )\arrow[r] & H^0( {\cS}_t, \cL ) \arrow[r] &H^0( {\cY}_t, \cL )  \\
\arrow[r] & H^1( {\cS}_t,  {\cL} \otimes I_{\cY} )\arrow[r] & H^1( {\cS}_t, \cL ) \arrow[r] &H^1( {\cY}_t, \cL ) =  0,
\end{tikzcd}
$$
 it suffices to show that  $t \mapsto h^0( {\cS}_t,  {\cL} \otimes I_{\cY} )$  is constant.

If $h^0\(S, {\cL}  \otimes I_{\cY} \) = 0$, then up to shrinking $\gD$,
$t \mapsto h^0(\cS_t,\cL \otimes I_{\cY}) $
is constant by upper semi-continuity. Assume  that $h^0\(S, {\cL}  \otimes I_{\cY} \)  > 0$. 
Since $a(S) = 0$, by~\cite[Proposition IV.8.1]{Barth} $1 \le h^0\(S, {\cL}  \otimes I_{\cY } \) \le h^0\(S, {\cL} \) \le 1$. 
So $\cL_{|S} \simeq \cO_S(D)$ for some divisor $D \subset S$ containing $S \cap \cY = Y$.
Let 
$$\cD \cnec \Phi^{-1}(D \times \gD) 
\subset \cU \subset \cS,$$
which is a divisor on $\cS$ flat over $\gD$.
As $\cO_\cS(\cD)_{|S} = \cO_S(D) \simeq \cL_{|S}$,
the vanishing $H^1(\cS_t,\cO_{\cS_t}) = 0$ implies that 
$\cO_\cS(\cD)_{|\cS_t}  \simeq \cL_{|\cS_t}$ for all $t \in \gD$.
It follows that $\cO_\cS(\cD) \simeq \cL \otimes \Pi^*\cM$ for some line bundle $\cM$ on $\gD$,
so up to shrinking $\gD$, we have $\cO_\cS(\cD)  \simeq \cL$.
Let $\gs \in H^0(\cS,\cO_\cS(\cD))$ such that $\Div(\gs) = \cD$. 
Since $Y \subset D$, we have
$$\cY = \Phi^{-1}(Y \times \gD) \subset \Phi^{-1}(D \times \gD) = \cD = \Div(\gs).$$
Thus $\gs_{|\cS_t} \ne 0$ and $\Div(\gs_{|\cS_t}) \supset (\cY \cap \cS_t)$, 
showing that $h^0( {\cS}_t,  {\cL} \otimes I_{\cY} ) =  h^0(\cS_t,\cO_\cS(\cD) \otimes I_\cY)  \ge 1$ for all $t \in \gD$.
Finally since $1 = h^0\(S, {\cL}  \otimes I_{\cY } \)$, up to further shrinking $\gD$
it follows that 
$h^0( {\cS}_t,  {\cL} \otimes I_{\cY} )$ is constant in $t \in \gD$ by upper semi-continuity.
\end{proof}

Now we prove Proposition~\ref{pro-liftP1} for the case 
where $\ti{E}$ does not contain any $\ga$-twisted subsheaf of rank 1.
As $\ga$ is a $2$-cocycle with coefficients in $\mu_2 \subset \cO_{\ti{S}}^\times$,  the square $\ti{E}^{\otimes 2}$ is isomorphic to an untwisted locally free sheaf. So by~\eqref{eqn-tensAt} and~\eqref{eqn-Atc1}
$$4\cdot \Tr(\At(\ti{E})) = \Tr(\At(\ti{E}^{\otimes 2})) = - c_1(\ti{E}^{\otimes 2}) \in \NS(\ti{S}).$$  
By assumption, $\ti{\Pi} : \ti{\cS} \to \gD$ preserves $\eta^*\NS(S)$ as well as every $(-1)$-curve in $\ti{S}$, so $\ti{\Pi}$ also preserves  $\NS(\ti{S})$. 
In particular, $\Tr(\At(\ti{E}))$ remains of type $(1,1)$ under the parallel transport by the Gauss-Manin connection along $\gD$. Thus by Corollary~\ref{cor-semireg}, it suffices to show that the trace map $\tr : \Ext^2(\ti{E},\ti{E}) \to H^2(\tilde{S},\cO_{\tilde{S}})$ is injective. 

We show that the dual  $\tr^\vee$ of the trace map is an isomorphism.
By Serre duality, $\tr^\vee$  is isomorphic to 
$$H^0(\tilde{S},K_{\tilde{S}}) \to H^0(\ti{S}, \ti{E}^\vee \otimes \ti{E} \otimes K_{\tilde{S}}) = \Hom(\ti{E},\ti{E} \otimes K_{\tilde{S}}) $$ 
induced by $\cO_{\ti{S}} \hto \ti{E}^\vee \otimes \ti{E}$. So $\tr^\vee$ is injective. It remains to show that $h^0(\tilde{S},K_{\tilde{S}}) \ge  \dim\Hom(\ti{E},\ti{E} \otimes K_{\tilde{S}})$. On the one hand since $\tilde{S}$ is a sequence of blow-ups of a K-trivial surface, we have $h^0(\ti{S}, \go^{\otimes i}_{\ti{S}}) = 1$ for every $i \ge 0$. On the other hand since $\rk(\ti{E}) = 2$ and $\ti{E}$ contains no subsheaf of rank 1, for every $\phi \in \Hom(\ti{E},\ti{E} \otimes K_{\tilde{S}})$ we have $\ker(\phi) \ne 0$ if and only if $\phi = 0$, so the determinant
$$\det : \Hom(\ti{E},\ti{E} \otimes K_{\tilde{S}}) \to H^0(\tilde{S},K_{\tilde{S}}^{\otimes 2})$$ 
satisfies $\det^{-1}(0) = \{0\}$. As $H^0(\tilde{S},K_{\tilde{S}}^{\otimes 2}) \simeq \bC$ and $\det$ is a polynomial function, necessarily 
$$\dim\Hom(\ti{E},\ti{E} \otimes K_{\tilde{S}}) \le 1 = h^0(\tilde{S},K_{\tilde{S}}). $$
Hence $\tr^\vee$ is an isomorphism.
\end{proof}

\begin{rem}
	When $\NS(S) = 0$, it follows from Proposition~\ref{pro-liftP1} that there is no obstruction to deforming a $\bP^1$-bundle over $S$ along any small deformation of $S$. 
\end{rem}

Now we can prove Proposition~\ref{pro-sec} for good $\bP^1$-bundles over a surface of algebraic dimension $0$.

\begin{pro}\label{pro-aaa=0}
Let $f: X\to \ti{S}$ be a good $\bP^1$-bundle (see Definition~\ref{Def-bonP1})
over a smooth compact K\"ahler surface $\ti{S}$ of algebraic dimension 0. There exists a deformation 
$$\Pi : \cX \xto{\ti{f}} \ti{\cS} \to \gD$$ 
of $f$ such that 
$\Pi$ is an algebraic approximation of $X$.
Moreover, up to shrinking $\gD$, the underlying deformation $\Pi$ of $X$ is 
$f^{-1}(Y)$-locally trivial for every subvariety $Y \subsetneq \ti{S}$.
\end{pro}

\begin{proof} 

Let $\mu : \ti{S} \to S$ be the map from $\ti{S}$ to its minimal model $S$.
By Lemma~\ref{lem-2toreaa}, there exists an algebraic approximation $\cS \to \gD$ of $S$ 
which preserves $\NS(S)$
and is $\gS$-locally trivial for every $\gS \subsetneq S$.  
By Lemma~\ref{lem-lifta=0frommin} and Proposition~\ref{pro-liftP1},
up to shrinking $\gD$
the family $\cS \to \gD$ can be lifted to a deformation $\ti{\cS} \to \gD$ of $\ti{S}$
which is $Y$-locally trivial for every subvariety $Y\subsetneq \ti{S}$, and $\ti{\cS} \to \gD$
can be further lifted to a deformation 
$$\cX \xto{\ti{f}} \ti{\cS} \to \gD$$ 
of  $f : X \to \ti{S}$ with $\ti{f}$ being a $\bP^1$-bundle.
As $\cS \to \gD$ is an algebraic approximation of $S$, the family $\ti{\cS} \to \gD$ is an algebraic approximation of $\ti{S}$. 
It follows from Corollary~\ref{cor-critprojP1} that $\cX \to \ti{\cS} \to \gD$ is an algebraic approximation of $X$. 
It remains to show that $\Pi : \cX \to \gD$ is
$f^{-1}(Y)$-locally trivial for every subvariety $Y \subsetneq \ti{S}$.

Let $\eta : \ti{S} \xto{\mu} S \xto{\nu} S_\can$ be the composition of $\mu$ with the contraction of all the $(-2)$-curves of $S$. By Lemma~\ref{lem-fini-2}, the exceptional divisor $D \subset \ti{S}$ of $\eta$ is the union of all the curves of $\ti{S}$. 
To finish the proof,
it suffices to show that $\Pi : \cX \to \gD$ is 
$f^{-1}(C)$-locally trivial for every connected component $C \subset D$ of $D$. 
Indeed, this implies that $\Pi : {\cX}  \to \gD$ is $f^{-1}(D)$-locally trivial by Lemma~\ref{lem-loctrivsm}.\emph{ii)},
so that for every $Y \subsetneq \ti{S}$, if $Y' \cnec Y \cap D$, 
then $\Pi : {\cX}  \to \gD$ is $f^{-1}(Y')$-locally trivial. 
Since $D$ contains all the curves of $\ti{S}$, $Y'' \cnec Y \bss D$ has dimension $0$, 
so $\Pi$ is also $f^{-1}(Y'')$-locally trivial by Lemma~\ref{lem-loctrivsm}.\emph{iv)} and Lemma~\ref{lem-loctrivsm}.\emph{ii)}. 
As $Y = Y' \sqcup Y''$, 
it follows again from Lemma~\ref{lem-loctrivsm}.\emph{ii)} that $\Pi$ is $f^{-1}(Y)$-locally trivial.

Let $C \subset D$ be a connected component of $D$. 
Since $S_\can$ has no curve, $\eta(C)$ is a point $p \in S_\can$.
Let $U_\can \subset S_{\can}$ be a Stein neighborhood of $p$ and 
let $U  \cnec \nu^{-1}(U_\can)$ and $\ti{U}  \cnec \eta^{-1}(U_\can)$. 
Since $\ti{\cS} \to \gD$ is $C$-locally trivial,   
by Lemma~\ref{lem-Yltexpl}
there exist an open subset $\ti{\cU} \subset \ti{\cS}$ 
and an isomorphism $\phi : \ti{\cU} \eto \ti{U} \times \gD$ over $\gD$
such that $\ti{\cU} \cap \ti{S} = \ti{U}$
and $\phi_{|\ti{U}} : \ti{U} \to \ti{U}$ is the identity. 
Let $\ti{\cC} \cnec \phi^{-1}(\ti{C} \times \gD) \subset \cU \subset \ti{\cS}$, 
namely the trivial deformation of $\ti{C}$ along $\gD$ induced by $\phi$.
To show that $\Pi : \cX \to \gD$ is 
$f^{-1}(C)$-locally trivial,
we will show that up to shrinking $U_\can$, 
$\ti{f}^{-1}(\ti{\cU}) \to \gD$ is a trivial deformation of $f^{-1}(\ti{U})$.

We have $H^2(\ti{U} \times \gD,\cO^\times_{\ti{U} \times \gD}) = 0$. Indeed,
as $\ti{U} \times \gD$ has only rational singularities, 
we have  
$H^2(\ti{U} \times \gD , \cO_{\ti{U} \times \gD}) \simeq 
H^2\(U_\can \times \gD,(\eta \times \Id_{\gD})_*\cO_{\ti{U} \times \gD}\)$,
which vanishes because $U_\can \times \gD$ is Stein. 
As $\ti{U} \times \gD$ deformation retracts to $C$, we have 
$H^3(\ti{U} \times \gD, \bZ) \simeq H^3(C, \bZ) = 0$. 
Hence $H^2(\ti{U} \times \gD,\cO^\times_{\ti{U} \times \gD}) = 0$ 
by the exponential exact sequence.
It follows that there exists a locally free sheaf $\cE$ over $\ti{\cU} \simeq \ti{U} \times \gD$
such that $\bP(\cE) \simeq \ti{f}^{-1}(\ti{\cU})$ over $\ti{\cU}$.

Through the isomorphism $\phi : \ti{\cU} \eto \ti{U} \times \gD$,
we regard $\cE$ as the deformation over $\gD$ of 
the locally free sheaf $E \cnec \cE_{|\ti{U}}$ on $\ti{U}$.
To show that $\ti{f}^{-1}(\ti{\cU}) \to \gD$ is a trivial deformation,
it then suffices to show that any deformation of $E$ along $\gD$ is trivial
(up to shrinking $\ti{U}$ and $\gD$).
Since up to shrinking $\ti{U}$, 
there exists a semi-universal deformation of $E$~\cite[Satz 4 in \S3]{MR630650},
it suffices to show that $\Ext^1(E,E) = 0$. 

Since $X \to \ti{S}$ is a good $\bP^1$-bundle, by definition
up to shrinking $U_\can$ there exists a reflexive sheaf $F_\can$ on $U_\can$
such that $\bP(E) \simeq \bP(F)$ over $\ti{U}$ where $F \cnec \mu^*\((\nu^*F_\can)/\torsion\)$.
Since $E \simeq F \otimes L$ for some line bundle $L$ over $\ti{U}$, it suffices to show that $\Ext^1(F,F) = 0$.
By Lemma~\ref{lem-Esnault}, we have a surjective morphism
$\cO_{\ti{U}}^N \tto F$ and since $F$ is locally free, this induces a surjective morphism
$$\cHom(F,\cO_{\ti{U}}^N) \tto \cHom(F,F).$$
For any coherent sheaf $\cF$ on $\ti{U}$,
since $R^2\eta_*\cF = 0$ and $H^i(U_\can,R^j\eta_*\cF)$ whenever $i > 0$, we have
 $H^2(\ti{U},\cF) = 0$, so we have a surjective morphism
 $$H^1\(\ti{U},\cHom(F,\cO_{\ti{U}}^N)\) \tto H^1\(\ti{U}, \cHom(F,F)\).$$
 Finally, as $H^1\({U}_\can,\eta_*\cHom(F,\cO_{\ti{U}}^N)\) = 0$ and
 $R^1\eta_*\cHom(F,\cO_{\ti{U}}) = 0$ by Lemma~\ref{lem-Esnault},
 we have $H^1\(\ti{U},\cHom(F,\cO_{\ti{U}}^N)\) = 0$. 
 Hence $\Ext^1(F,F) = H^1\(\ti{U}, \cHom(F,F)\) = 0$. 
\end{proof}

\section{Algebraic approximations of $\bP^1$-fibrations over an elliptic surface}\label{sec-a1}

\subsection{Lifting the $G$-equivariant tautological family}
\hfill

The main result of \S\ref{sec-a1} is 
Corollary~\ref{cor-ellipsuraa}, 
about the existence of algebraic approximations of
 compact K\"ahler $\bP^1$-fibrations $f : X \to S$
over an
elliptic surface $p:S \to B$. 
When $p$ has local sections at every point of $B$, 
by Proposition~\ref{pro-def-existfam} the elliptic surface $S$ 
has an algebraic approximation realized by its associated tautological family,  
and we will construct an algebraic approximation of $f$ by lifting the tautological family associated to $p$.
The $G$-equivariant version 
of this statement is formulated in Proposition~\ref{pro-ellipsuraa}, 
which we prove first, 
then Corollary~\ref{cor-ellipsuraa} will follow 
as a consequence.

\begin{pro}\label{pro-ellipsuraa}
	Let $f : X \to S$ be a $\bP^1$-fibration 
	over an elliptic surface $p:S \to B$, itself over a smooth projective curve $B$. 
	Let $G$ be a finite group acting faithfully on $X$, $S$, and $B$ such that both $f$ and $p$ are $G$-equivariant. Assume the following conditions:
	\begin{itemize}
		\item 
	    $X$, $S$ are compact K\"ahler manifolds.	
		\item There exist $N \in \{I,II,III\}$ and a Zariski dense open subset $B^\circ \subset B$ 
		such that for every $b \in B^\circ$, if $X_b \cnec (p \circ f)^{-1}(b)$ and $S_b = p^{-1}(b)$, then $S_b$ is smooth and $X_b  \to S_b$ is a ruled surface  of type $N$ 
		(see Definition~\ref{Def-typee=0}).
		\item $p:S \to B$ has local holomorphic sections at every point of $B$.
	\end{itemize}
	Then the $G$-equivariant tautological family (see \S\ref{ssec-famtautfibellip} and Proposition~\ref{pro-def-existfam})
	$$ \Pi_S :  \cS \xto{q}  B \times V \to V \cnec H^1(B,\bar{\cE}_{\bH/B})^G $$ 
	associated to the $G$-equivariant elliptic fibration $p$ can be lifted to a
	$G$-equivariant deformation 
	$$\Pi_X : \cX \to \cS \to B \times V \to V$$
	of $f$ 
	such that the underlying deformation of $\pi \colonec p \circ f: X \to B$ is $G$-equivariantly locally trivial over $B$. 
\end{pro}

The strategy of the proof
is simple. 
Recall from \S\ref{ssec-famtautfibellip} that $\Pi_S$ is constructed by
a family of \v{C}ech 1-cocycles $V \ni t \mapsto \{ \exp(\xi_{ij}(t))\}$ with coefficients in $\cAut^0_B(S)$. 
The heart of the proof consists of
lifting $t \mapsto \{\exp(\xi_{ij}(t))\}$  
to a deformation of 1-cocycles with coefficients in $\cAut^0_B(X)$,
so that we can glue the local fibrations $f^{-1}(U_i) \cnec X_i \to U_i$ using these 1-cocycles and obtain a deformation of $f$ lifting $\Pi_S$. 
The construction of such a lifting depends on
whether $\Aut^0(X_b)$ is commutative, which is the case
if the ruled surface $X_b$ is of type II or III.
We therefore separate the case where $X_b$ is type I from type II or II in the proof.
In the former case, $X_b$ is simply a trivial $\bP^1$-bundle and
we construct the lifting by hand; 
in the latter case, we use cohomological methods to construct the lifting.

\begin{proof}[Proof of Proposition~\ref{pro-ellipsuraa}]

First we fix some notations. 
For every $b \in B$, let $X_b \cnec \pi^{-1}(b)$ and $S_b \cnec p^{-1}(b)$. 
By assumption, there exist $N \in \{I,II,III\}$ and a Zariski dense open subset $B^\circ \subset B$ such that for every $b \in B^\circ$, 
$S_b$ is a smooth elliptic curve and
$X_b$ is a ruled surface over $S_b$ of type $N$. 
Let $r : B \to B/G =: \bar{B}$ be the quotient.
Up to shrinking $B^\circ$, we can assume
that $B^\circ$ is $G$-stable (so that $r^{-1}(r(B^\circ)) = B^\circ$)
and that $\bar{B}^\circ \cnec r(B^\circ)$ does not contain any branch point of $r$. 
Let ${p}^\circ : {S}^\circ \to {B}^\circ$ and ${f}^\circ : {X}^\circ \to {S}^\circ$ denote the restriction of ${p}$ and ${f}$ to ${S}^\circ \cnec {p}^{-1}({B}^\circ)$ and ${X}^\circ \cnec {f}^{-1}({S}^\circ)$ respectively. Let ${\pi}^\circ \cnec {p}^\circ \circ {f}^\circ : {X}^\circ \to {B}^\circ$. 
We also define the maps
\begin{equation}\label{diag-XSUV}
\begin{tikzcd}[cramped, row sep = 15, column sep = 10]
X^\circ \times V \ar[dr, swap, "\pi^\circ_V \cnec \pi^\circ \times \Id_V"] \arrow[rr, "f^\circ_V \cnec  f^\circ \times \Id_V"] & & S^\circ \times V \ar[dl, "p^\circ_V \cnec p^\circ \times \Id_V"]   \\
& {B^\circ} \times V.  &    
\end{tikzcd}
\end{equation}

The proof of Proposition~\ref{pro-ellipsuraa} consists of several steps.

\noindent \textbf{Step 0: Construction of a good open cover of $B$.}
	
	In this step, we will construct a good open cover of $B$
	satisfying the properties listed in Lemma~\ref{lem-bonrecG}.
	
	\begin{lem}\label{lem-bonrec}
	There exists a finite good open cover $\fV  = \{V_i\}_{i \in \bar{I}}$ of $\bar{B}$ satisfying the following properties.
\begin{enumerate}[label = \roman{enumi})]
		\item Each $V_i$ contains at most one branch point of  $r : B \to \bar{B}$. 
		\item For each pair of indices $i, j \in \bar{I}$ such that $i \ne j$, we have ${V}_{ij} \cnec {V}_i \cap {V}_j \subset \bar{B}^\circ = r(B^\circ)$.
		In particular, if $V_i$ and $V_j$ both contain  branch points of $r$, then either $V_i \cap V_j = \emptyset$ or $i = j$. 
		\item The fibration $p : S \to B$ has local sections over each $r^{-1}(V_i) \subset B$.
	\end{enumerate}
	\end{lem} 
\begin{proof}
	
	Since $r : B \to \bar{B}$ is finite, for every $b \in \bar{B}$ and every $b' \in r^{-1}(b)$, every neighborhood of $b'$ contains a connected component of $r^{-1}(V_b)$ for some neighborhood $V_b \subset \bar{B}$ of $b$. Thus, since $p : S \to B$ has local sections around every point of $B$, we obtain the following assertion.
	\begin{claim}\label{claim-recii)}
		For every $b \in \bar{B}$, we can find a neighborhood $V_b \subset \bar{B}$ of $b$ such that the fibration $p : S \to B$ has local sections over $r^{-1}(V_b)$.
	\end{claim}
	
	    Let $\fT$ be a finite triangulation of $\bar{B}$ and let $\bar{I}$ be the set of vertices of $\fT$.   For every $i \in \bar{I}$, let $V_i \subset \bar{B}$ be the union of (open) strata of the CW-complex defined by $\fT$ which are adjacent to $i$. Up to refining $\fT$,  we can assume that $\fT$ satisfies the following properties.
	       
	    \begin{itemize}
	    	\item The finite set $\gS \cnec \bar{B} \bss \bar{B}^\circ$ is contained in $\bar{I}$.
	    	\item For every $i \in \bar{I}$, the fibration $p : S \to B$ has local sections over $r^{-1}(V_i)$ (by Claim~\ref{claim-recii)}).
		\end{itemize}
	   By construction, the collection $\fV \cnec \{V_i\}_{i \in \bar{I}}$ forms a good open cover of $\bar{B}$ satisfying the properties listed in Lemma~\ref{lem-bonrec}.
\end{proof}

	Let  $\fV = \{V_j\}_{j \in \bar{I}}$  be a good open cover of $\bar{B}$ as in Lemma~\ref{lem-bonrec} and let $\fU =\{U_i\}_{i \in {I}}$ be the open cover of ${B}$ consisting of connected components of $r^{-1}({V_j})$ with $j$ runs through $\bar{I}$.
	
	\begin{lem}\label{lem-bonrecG}
	The collection $\fU = \{U_i\}_{i \in {I}}$ is a finite good Stein open cover of ${B}$ satisfying the following properties.
	\begin{enumerate}[label = \roman{enumi})]
		\item The fibration ${p} : {S} \to {B}$ has local sections over each $U_i$.
		\item For each pair of indices $i_1, i_2 \in I$ such that $i_1 \ne i_2$, we have $U_{i_1i_2} \cnec U_{i_1} \cap U_{i_2} \subset r^{-1}(r(B^\circ)) = B^\circ$.
		\item The open cover $\fU$ is $G$-invariant. Moreover for every $g \in G$ and every pairwise distinct indices $i_1 ,\ldots , i_j \in {I}$ with $j \ge 2$, we have $g(U_{i_1 \cdots  i_j}) \cap U_{i_1 \cdots  i_j} = \emptyset$ as long as $g \ne \Id_{{B}}$ where $U_{i_1 \cdots  i_j} \cnec \cap_{l=1}^j U_{i_l}$.
	\end{enumerate}
\end{lem} 

\begin{proof}
	By Lemma~\ref{lem-bonrec} and the definition of $\fU$, it follows easily that $\{U_i\}_{i \in {I}}$ is a $G$-invariant open cover of ${B}$ with $|{I}| < \infty$ 
	and satisfying \emph{i)} and \emph{ii)}.  
	By Lemma~\ref{lem-bonStein}, $\fU$ is Stein. It remains to prove \emph{iii)} and that $\fU$ is good.

	By Lemma~\ref{lem-bonrec}.\emph{i)} and the definition of $\fU$, 
	$U_l \to r(U_l)$ is a finite cover with at most one ramification point 
	for every $l \in {I}$ and $r(U_l)$ is contractible. 
	So a deformation retraction of $r({U_i})$ to the branch point (or to any point if there is none) yields a contraction of ${U_i}$ to a point, 
	which shows that $U_i$ is contractible.
	Given two distinct open subsets $U_{i_1}$ and $U_{i_2}$ in $\fU$ such that $U_{i_1} \cap U_{i_2} \ne \emptyset$, 
	we have $r(U_{i_1}) \cap r(U_{i_2}) \ne \emptyset$. 
	Since $U_l$ is a connected component of $r^{-1}(r(U_l))$ for every $l \in {I}$,  $U_{i_1} \cap U_{i_2} \ne \emptyset$ implies that $r(U_{i_1}) \ne r(U_{i_2})$.
	So by the second statement of
	Lemma~\ref{lem-bonrec}.\emph{ii)}, either $U_{i_1} \to r(U_{i_1})$ or $U_{i_2} \to r(U_{i_2})$ is biholomorphic. 
	Therefore if $U_{i_1},\ldots,U_{i_j}$ are pairwise distinct open subsets in $\fU$ with $j \ge 2$, 
	then $U_{i_1 \cdots  i_j}  \to r(U_{i_1 \cdots  i_j} )$ is biholomorphic. 
	It follows that $U_{i_1 \cdots  i_j} \simeq r(U_{i_1 \cdots  i_j})$ is contractible (because $\fV$ is a good open cover), so $\fU$ is a good open cover of $B$. 
	It also follows that each $U_{i_1 \cdots  i_j}$ is a connected component of $r^{-1}(r(U_{i_1 \cdots  i_j}))$. 
	As $r$ is the quotient by the faithful $G$-action, 
	it follows that 
	$g(U_{i_1 \cdots  i_j}) \cap U_{i_1 \cdots  i_j} = \emptyset$ 
	provided $g \ne \Id_{{B}}$, which proves \emph{iii)}. 
	\end{proof}

From now on, we fix an open cover $\fU = \{U_i\}_{i \in I}$ of $B$ 
as in Lemma~\ref{lem-bonrecG}.
For every multi-index $i_1,\ldots,i_l \in I$, we define ${S}_{i_1\cdots i_l} \cnec {p}^{-1}({U}_{i_1\cdots i_l})$ and ${X}_{i_1\cdots i_l} \cnec {\pi}^{-1}({U}_{i_1\cdots i_l})$. 
The restriction of $p$ (resp. $\pi$ and $f$)  to ${S}_{i_1\cdots i_l}$ 
(resp. ${X}_{i_1\cdots i_l}$)
are denoted by $p_{i_1\cdots i_l} : {S}_{i_1\cdots i_l} \to {U}_{i_1\cdots i_l}$
(resp. $f_{i_1\cdots i_l} : {X}_{i_1\cdots i_l} \to {S}_{i_1\cdots i_l}$ and 
$\pi_{i_1\cdots i_l} : {X}_{i_1\cdots i_l} \to {U}_{i_1\cdots i_l}$).

\noindent \textbf{Step 1: Tautological family.}

In this step, we recall from  \S\ref{ssec-famtautfibellip} 
the construction of the $G$-equivariant 
tautological family $\Pi_S$ in Proposition~\ref{pro-ellipsuraa}
in order to fix the notations.

Let $J \to B^\circ$ be the Jacobian fibration associated to the smooth elliptic fibration $p^\circ : S^\circ \to B^\circ$ and let ${\bH} \cnec R^1{p}^\circ_*\bZ$. 
Let $ \cJ$ be the sheaf of local sections of the Jacobian fibration $J \times V \to B^\circ \times V$ associated to $p^\circ_V : S^\circ \times V \to B^\circ \times V$.
Since $R^1(p^\circ_V)_*\bZ \simeq \pr_1^{-1}\bH$ 
where $\pr_1 : B^\circ \times V \to B^\circ$ is the projection to the first factor, we have
$$\cJ \simeq \cE_{\pr_1^{-1}\bH/B^\circ \times V}/\pr^{-1}\bH, $$
where we recall that $\cE_{\pr_1^{-1}\bH/B^\circ \times V}$ is defined by~\eqref{def-EHB} and satisfies 
$$\cE_{\pr_1^{-1}\bH/B^\circ \times V} \simeq (p^\circ_V)_*T_{S^\circ \times V / B^\circ \times V}  
\simeq \pr_1^*p^\circ_*T_{S^\circ / B^\circ }
\simeq \pr_1^*\cE_{\bH/B^\circ}.$$
We also recall  
from Example~\ref{ex-fibab} that 
$\cJ \simeq \cAut^0_{B^\circ \times V}(S^\circ\times V)$,
and together with
$\cE_{\pr_1^{-1}\bH/B^\circ \times V} 
\simeq (p^\circ_V)_*T_{S^\circ \times V / B^\circ \times V}$ mentioned above,
these isomorphisms commute with the exponential maps:
\begin{equation}
\begin{tikzcd}[
cramped, row sep = 8, column sep = 10,
]
\pr_1^*\cE_{\bH/B^\circ} \simeq \cE_{\pr_1^{-1}\bH/B^\circ \times V} \ar[d,no head,"\wr"] \arrow[r, "\exp_{\pr_1^{-1}\bH/B^\circ \times V}"] & \cJ \ar[d,no head ,"\wr"]   \\
(p^\circ_V)_*T_{S^\circ \times V / B^\circ \times V} \ar[r, "\exp"] & \cAut^0_{B^\circ \times V}(S^\circ\times V).    
\end{tikzcd}
\end{equation}

Recall from \S\ref{ssec-famtautfibellip} that 
the map $q :  {\cS} \to  {B} \times V$
defining the $G$-equivariant tautological family $\Pi_S$ in Proposition~\ref{pro-ellipsuraa},
is obtained by gluing  the local fibrations $\Set{p_i \times \Id_V : {S}_i \times V \to  {U}_i \times V} $ 
using the $G$-equivariant 1-cocycle of translations 
$\Set{e_{ij} :  {S}_{ij} \times V \to  {S}_{ij} \times V}$ defined by 
$$e_{ij} \cnec \exp( {U}_{ij} \times V)(\xi_{ij}) \in \cJ({U}_{ij} \times V) \subset \Aut(S_{ij} \times V),$$
for some $G$-invariant 1-cocycle 
$\Set{\xi_{ij} \in (\pr_1^*\cE_{\bH/B^\circ})(U_{ij} \times V) }$.
As in \S\ref{ssec-famtautfibellip}, 
we regard each $\xi_{ij} \in (\pr_1^*\cE_{\bH/B^\circ})(U_{ij} \times V)$ 
as a map $\xi_{ij} : V \to {\cE}_{\bH/B^\circ} (U_{ij})$. 
So for every $v \in V$, the restriction $e_{ij| S_{ij} \times \{v\}}$ is equal to the translation
$$e_{ij| S_{ij} \times \{v\}} =  \tr\(\exp_{\bH/B^\circ}( {U}_{ij})(\xi_{ij}(v))\) : S_{ij} \times \{v\} \eto S_{ij} \times \{v\}$$
where $\exp_{\bH/B^\circ} : \cE_{\bH/B^\circ} \to \cJ_{\bH/B^\circ}$ is the quotient as in~\eqref{SE-Jac}. 
As we mentioned in \S\ref{ssec-famtautfibellip}, 
the 1-cocycle $\{\xi_{ij}\}$ can be chosen to satisfy
$\xi_{ij}(0) = 0$ (see~\cite[Lemma 4.19]{HYLbimkod1}), 
so
\begin{equation}\label{eqn-ideij0}
e_{ij| {S}_{ij} \times \{0\}} = \Id_{ {S}_{ij} \times \{0\}}.
\end{equation}

\noindent \textbf{Step 2: Case type I.}

In this step, 
we prove Proposition~\ref{pro-ellipsuraa} in the case where $X_b$ is of type I  for every $b \in B^\circ$, 
namely ${X}_b \simeq  {S}_b \times \bP^1$ over $ {S}_b$ by Definition~\ref{Def-typee=0}. 
We first prove the following.

\begin{lem}\label{lem-nonobsProjellpaff}
	Let $g : \sS \to \gD$ be a smooth family of curves over a contractible Stein space $\gD$. Every $\bP^n$-bundle over ${\sS}$ is the projectivization of some vector bundle.
\end{lem}

\begin{proof}
	Since the obstruction to lifting a projective bundle to a vector bundle is an element in $H^2({\sS},\cO_{\sS}^\times)$, it suffices to show that $H^2({\sS},\cO_{\sS}^\times) = \{0\}$. Since the fibers of $g$ are curves, we have $R^ig_*\cO_{\sS} = 0$ for every $i \ge 2$. We also have $H^p(\gD,R^qg_*\cO_{\sS}) = 0$ for every $p \ge 1$ because $\gD$ is Stein, so the Leray spectral sequence implies $H^i({\sS},\cO_{\sS}) =0$ for $i = 2$ or $3$. It follows from the exponential sheaf sequence that $H^2({\sS},\cO_{\sS}^\times) \simeq H^3({\sS},\bZ)$. As $\gD$ is contractible, $\sS$ is homotopy equivalent to a curve. Hence $H^2({\sS},\cO_{\sS}^\times) \simeq H^3({\sS},\bZ) = 0$.
\end{proof}

Based on the assumption that ${X}_b \simeq  {S}_b \times \bP^1$ over $ {S}_b$ for every $b \in B$, we can prove the following.

\begin{lem}\label{lem-trivtypeI}
	We have ${X}_{ij} \simeq {S}_{ij} \times \bP^1$ over ${S}_{ij}$.
\end{lem}

\begin{proof}
 
	As $U_{ij}$ is contractible and Stein, 
	by Lemma~\ref{lem-nonobsProjellpaff} there exists a locally free sheaf $\cE_{ij}$ over ${S}_{ij}$ 
	such that ${X}_{ij} \simeq \bP(\cE_{ij})$ over ${S}_{ij}$. 
	For every $b \in U_{ij}$, let $\cE_b \cnec \cE_{ij| S_b}$.
	Since $ {X}_b \simeq  {S}_b \times \bP^1$ over $ {S}_b$, we have $\cE_b \simeq \cL_b \oplus \cL_b$ for some invertible sheaf $\cL_b$ on $S_b$. Fix $o \in U_{ij}$. Since $U_{ij}$ is contractible, 
	there exists an invertible sheaf $\cL$ on $S_{ij}$ such that $\cL^{\otimes 2} \simeq \det\cE_{ij}$ and $\cL_{|S_o} \simeq \cL_o$. 
	Let $\cE' \cnec \cE_{ij} \otimes \cL^\vee$.
	\begin{claim}\label{claim-oc}
		The condition $\cE'_b \cnec \cE'_{| {S}_{b}} \simeq \cO_{ {S}_{b}} \oplus \cO_{ {S}_{b}}$ is open and closed in $b \in U_{ij}$.
	\end{claim}
	
	\begin{proof}
	As $U_{ij}$ is contractible, 
	there exist $\cT_0 \cnec \cO_{S_{ij}}, \cT_1,\cT_2,\cT_3 \in \Pic^0(S_{ij})$ 
	such that $\cO_{S_{b}}, {\cT_1}_{|S_b},{\cT_2}_{|S_b},{\cT_3}_{|S_b}$ are all the
	2-torsion line bundles on the elliptic curve $S_b$.
	Since $\det\cE'_b \simeq \cO_{S_b}$ and $\cE'_b \simeq \cL'_b \oplus \cL'_b$ for some line bundle $\cL'_b$, we have $\cE'_b \simeq {\cT_i}_{|S_b} \oplus {\cT_i}_{|S_b}$ for some $i$.
	As $H^0(S_b,\cT_i) \ne 0$ if and only if $i = 0$ and 
	$b \mapsto H^0(S_b,\cE'_b)$ is upper semi-continuous,
	the property $\cE'_b  \simeq \cO_{ {S}_{b}} \oplus \cO_{ {S}_{b}}$
	is thus closed in $b$.
	
	Next, assume to the contrary that $\cE'_b  \simeq \cO_{ {S}_{b}} \oplus \cO_{ {S}_{b}}$ is not open at $b \in U_{ij}$, 
	then there exist $i \in \{1,2,3\}$
	and a sequence $b_\ell \in U_{ij} $ converging to $b$
	such that $\cE'_{b_\ell}  \simeq {\cT_i}_{|S_{b_\ell}} \oplus {\cT_i}_{|S_{b_\ell}}$ for all $\ell$.
	We then have $H^0(S_b,\cE'\otimes \cT_i) =H^0(S_b, \cT_i)^{\oplus 2} =  0$ but 
	$H^0(S_{b_\ell},\cE'\otimes \cT_i) = H^0(S_{b_\ell},\cO_{S_{b_\ell}})^{\oplus 2} \ne 0$,
	which violates the upper semi-continuity.
	\end{proof} 
	
	As $\cE'_o = \cE'_{| {S}_{o}} \simeq \cO_{ {S}_{o}} \oplus \cO_{ {S}_{o}}$
	and $U_{ij}$ is connected,
	Claim~\ref{claim-oc} implies that
	$\cE'_b \simeq \cO_{ {S}_{b}} \oplus \cO_{ {S}_{b}}$ for every $b \in U_{ij}$. 
	It follows that $({p}_{|S_{ij}})_*\cE'$ is locally free of rank 2 over $U_{ij}$ by Grauert's base change, so $({p}_{|S_{ij}})_*\cE' \simeq \cO_{U_{ij}} \oplus \cO_{U_{ij}}$ because $U_{ij}$ is biholomorphic to a disc. 
	Therefore $\cE'$ has two global sections $\gs$ and $\gs'$ 
	which are linearly independent everywhere over ${S}_{ij}$, 
	so $\cE'$ is trivial. Hence ${X}_{ij} \simeq \bP(\cE') \simeq {S}_{ij} \times \bP^1$ over ${S}_{ij}$.
\end{proof}

For every $i \ne j \in I$, we fix an isomorphism $\imath_{ij} : {X}_{ij} \eto {S}_{ij} \times \bP^1$ over $S_{ij}$ as in Lemma~\ref{lem-trivtypeI}. 
Since $G$ acts freely on $\Set{U_{ij} | i \ne j \in I }$ by Lemma~\ref{lem-bonrecG}.\emph{iii)},
we can choose $\imath_{ij}$ in such a way that the collection $\{\imath_{ij}\}$ is $G$-invariant. Let $\ti{e}'_{ij} \cnec e_{ij} \times \Id_{\bP^1} \in \Aut({S}_{ij} \times V  \times \bP^1)$ and 
$$\tilde{\imath}_{ij} = \imath_{ij} \times \Id_V : {X}_{ij} \times V 
\eto {S}_{ij} \times \bP^1  \times V
\simeq {S}_{ij} \times V  \times \bP^1.$$ 
Let  $i,j,k \in I$ and let $b \in U_{ijk}$. 
Since  $\imath_{ki|{X}_{b}} \circ \imath_{jk|{S}_{b} \times \bP^1}^{-1}$ is an automorphism of ${S}_{b} \times \bP^1$ over ${S}_{b}$, we have $\imath_{ki|{X}_{b}} \circ \imath_{jk|{S}_{b} \times \bP^1}^{-1} \in \Id_{S_b} \times \Aut(\bP^1)$
(because $S_b$ is proper and $\Aut(\bP^1)$ is affine). 
It follows that $\tilde{\imath}_{ki} \circ \tilde{\imath}_{jk}^{-1}$ commutes with $\ti{e}'_{jk}$ 
(regarded as automorphisms of ${S}_{ijk} \times V  \times \bP^1$). 
For the same reason,  $\tilde{\imath}_{ki} \circ \tilde{\imath}_{ij}^{-1}$ commutes with $\ti{e}'_{ij}$. 
Thus if we define $\ti{e}_{ij} \cnec \tilde{\imath}^{-1}_{ij} \circ \ti{e}'_{ij} \circ \tilde{\imath}_{ij} \in \Aut(X_{ij} \times V)$,
then 
$$\ti{e}_{ki} \circ \ti{e}_{jk}\circ \ti{e}_{ij} = 
(\tilde{\imath}^{-1}_{ki} \circ \ti{e}'_{ki} \circ \tilde{\imath}_{ki})
\circ (\tilde{\imath}^{-1}_{jk} \circ \ti{e}'_{jk} \circ \tilde{\imath}_{jk}) 
\circ (\tilde{\imath}^{-1}_{ij} \circ \ti{e}'_{ij} \circ \tilde{\imath}_{ij}) = \tilde{\imath}^{-1}_{ki} \circ (\ti{e}'_{ki}\circ \ti{e}'_{jk}\circ \ti{e}'_{ij}) \circ \tilde{\imath}_{ki} = \tilde{\imath}^{-1}_{ki} \circ \tilde{\imath}_{ki} = \Id_{{X}_{ijk}},$$
where the second equality follows from the above commutative relations,
and the third  equality follows from the property that $\Set{\te'_{ij} = e_{ij} \times \Id_{\bP^1} }$ is a 1-cocycle (because $\Set{e_{ij} \in \Aut({S}_{ij} \times V) }$ is so).
Therefore we can glue the collection of fibrations 
$\Set{ f_i \times \Id_V : {X}_i \times V \to {S}_i \times V}_{i \in I}$ 
along $\Set{f_{ij} \times \Id_V : {X}_{ij} \times V \to {S}_{ij} \times V}$ 
using the $G$-invariant 1-cocycle of automorphisms $\Set{\ti{e}_{ij} \in \Aut({X}_{ij} \times V) }$ 
which lifts $\Set{e_{ij} \in \Aut({S}_{ij} \times V) }$,
and obtain a $G$-equivariant map ${\cX} \to {\cS}$. 
Let 
$$\Pi_X : {\cX} \to {\cS} \to {B} \times V \to V$$ 
be the composition of $\cX \to \cS$ with $\Pi_S$. 
By~\eqref{eqn-ideij0}, $e_{ij| {S}_{ij} \times \{0\}} = \Id_{ {S}_{ij} \times \{0\}}$ 
and $\ti{e}_{ij |{X}_{ij} \times \{0\}} = \Id_{{X}_{ij}}$ for every $i$ and $j$, so the central fiber of the family $\Pi_X$ of fibrations over $B$ is $\pi : X \to B$.
That $\Pi_X$ is
a deformation of $\pi$ 
which is $G$-equivariantly locally trivial over ${B}$ follows from the construction
and the property that $\Pi_S$ is $G$-equivariantly locally trivial over ${B}$
(Proposition~\ref{pro-def-existfam}).

\noindent \textbf{Step 3: Case type II and III.}

From now on, we assume that 
the ruled surface $X_b$ is purely of type II or III for every $b \in B^\circ$.

\begin{lem}\label{lem-condsurjlisse}
	The commutative diagram~\eqref{diag-XSUV}
	satisfies Assumption~\ref{assump-autrel}.
\end{lem}
\begin{proof}
	
	Since~\eqref{diag-XSUV} is the base change of $X^\circ \xto{f^\circ} S^\circ \xto{p^\circ} B^\circ$ under the projection $B^\circ \times V \to B^\circ$, it suffices to prove that
	$X^\circ \xto{f^\circ} S^\circ \xto{p^\circ} B^\circ$ satisfies Assumption~\ref{assump-autrel}.
	
	Both ${p}^\circ : {S}^\circ \to {B}^\circ$ and ${f}^\circ : {X}^\circ \to {S}^\circ$ are smooth by assumption. Let $b \in B^\circ$.  
	As $S_b$ is an elliptic curve, we have $h^0(S_b,T_{S_b}) = 1$. 
	Since $X_b$ is a ruled surface of type II or III, the function $B^\circ \ni b \mapsto h^0(X_b,T_{X_b})$ is constant by Lemma~\ref{lem-SeilerDefruled},  which proves
	that $X^\circ \xto{f^\circ} S^\circ \xto{p^\circ} B^\circ$ satisfies Assumption~\ref{assump-autrel}.\emph{i)}.

	As $S_b$ is an elliptic curve, $\Aut(S_b)$ is commutative.
	Since $X_b$ is a ruled surface of type  II or III over an elliptic curve, 
	the commutativity of $\Aut^0(X_b)$ follows from~\cite[Theorem 3.(2) or (3)]{MaruyamaRuled}, according to whether $X_b$ is of type II or III respectively.
	Hence $X^\circ \xto{f^\circ} S^\circ \xto{p^\circ} B^\circ$ satisfies Assumption~\ref{assump-autrel}.\emph{ii}).
\end{proof}

Recall from Step 1 that we have an isomorphism between the exponential maps
$$ \exp_{\pr_1^{-1}\bH/B^\circ \times V} : \cE_{\pr_1^{-1}\bH/B^\circ \times V} \to \cJ$$
and
$$ \exp : (p^\circ_V)_*T_{S^\circ \times V / B^\circ \times V} \to \cAut^0_{B^\circ \times V}(S^\circ\times V).$$
For simplicity, let 
$$ \cT \cnec ({\pi}^\circ_V)_*T_{X^\circ \times V/ S^\circ \times V} \simeq \pr_1^*\pi^\circ_*T_{X^\circ / S^\circ},     \ \ \ \ \cT' \cnec ({\pi}^\circ_V)_*T_{X^\circ \times V/ B^\circ \times V}\simeq \pr_1^*\pi^\circ_*T_{X^\circ / B^\circ},  
$$
$$
\cE \cnec   \cE_{\pr_1^{-1}\bH/B^\circ \times V}
\simeq  (p^\circ_V)_*T_{S^\circ \times V / B^\circ \times V} \simeq \pr_1^*p^\circ_*T_{S^\circ / B^\circ}
\simeq \pr_1^*\cE_{\bH/B^\circ},
$$
where $\pr_1 : B^\circ \times V \to B^\circ$ is the projection to the first factor, and (with the notations introduced in \S\ref{ssec-autrel})  let
$$ \cA \cnec \(({p}^{\circ}_V)_* \cAut_{S^\circ \times V}(X^\circ \times V)\) \cap \cA',  \ \ \ \ \ \ \ \ \cA' \cnec  \cAut_{B^\circ \times V}^0(X^\circ \times V),  \ \ \ \ \ \ \ \ \cJ =   \cAut^0_{B^\circ \times V}(S^\circ\times V). $$ 
Since~\eqref{diag-XSUV} satisfies Assumption~\ref{assump-autrel} by Lemma~\ref{lem-condsurjlisse}, we have the commutative diagram of morphisms of sheaves of abelian groups
\begin{equation}\label{suitex-extautregleeJ}
\begin{tikzcd}[cramped, row sep = 20, column sep = 40]
0  \ar[r] & \cT   \ar[d,"\exp_1"] \ar[r] & \cT' \ar[d,"\exp_2"]\ar[r, ""] & \cE  \ar[d,"\exp_3 = \exp"] \ar[r] & 0 \\
0  \ar[r] & \cA \ar[r] & \cA' \ar[r, "\Psi"] & \cJ \ar[r] & 0,
\end{tikzcd}
\end{equation}
which is~\eqref{dc-faisTA} defined for~\eqref{diag-XSUV} instead of~\eqref{CD-XYZ}. 
Since each map in~\eqref{diag-XSUV} is $G$-equivariant, the $G$-actions on each variety in~\eqref{diag-XSUV} induce a $G$-sheaf structure on each sheaf in~\eqref{suitex-extautregleeJ} by conjugation and each morphism in~\eqref{suitex-extautregleeJ} is $G$-equivariant. 

\begin{lem}
	The commutative diagram~\eqref{suitex-extautregleeJ} 
	has exact rows.
\end{lem}
\begin{proof}
	By Lemma~\ref{lem-pPsisurj}, it suffices to show that the descent morphism $\Aut^0_{B^\circ \times V}(X^\circ \times V) \to \Aut^0_{B^\circ \times V}(S^\circ \times V)$ is surjective. Since~\eqref{diag-XSUV} is the base change of $X^\circ \xto{f^\circ} S^\circ \xto{p^\circ} B^\circ$ under the projection $B^\circ \times V \to B^\circ$, it suffices to show that the descent morphism $\Aut^0(X_b) \to \Aut^0(S_b)$ is surjective for every $b\in B^\circ$. 
	The latter follows from~\cite[Lemma 8]{MaruyamaRuled} 
	since $X_b$ is a ruled surface over $S_b$ of type II or III.
\end{proof}

For every open subset $U \subset B^\circ$ and $\gs \in \cT'(U \times V)$,
since $\cT' \simeq \pr_1^*\pi^\circ_*T_{X^\circ / B^\circ}$, 
we can regard $\gs$ as a map $\gs : V \to (\pi^\circ_*T_{X^\circ/ B^\circ})(U)$.
Similarly, all
$\gs_1 \in \cT(U \times V)$ and  $\gs_2 \in \cE(U \times V)$ 
are also be regarded as maps $\gs_1 : V \to (\pi^\circ_*T_{X^\circ/ S^\circ})(U)$ and
$\gs_2 : V \to  (p^\circ_*T_{S^\circ/ B^\circ})(U)$.

Recall that we introduced in Step 1 the $G$-invariant 1-cocycle
 $$\Set{{\xi}_{ij} \in  \cE(U_{ij} \times V)}  
 = \Set{\xi_{ij} : V \to  {\cE}_{ {\bH}/ {B^\circ}}( {U}_{ij}) \simeq p^\circ_*T_{S^\circ/ B^\circ}(U_{ij}) },$$   
which is used to construct the tautological family $\Pi_S$. 
 \begin{claim}\label{claim-tixirel}
 There exists a $G$-invariant 1-cochain 
 $$\Set{\ti{\xi}_{ij} \in \cT'({U}_{ij} \times V) }
 = \Set{\ti{\xi}_{ij} : V \to \pi^\circ_*T_{X^\circ/ B^\circ}(U_{ij}) }$$
 lifting the 1-cochain $\{\xi_{ij}\}$ 
 and satisfying ${\ti{\xi}}_{ij}(0) = 0$. 	
 \end{claim}
 \begin{proof}
 	 
Since $U_{ij}$ and ${U}_{ij} \times V$ are Stein, by~\eqref{suitex-extautregleeJ} we have the commutative diagram 
 		\begin{equation}
 \begin{tikzcd}[cramped, row sep = 15, column sep = 20]
 0 \ar[r] &  \cT({U}_{ij} \times V)   \ar[d, "\imath_1"] \ar[r] & \cT'({U}_{ij} \times V)   \arrow[r] \ar[d, "\imath_2"] & \cE({U}_{ij} \times V)  \ar[r] \ar[d,"\imath_3"]  & 0  \\
 0 \ar[r] & \cT_{|B^\circ \times 0}({U}_{ij}\times 0)  \ar[r] &  \cT'_{|B^\circ \times 0}({U}_{ij}\times 0) \arrow[r]  &  \cE_{|B^\circ \times 0}({U}_{ij}\times 0) \ar[r] & 0 
 \end{tikzcd}
 \end{equation}
  with exact rows, where the vertical arrows are the restrictions to 
  $B^\circ \times 0 \subset B^\circ \times V$. Since $\cT \simeq \pr_1^*\pi^\circ_*T_{X^\circ / S^\circ} \simeq \pr_1^* \cT_{|B^\circ \times 0} $, the map $\imath_1$  is surjective, so the induced map $\ker(\imath_2) \to \ker(\imath_3)$ is also surjective by the snake lemma. Therefore since $\xi_{ij}(0) = 0$ by assumption, we can find a lifting $\ti{\xi}_{ij}  \in \cT'({U}_{ij} \times V)$ of $\xi_{ij} \in \cE({U}_{ij} \times V)$ such that ${\ti{\xi}}_{ij}(0) = 0$. 
  As $\{{\xi}_{ij}\}$ is $G$-invariant,
  up to replacing the 1-cochain $\{{\ti{\xi}}_{ij}\}$ by 
  $\frac{1}{|G|}\sum_{g \in G} g \cdot \{{\ti{\xi}}_{ij}\}$, we can assume that the lifting $\{{\ti{\xi}}_{ij}\}$ of $\{{\xi}_{ij}\}$ is $G$-invariant.
\end{proof}

In general we do not know
whether the 1-cochain $\{\ti{\xi}_{ij}\}$ constructed above is a 1-cocycle. 
The next lemma  
shows that up to modifying $\{\ti{\xi}_{ij}\}$, 
a $G$-invariant 1-cocycle lifting $\{\xi_{ij}\}$ exists. 

\begin{lem}\label{lem-reltixi''}
	 There exists a $G$-invariant 1-cocycle 
	$\Set{\ti{\xi}'_{ij} \in \cT'({U}_{ij} \times V) }
	= \Set{\ti{\xi}'_{ij} : V \to \pi^\circ_*T_{X^\circ/ B^\circ}(U_{ij}) }$
	lifting the 1-cocycle $\{\xi_{ij}\}$ 
	and satisfying ${\ti{\xi}}'_{ij}(0) = 0$. 
\end{lem}

\begin{proof}
	Let $\{\ti{\xi}_{ij}\}$ be a 1-cochain  lifting 
	the 1-cocycle $\{\xi_{ij}\}$  constructed in Claim~\ref{claim-tixirel}
	and let $\ti{\xi}_{ijk} \cnec \ti{\xi}_{ij} +\ti{\xi}_{jk}+\ti{\xi}_{ki}$. 
	Then $\{\ti{\xi}_{ijk}\}$ is a 2-cocycle with coefficients in $\cT$ by~\eqref{suitex-extautregleeJ}.
	Define $\bar{\cT} \cnec \pr_1'^*\pi_*T_{X / S}$ 
	where $\pr_1' : B \times V \to B$ is the projection to the first factor. 
	We have $\bar{\cT}_{|B^\circ \times V} = \cT$,
	so by Lemma~\ref{lem-bonrecG}.\emph{ii)}, we can also regard $\{\ti{\xi}_{ijk}\}$ as a 
	2-cocycle  with coefficients in $\bar{\cT}$
	(with respect to the open cover $\{U_i \times V\}_{i \in I}$).
	Since $\dim B = 1$ and $V$ is Stein, investigating the
	Leray spectral sequence associated to the projection $B \times V \to V$
	shows that $H^2(B \times V,\bar{\cT}) = 0$.
	Thus, as $\{U_i \times V\}_{i \in I}$ is a good Stein cover of $B \times V$, 
it follows that $\{ \ti{\xi}_{ijk}\}$ is a 2-coboundary 
with coefficients in $\bar{\cT}$, 
and therefore with coefficients in ${\cT}$ 
(again by Lemma~\ref{lem-bonrecG}.\emph{ii)}
together with $\bar{\cT}_{|B^\circ \times V} = \cT$). 

Consider the commutative diagram
\begin{equation}
\begin{tikzcd}[cramped, row sep = 15, column sep = 20]
0 \ar[r] &  Z^1\(\{U_i \times V\}_{i\in I}, \cT\)   \ar[d, "\jmath_1"] \ar[r] & C^1\(\{U_i \times V\}_{i\in I}, \cT\)   \arrow[r] \ar[d, "\jmath_2"] & B^2\(\{U_i \times V\}_{i\in I}, \cT\)  \ar[r] \ar[d,"\jmath_3"]  & 0  \\
0 \ar[r] &  Z^1\(\{U_i\}_{i\in I}, \cT_{| B^\circ \times 0}\)  \ar[r] &  C^1\(\{U_i\}_{i\in I}, \cT_{| B^\circ \times 0}\)  \arrow[r]  &  B^2\(\{U_i\}_{i\in I}, \cT_{| B^\circ \times 0}\)  \ar[r] & 0 
\end{tikzcd}
\end{equation}
with exact rows,  where $Z^i$, $C^i$, $B^i$ are the spaces of
\v{C}ech $i$-cocycles, $i$-cochains, $i$-coboundaries. 
Since $\jmath_1$ is surjective, the induced map $\ker(\jmath_2) \to \ker(\jmath_3)$ is also surjective by the snake lemma. 
Therefore as $\{ \ti{\xi}_{ijk}\} \in B^2\(\{U_i \times V\}_{i\in I}, \cT\)$ 
satisfies $\ti{\xi}_{ijk}(0) = 0$ by construction,
there exists $\{\gb_{ij}\} \in C^1\(\{U_i \times V\}_{i\in I}, \cT\) $ such that 
$\gb_{ij} + \gb_{jk}+\gb_{ki} = \ti{\xi}_{ijk} $ and $\gb_{ij}(0) = 0$. 
Thus $\{\ti{\xi}'_{ij} \cnec \ti{\xi}_{ij} - {\gb}_{ij} \}$ is a 1-cocycle
with coefficients in $\cT'$ satisfying $\ti{\xi}'_{ij}(0) = 0$. 

Since $\{\ti{\xi}'_{ij}\} = \{\ti{\xi}_{ij}\} - \{ {\gb}_{ij} \}$, 
and the images of $\{\ti{\xi}_{ij}\}$ and $\{ {\gb}_{ij} \}$ in $C^1(\{U_i \times V\}_{i\in I}, \cE)$ 
induced by $\cT' \to \cE$ in~\eqref{suitex-extautregleeJ}
are $\{\xi_{ij}\}$ and $0$ respectively 
(because $\{ {\gb}_{ij} \} \in C^1(\{U_i \times V\}_{i\in I}, \cT)$), 
the 1-cocycle $\{\ti{\xi}'_{ij}\}$ is a   
lifting of $\{\ti{\xi}_{ij}\}$. 
As $\{{\xi}_{ij}\}$ is $G$-invariant,
up to replacing the 1-cocycle $\{{\ti{\xi}}'_{ij}\}$ by 
$\frac{1}{|G|}\sum_{g \in G} g \cdot \{{\ti{\xi}}'_{ij}\}$, we can assume that the lifting $\{{\ti{\xi}}'_{ij}\}$ of $\{{\xi}_{ij}\}$ is $G$-invariant.
\end{proof}

Recall that the exponential map $\exp : \cE \to \cJ$ 
in~\eqref{suitex-extautregleeJ} is 
isomorphic to 
$\exp_{\pr_1^{-1}\bH/B^\circ \times V} : 
\cE_{\pr_1^{-1}\bH/B^\circ \times V} \to \cJ.$
Thus if $\Set{\ti{\xi}'_{ij} \in \cT'({U}_{ij} \times V) }$ 
is a 1-cocycle as in Lemma~\ref{lem-reltixi''} and
$$\te_{ij}  \cnec \exp_2(U_{ij} \times V)(\ti{\xi}'_{ij}) : {X}_{ij} \times V \eto {X}_{ij} \times V$$ 
(where $\exp_2$ is introduced in~\eqref{suitex-extautregleeJ}), then 
$\Set{\ti{e}_{ij} \in \Aut({X}_{ij} \times V) }$
is a 1-cocycle of maps and 
is a $G$-invariant lifting of $\Set{e_{ij} \in \Aut({S}_{ij} \times V) }$ 
 by~\eqref{suitex-extautregleeJ} and the definition of $\{e_{ij}\}$.
So we can glue the collection of fibrations 
$\Set{ f_i \times \Id_V : {X}_i \times V \to {S}_i \times V }_{i \in I}$ 
along $\Set{f_{ij} \times \Id_V : {X}_{ij} \times V \to {S}_{ij} \times V }$ 
using $\{\ti{e}_{ij} \}$ and $\{{e}_{ij} \}$ 
 and obtain a $G$-equivariant map ${\cX} \to {\cS}$. Let
$$\Pi_X : {\cX} \to {\cS} \to {B} \times V \to V$$ 
be the  composition of $\cX \to \cS$ with $\Pi_S$.
Recall that $e_{ij| {S}_{ij} \times \{0\}} = \Id_{ {S}_{ij} \times \{0\}}$ by~\eqref{eqn-ideij0}. 
Since $\ti{\xi}'_{ij}(0) = 0$  by construction,  
we  have  ${\te}_{ij|{X}_{ij} \times \{0\}} = \Id_{{X}_{ij} \times \{0\}}$, 
so the central fiber of $\Pi_X$ is $\pi : X \to B$. 
That $\Pi_X$ is a deformation of ${\pi}$ which is
$G$-equivariantly locally trivial over ${B}$ follows from the construction
and the property that $\Pi_S$ is $G$-equivariantly locally trivial over ${B}$.
This terminates the proof of Proposition~\ref{pro-ellipsuraa}.
\end{proof}

\subsection{Algebraic approximations of $\bP^1$-fibrations over an elliptic surface}
\hfill

\begin{cor}\label{cor-ellipsuraa}
	Let $f : X \to S$ be a $\bP^1$-fibration over 
	an  
	elliptic surface $p:S \to B$
	satisfying the description in 
	Proposition~\ref{pro-classuniregl}.(ii). 
	There exists a bimeromorphic modification
	$X' \xto{f'} S' \xto{p'} B$ of $X \xto{f} S \xto{p} B$
	 which satisfies the following properties:
	\begin{enumerate}[label = \roman{enumi})]
		\item $X'$ is normal and the underlying bimeromorphic map $X' \dto X$
		is an isomorphism over a nonempty Zariski open $U \subset B$.
		\item There exists an algebraic approximation 
		$$\Pi_{S'} : \cS' \to B\times V \to V$$ 
		of $p$ which can be lifted to an algebraic approximation 
		$$\Pi_{X'} : \cX' \to \cS' \to B \times V \to V$$
		of $f'$ such that the underlying deformation of 
		$\pi' \colonec p' \circ f': X' \to B$ is locally trivial over $B$. 
		In particular, 
		$\Pi_{X'}$ is an $f'^{-1}(C)$-locally trivial algebraic approximation of $X'$ 
		for every subvariety $C \subsetneq S'$.
	\end{enumerate}
\end{cor}

\begin{proof}

By~\cite[Proposition 3.11]{ClaudonHorpi1}, there exists a Galois cover $r : \tilde{B} \to B$ of $B$ such that the elliptic fibration $S \times_B \tilde{B} \to \tilde{B}$ has local sections at every point of $\ti{B}$.	
Let $G \colonec \Gal(\tilde{B}/B)$. 
Let $\nu_S: \ti{S} \to  S \times_B \tilde{B}$ 
and $\nu_X: \ti{X} \to  X \times_S \tilde{S}$ be 
$G$-equivariant K\"ahler  desingularizations   
of $S \times_B \tilde{B}$ and  $X \times_S \tilde{S}$ respectively,
such that the fundamental loci of $\nu_S^{-1}$ and $\nu_X^{-1}$ 
are contained in the singular loci of $S \times_B \tilde{B}$ and $X \times_S \tilde{S}$ respectively
(see \eg~\cite[Theorem 2.5]{HYLbimkod1}).  
Let $\ti{p} : \ti{S} \xto{\nu_S} S \times_B \ti{B} \to \ti{B}$
(resp. $\ti{f} : \ti{X} \xto{\nu_X} X \times_S \ti{S} \to \ti{S}$) 
denote the induced elliptic surface (resp. $\bP^1$-fibration). 
Let $U \subset B$ be a nonempty Zariski open subset
such that both $p : S \to B$ and $\pi = (p \circ f) : X \to B$ are smooth over $U$.

\begin{claim}\label{claim-rescb}
	
	The restriction of
	$\ti{\pi} : \ti{X} \xto{\ti{f}} \ti{S} \xto{\ti{p}} \ti{B}$ over 
	$\ti{U} \cnec r^{-1}(U) \subset \ti{B}$
	is the base change of $X \xto{f} S \xto{p} B$ by $\ti{U} \xto{r} B$.
	In particular, the maps $X' \xto{f'} S' \xto{p'} B$ defined by
	$\ti{X}/G \to \ti{S}/G \to B$ satisfy i) in Corollary~\ref{cor-ellipsuraa}.
\end{claim}

\begin{proof} 
	By construction, $S \times_B \tilde{B} \to \ti{B}$ is smooth over $\ti{U}$.
	It follows from the property of the desingularization 
	$\nu_S : \ti{S} \to S \times_B \tilde{B}$ that
	necessarily, the restriction of $\ti{S} \to \ti{B}$ over $\ti{U}$
	is isomorphic to the base change $S \times_B \ti{U} \to \ti{U}$ of $p : S \to B$ 
	by $\ti{U} \xto{r} B$.
	It also follows that
	$X \times_S \tilde{S} \to \ti{B}$ is isomorphic to
	$X \times_B \tilde{B} \to \ti{B}$ over $\ti{U} \subset \ti{B}$.
	Likewise, since $X \times_B \tilde{B} \to \ti{B}$ is smooth over $\ti{U}$,
	it follows from the property of 
	the desingularization $\nu_X : \ti{X} \to X \times_S \tilde{S}$
	that the restriction of $\ti{X} \to \ti{B}$ over $\ti{U}$
	is isomorphic to the base change $X \times_B \ti{U} \to \ti{U}$ of $\pi : X \to B$ 
	by $\ti{U} \xto{r} B$. This proves the main statement of Claim~\ref{claim-rescb}.
	Finally, since $X'$ is the quotient of a complex manifold by a finite group,
	$X'$ is normal. 
\end{proof}

By Claim~\ref{claim-rescb}, the $\bP^1$-fibration $\ti{f}$
still satisfies the description in 
Proposition~\ref{pro-classuniregl}.\emph{(ii)}. 
As $\ti{p} : \ti{S} \to \ti{B}$ has 
local sections at every point of $\ti{B}$,
the $G$-equivariant $\bP^1$-fibration $\ti{f}$ and the elliptic surface $\ti{p}$ satisfy the conditions in Proposition~\ref{pro-ellipsuraa}. Therefore by Proposition~\ref{pro-ellipsuraa},
the $G$-equivariant tautological family
\begin{equation}\label{defo-surfellp}
\Pi_{\ti{S}} : \tilde{\cS}  \to  \tilde{B} \times V \to V 
\end{equation}
associated to $\ti{p}$ 
can be lifted to a deformation 
$$\Pi_{\ti{X}} : \ti{\cX} \to \ti{\cS} \to \ti{B} \times V \to V$$
of $\ti{f}$ preserving the $G$-action  
such that the underlying deformation of $\ti{\pi} \colonec \ti{p} \circ \ti{f} : \ti{X} \to \ti{B}$ is $G$-equivariantly locally trivial over $\ti{B}$.
By Lemma~\ref{lem-Gquotloctriv}, the quotient 
$$\Pi_{X'} : \cX' \colonec \ti{\cX}/G \to {\cS}' \colonec \tilde{\cS}/G \to  B \times V \to V$$
of $\Pi_{\ti{X}}$ by $G$ is a deformation of $\pi' : X' \to B$ which is locally trivial over $B$. 
As the fiber $\cX'_t$ of $\Pi_{X'}$ over $t$ is a $\bP^1$-fibration over $\cS'_t$, by Corollary~\ref{cor-critprojP1} $\cX'_t$ is algebraic if $\cS'_t$ is algebraic. 
Since $\Pi_{\ti{S}} : \ti{\cS} \to V$ is an algebraic approximation of $\ti{S}$
by Proposition~\ref{pro-def-existfam}, 
it follows that $\cX' \to  V$ is an algebraic approximation of $X'$. 
For the last statement, as  $S'$ is bimeromorphic to $S$ which is non-algebraic, every subvariety $C \subsetneq S'$ is contained in a finite union of fibers of $p' : S' \to B$ by Corollary~\ref{cor-multsecMoibase}. 
Therefore since $\Pi_{X'}$ is locally trivial over $B$, 
it is $f'^{-1}(C)$-locally trivial for every $C \subsetneq S'$. 
\end{proof}

\subsection{Proof of Propositions~\ref{pro-main3unir} and~\ref{pro-sec}}
\hfill

Now Propositions~\ref{pro-main3unir} and~\ref{pro-sec} 
follow easily from 
the results we proved in \S\ref{sec-a0K3} and \S\ref{sec-a1}.

\begin{proof}[Proof of Propositions~\ref{pro-main3unir} and~\ref{pro-sec}]
	First we prove Proposition~\ref{pro-sec}.
Let $f : X \to S$ be a $\bP^1$-fibration as in Proposition~\ref{pro-classuniregl}. According to whether $f$ is in the first or the second case of Proposition~\ref{pro-classuniregl}, 
we use either Proposition~\ref{pro-aaa=0} or Corollary~\ref{cor-ellipsuraa} to conclude.

Proposition~\ref{pro-main3unir} 
follows from Proposition~\ref{pro-classuniregl} proven in \S\ref{sec-bim}
and Proposition~\ref{pro-sec}.
\end{proof}

\section{Conclusion}\label{sec-conclf}

\ssec{Proof of Theorem~\ref{thm-mainC}}
\hfill

Assembling the results proven previously in this article, 
we finally conclude the proof of Theorem~\ref{thm-mainC}.

\begin{proof}[Proof of Theorem~\ref{thm-mainC}]
	Let $X$ be a compact threefold in the Fujiki class $\cC$ with at worst rational singularities and 
	let $a(X) \in  \{0,1,2,3\}$ be its algebraic dimension. 
If $a(X) = 3$, then $X$ is Moishezon by definition. 
If $a(X) = 2$, then Theorem~\ref{thm-mainC}  is covered by Theorem~\ref{thm-HYLkodfibellip}. 
Finally if $a(X) \le 1$, we apply Proposition~\ref{pro-main3gk1} or~\ref{pro-main3unir} to $X$, 
and conclude by Proposition~\ref{pro-red}  that $X$ has an algebraic approximation.
\end{proof}

\ssec{Some open problems}
\hfill

We now know that compact K\"ahler manifolds $X$ with either $\dim X = 3$ or $a(X) \ge \dim X - 1$ have algebraic approximations (Theorem~\ref{thm-main3} and~\cite[Theorem 1.1]{HYLkodfibellip}). As current known examples of compact K\"ahler manifolds $X$ answering negatively the Kodaira problem all satisfy $a(X) \le \dim X - 4$, the following question remains open.

\begin{que}
Does there exist a compact K\"ahler manifold $X$ of algebraic dimension $a(X) \ge \dim X - 3$ which does not have any algebraic approximation?
\end{que}

Due to Voisin's examples~\cite{Voisincs}, uniruled compact K\"ahler manifolds can fail to admit algebraic approximations starting from dimension 5 and on. Indeed, let $X$ be one of Voisin's examples of dimension 4 and consider the product $X \times \bP^1$. Since every deformation of $X \times \bP^1$ induces a deformation of the projection $X \times \bP^1 \to X$~\cite[Theorem 2.1]{RanStabMap}, $X \times \bP^1$ does not have any algebraic approximation. For uniruled fourfolds, the Kodaira problem is still open.

\begin{que}
Does there exist a uniruled compact K\"ahler fourfold $X$ which does not have any algebraic approximation, or even the homotopy type of a projective manifold?
\end{que}

\section*{Acknowledgment}

This article was written under the support of the SFB/TR 45 
"Periods, Moduli Spaces and Arithmetic of Algebraic Varieties" of the DFG (German Research Foundation) at the University of Bonn, Taiwan Ministry of Education Yushan Young Scholar Fellowship (NTU-110VV006),
and National Science and Technology Council (110-2628-M-002-006-). 
The author would like to thank 
D.~Huybrechts for helpful discussions on twisted sheaves,
J.-P.~Demailly, F.~Gounelas, A.~H\"oring, K.~Oguiso, Y.~Prokhorov, S.~Schreieder,
C. Shramov, A.~Soldatenkov, 
and C.~Voisin for general discussions and suggestions related to this work,
as well as the referees for 
the constructive comments and questions.

\bibliographystyle{plain}
\bibliography{aa3s}

\end{document}